\documentclass[a4paper,11pt]{article}
\usepackage[english]{babel}
\usepackage{amsfonts,amsthm}
\usepackage{amssymb,amsmath}
\usepackage{dsfont}
\usepackage[dvips]{graphicx}

\input{macros.sty}

\newcommand{\superdetail}[1]{}
\newcommand{\detail}[1]{#1}
\newcommand{\nodetail}[1]{}
\renewcommand{\detailtag}{\notag}
\newtheorem{postita}{Post-it}

\begin{document}

\title{Data-driven calibration of linear estimators\\
 with minimal penalties}

\author{
Sylvain Arlot
\thanks{http://www.di.ens.fr/$\sim$arlot/} \\
CNRS ; Sierra project-team\\
Laboratoire d'Informatique de \\
l'Ecole Normale Superieure \\
(CNRS/ENS/INRIA UMR 8548)\\
23, avenue d'Italie, CS 81321 \\
75214 Paris Cedex 13, France\\
\texttt{sylvain.arlot@ens.fr} \\
\and
Francis Bach
\thanks{http://www.di.ens.fr/$\sim$fbach/} \\
INRIA ;  Sierra project-team\\
Laboratoire d'Informatique de \\
l'Ecole Normale Superieure \\
(CNRS/ENS/INRIA UMR 8548)\\
23, avenue d'Italie, CS 81321 \\
75214 Paris Cedex 13, France\\
\texttt{francis.bach@ens.fr} \\
}

\date{\today}

\maketitle

\begin{abstract}
This paper tackles the problem of selecting among several linear estimators in non-parametric regression; this includes model selection for linear regression, the choice of a regularization parameter in kernel ridge regression, spline smoothing or locally weighted regression, and the choice of a kernel in multiple kernel learning.
We propose a new algorithm which first estimates consistently the variance of the noise, based upon the concept of minimal penalty, which was previously introduced in the context of model selection.
Then, plugging our variance estimate in Mallows' $C_L$ penalty is proved to lead to an algorithm satisfying an oracle inequality.
Simulation experiments with kernel ridge regression and multiple kernel learning show that the proposed algorithm often improves significantly existing calibration procedures such as generalized cross-validation.
\end{abstract}

\section{Introduction} \label{sec.intro}
Smoothing splines or kernel-based methods are now well-established tools for supervised learning, allowing to perform various tasks, such as regression or binary classification, with linear and non-linear predictors~\cite{Shawe-Taylor04,scholkopf-smola-book}. A central issue common to all regularization frameworks is the choice of the regularization parameter: while most practitioners use cross-validation procedures to select such a parameter,  data-driven procedures not based on cross-validation are rarely used. The choice of the kernel, a seemingly unrelated issue, is also important for good predictive performance: several techniques exist, either based on cross-validation, Gaussian processes or multiple kernel learning~\cite{chapelle,GP,grouplasso}.

In this paper, we consider least-squares regression and cast these two problems as the problem of selecting among several \emph{linear estimators}, where the goal is to choose an estimator with a quadratic risk which is as small as possible.
As shown in  Section~\ref{sec.linear}, this problem includes for instance model selection for linear regression, the choice of a regularization parameter in kernel ridge regression, spline smoothing, or locally weighted regression, the choice of a kernel in multiple kernel learning, the choice of $k$ in $k$-nearest-neighbors regression, and the choice of a bandwidth of Nadaraya-Watson estimators. 

Another motivation for studying linear estimators is their good theoretical properties. 
For instance, when the signal belongs to a Sobolev ball, it is known the Pinsker estimator (which is linear) is asymptotically minimax up to the optimal constant, while the best projection estimator is only rate-minimax \cite{Efr_Pin:1996,Tsy:2009}. 
Furthermore, the set of signals that are well estimated by linear estimators is very rich: it contains, for instance, sampled smooth functions, sampled modulated smooth functions and sampled harmonic functions \cite{Jud_Nem:2009}.
Finally, convergence rates of some linear spectral estimators have recently been proved optimal \cite{Cap_DeV:2007,Bau_Per_Ros:2007}.

The main contribution of the paper is to extend the notion of \emph{minimal penalty} \cite{Bir_Mas:2006,Arl_Mas:2009:pente} presented in Section~\ref{sec.linear} to all discrete classes of linear operators, and to use it for defining a fully data-driven selection algorithm satisfying a non-asymptotic oracle inequality. 
Our new theoretical results presented in
Section~\ref{sec.algo} extend similar results which were limited to unregularized least-squares regression (i.e., projection operators). We also tackle continuously parameterized families of linear estimators which are typical in ridge regression and spline smoothing (where the one-dimensional parameter to be estimated is the regularization parameter). In order to do, we derive novel concentration inequalities which may be useful in other contexts (Section~\ref{A.sec.conc}).
Our results also enlighten the classical elbow heuristics based algorithms---e.g., ``L-curve maximum curvature criterion'' for Tikhonov \cite{Rez_Hos:2009} and several others regularization problems \cite{Han_OLe:1993,Gro_Wol:2009})---by providing theoretical grounds to another L-curve based calibration algorithm. 
Finally, in Section~\ref{sec.simu}, we show that our algorithm improves the performances of classical selection procedures, such as GCV \cite{Cra_Wah:1979}, for kernel ridge regression, nearest-neighbor regression or locally weighted regression.

\section{Linear estimators} \label{sec.linear}
In this section, we define the problem we aim to solve and give several examples of linear estimators.

\subsection{Framework and notation} \label{sec.heur.fram}
Let us assume that one observes
\[ Y_i = f(x_i) + \varepsilon_i \in \R \qquad \mbox{for} \quad  i = 1, \ldots , n \enspace ,\]
where $\varepsilon_1, \ldots, \varepsilon_n$ are i.i.d.~centered random variables with $\E[\varepsilon_i^2]=\sigma^2$ unknown, $f$ is an unknown measurable function $\X \mapsto \R$ and $x_1, \ldots, x_n \in \X$ are deterministic design points. No assumption is made on the set $\X\,$. 
The goal is to reconstruct the signal $F = (f(x_i))_{1 \leq i \leq n} \in \R^n\,$, with some estimator $\Fh \in \R^n\,$, depending only on $(x_1,Y_1), \ldots, (x_n,Y_n)\,$, and having a small quadratic risk 
$n^{-1} \snorm{\Fh - F}_2^2 \, $,
where $\forall t \in \R^n\,$, we denote by $\snorm{t}_2$  the $\ell_2$-norm of $t\,$, defined as $\snorm{t}_2^2 \egaldef \sum_{i=1}^n t_i^2\,$.

In this paper, we focus on \emph{linear estimators} $\Fh$ that can be written as a linear function of $Y = (Y_1, \ldots, Y_n) \in \R^n\,$, that is,
$\Fh = A Y \,$, 
for some (deterministic) $n\times n$ matrix $A\,$.
Here and in the rest of the paper, vectors such as $Y$ or $F$ are assumed to be column-vectors.
We present in Section~\ref{sec.heur.ex} several important families of estimators   of this  form. The matrix $A$ may depend on $x_1, \ldots, x_n$ (which are known and deterministic), but not on $Y$, and may be parameterized by certain quantities---usually regularization parameter or kernel combination weights.

Let us also define, for any matrix $A \in \M_n(\R)\,$,  the largest singular value of $A$:
\[ \normmat{A} \egaldef \sup_{t \in \R^n, \, t \neq 0} \set{ \frac{ \norm{A t}_2} { \norm{t}_2 }} \enspace . \]

\subsection{Examples of linear estimators} \label{sec.heur.ex}
In this paper, our theoretical results apply to  matrices $A$ such that 
\begin{equation}
\label{eq.hyp.Al.0} 
A \in \M_n(\R) \quad \normmat{A} \leq \majnormAl \quad \tr(A^{\top} A) \leq (2-\cAl) \tr(A) \mbox{ with } \cAl \in (0,2) \enspace ,
\end{equation}
for some constants $\majnormAl$ and $\cAl$.
The main examples we have in mind are the following.

\paragraph{Ordinary least-squares regression / model selection.}   If we consider linear predictors (here, linear in the inputs $x_1,\dots,x_n$) from a design matrix $X \in \R^{n \times p}\,$, then $\Fh = AY$ with $A = X(X^\top X)^{-1} X^\top\,$, which is a projection matrix (i.e., $A^\top A= A$); 
$\Fh = A Y$ is often called a {\em projection estimator}. 
In the variable selection setting, one wants to select a subset $J \subset \{1,\dots,p\}\,$, and matrices $A$ are parameterized by $J\,$. If we denote $X_J$ the matrix of size $n \times |J|$ composed of the columns of $X$ indexed by $J$, then the matrix $A_J$ is equal to
$X_J(X_J^\top X_J)^{-1} X_J^\top$. For this matrix, we have $\tr A_J^\top A_J  = \tr A_J^2 = \tr A_J$.

\paragraph{Kernel ridge regression / spline smoothing.}   %
We assume that a positive definite kernel $k:\X \times \X \to \R$ is given, and we are looking for a function $f: \X \to \R$ in the associated reproducing kernel Hilbert space (RKHS) $\mathcal{F}\,$, with norm $\| \cdot \|_{\mathcal{F}}\,$. If $K$ denotes the $n\times n$ kernel matrix, defined by $K_{ab} = k(x_a,x_b)\,$, then the ridge regression estimator---a.k.a.~spline smoothing estimator for spline kernels~\cite{wahba}, or Tikhonov regularization 
\cite{Tik_Mor:1981}---is obtained by minimizing with respect to $f \in \mathcal{F}$~\cite{scholkopf-smola-book}:
\[ 
\frac{1}{n} \sum_{i=1}^n ( Y_i - f(x_i) )^2 + \lambda \| f\|_{\mathcal{F}}^2 \enspace .
\]
The unique solution is equal to $\widehat{f} = \sum_{i=1}^n \alpha_i k(\cdot,x_i)\,$, where $\alpha = ( K + n \lambda I)^{-1} Y \,$. This leads to the smoothing matrix $A_\lambda = K  ( K + n \lambda \Id_n)^{-1}\,$, parameterized by the regularization parameter $\lambda \in \R_+\,$. In this case, $A$ is symmetric positive semi-definite, and we have $\tr A^2 \leqslant \tr A$.

\paragraph{Multiple kernel learning / Group Lasso / Lasso.}   
We now assume that we have $p$ different kernels~$k_j\,$, feature spaces $\mathcal{F}_j$ and feature maps $\Phi_j : \mathcal{X} \to \mathcal{F}_j\,$, $j = 1,\dots,p\,$.
The group Lasso~\cite{grouped} and multiple kernel learning~\cite{Lan_etal:2004,grouplasso} frameworks consider the following objective function
\[
J(f_1,\dots,f_p) \!= \! \textstyle \frac{1}{n} \displaystyle
\sum_{i=1}^n  \big( y_i -\textstyle\sum_{j=1}^p \langle f_j, \Phi_j(x_i)\rangle \displaystyle \big)^2 \!+  {2\lambda}   \sum_{j=1}^p \| f_j\|_{\mathcal{F}_j} \!
=L(f_1,\dots,f_p) +  {2\lambda}   \sum_{j=1}^p \| f_j\|_{\mathcal{F}_j} \enspace .
\]
Note that when $\Phi_j(x)$ is simply the $j$-th coordinate of $x \in \R^p\,$, we get back the penalization by the $\ell^1$-norm and thus the regular Lasso~\cite{lasso}.

Following~\cite{pontil,simpleMKL}, by using  $a^{1/2} = \min_{ b \geqslant 0 }\frac{1}{2} \{ \frac{a}{b} + b \}  \,$, we obtain a variational formulation of the sum of norms $ 2 \sum_{j=1}^p \| f_j\|  
 =  \min_{ \eta \in \R_+^p } \sum_{j=1}^p \textstyle \left\{\frac{ \| f_j \|^2}{\eta_j} + \eta_j \right\}\,$. 
Thus, minimizing $J(f_1,\dots,f_p)$ with respect to $(f_1,\dots,f_p)$ is equivalent to minimizing  with respect
to $\eta \in \R_+^p$ (see~\cite{grouplasso} for more details):
\[
\min_{f_1,\dots,f_p }  L(f_1,\dots,f_p) + \lambda    \sum_{j=1}^p \frac{ \| f_j \|^2}{\eta_j}
+ \lambda \sum_{j=1}^p \eta_j 
=\frac{1}{n}   y^\top  \big(  \textstyle \sum_{j=1}^p \eta_j K_j  +  n \lambda \Id_n \displaystyle \big)^{-1} \! y + \lambda   \sum_{j=1}^p \eta_j \enspace ,
\]
where $\Id_n$ is the $n \times n$ identity matrix.
Moreover, given $\eta\,$,  this leads to a smoothing matrix of the form 
\begin{equation}
\label{eq.mkl}
A_{\eta,\lambda} =  \textstyle \sparen{ \sum_{j=1}^p \eta_j K_j }
 \sparen{ \sum_{j=1}^p \eta_j K_j + n \lambda \Id_n }^{-1} \enspace ,
\end{equation}
 parameterized by the regularization parameter $\lambda \in \R_+$ and the kernel combinations in $\R_+^p$---note that it depends only on $\lambda^{-1} \eta\,$, which can be grouped in a single parameter in $\R_+^p\,$. Note that it corresponds to a specific parameterization of the kernel matrix $K$ using $\eta$, and that it can be extended to other types of parameterization.
 
 Thus, the Lasso/group lasso can be seen as particular (convex) ways of optimizing over $\eta\,$. In this paper, we propose a non-convex alternative with better statistical properties (oracle inequality in Theorem~\ref{thm.algo}). Note that in our setting, finding the solution of the problem is hard in general since the optimization is not convex. However, while the model selection problem is by nature combinatorial, our optimization problems for multiple kernels are all differentiable and are thus amenable to gradient descent procedures---which only find local optima.

\paragraph{Nearest-neighbor regression.}   
If we assume that we are given any similarity measure $d:\X \times \X \to \mathbb{R}$ and $k$ a strictly positive integer, then, from $n$ observations $x_1,\dots,x_n$, we may for each $i=1,\dots,n$, find $k$-nearest neighbors of $x_i$, i.e., find any set $J_i$ of $k$ points $x_j$, $j \in \{1,\dots,n\} \backslash \{k\}$, which are among the $k$ closest to $x_i$ according to $d$ (this definition takes into account possible ties). We can then build an $n \times n$ matrix  $A$ of nearest neighbors which is equal to $1/k$ for all pairs $(i,j)$ such that $j \in J_i$ for all $i \in \{1,\dots,n\}$, and equal to zero otherwise.

\paragraph{Nadaraya-Watson estimators \cite{Nad:1964,Wat:1964}.} 
We know assume that we are given a ``window function'' (not to be confused with a positive definite kernel) $k: \X \times \X \to \mathbb{R}_+$, from  which we build the $n \times n$ matrix $W$ of pairwise evaluations. The estimator correspond to the matrix $A$ obtained by normalizing $W$ to have unit row-sums, i.e.,
$A = W D^{-1}$, where $D = {\rm Diag}(W \unmat)$ is the diagonal matrix of row sums. In this situation we have $\majnormAl \leq \sqrt{ \max_i {D_{ii}} / \min_i D_{ii}} $. A typical example is the matrix $W$ defined as $W_{ij} = \exp( - \alpha \| x_i - x_j \|^2 )$ where $x_i$, $i=1,\dots,n$, are the observed input data points, and $\alpha$ is the smoothing parameter to be learned.

\begin{table}
\begin{center}
\begin{tabular}{|l|c|c|}
\hline
Method & $A$ & parameter \\
\hline
Ridge regression & $K ( K  + \lambda \idm)^{-1}$ & $\lambda$ \\
Kernel learning & $K ( K  +  \idm)^{-1}$ & $K$ \\
Nadaraya-Watson & $W {\rm Diag}(W \unmat)^{-1}$ & $\alpha$ where $W_{ij} = \exp( - \alpha \| x_i - x_j \|^2 )$ \\
Nearest-neighbor & $A \in \{0,1/k\}^{N \times N}$ & $k$ \\
\hline
\end{tabular}
\end{center}
\caption{Examples of linear estimators.}
\end{table}

Except for the bound on $\normmat{A}$ in the $k$-nearest neighbor and Nadaraya-Watson examples, Eq.~\eqref{eq.hyp.Al.0}  holds true with $\majnormAl=1$ and $\cAl=1$ for all the examples mentioned, as shown by the following result.
\begin{proposition} \label{pro.hyp.Al}
For any $n \geq 1\,$, $ \tr(A^{\top} A) \leq \tr(A) \leq n$ for any matrix $A \in \M_n(\R)$ among the following examples\textup{:}
\begin{enumerate}
\item[(i)] if $A$ is symmetric with $\Sp(A) \subset [0,1]\,$, for instance:
\begin{enumerate}
\item[(ia)] Ordinary least-squares regression\textup{:} $A$ is an orthogonal projection matrix.
\item[(ib)] Kernel ridge regression, Multiple kernel learning\textup{:} $\exists x \in (0,+\infty)$ and $K \in \M_n(\R)$ symmetric positive semi-definite such that $A = K (K + x \Id)^{-1}\,$.
\end{enumerate}
\item[(ii)] if $\forall i,j\,$, $A_{i,i} \geq A_{i,j} \geq 0$ and $\sum_{k=1}^n A_{i,k}=1\,$, for instance:
\begin{enumerate}
\item[(iia)] Nadaraya-Watson regression
\item[(iib)] $k$-nearest-neighbor regression, for some integer $k \in [1,n]\,$\textup{:} 
\begin{equation} \tag{{\bf kNN}} \label{eq.kNN} \left\{ \begin{aligned}
\forall 1 \leq i,j \leq n \, , \quad A_{i,j} \in \set{0 , \frac{1}{k} } \quad \mbox{with } k \in \set{1, \ldots, n} \\
\forall 1 \leq i \leq n \, , \quad 
A_{i,i} = \frac{1}{k} \quad \mbox{and} \quad \sum_{j=1}^n A_{i,j} = 1 \enspace . 
\end{aligned} 
\right. 
\end{equation}
\end{enumerate}
\end{enumerate}
In example (i), $\normmat{A} \leq 1$. In examples (ia) and (iib), $\tr(A) = \tr(A^{\top} A)\,$.
\end{proposition}
Proposition~\ref{pro.hyp.Al} is proved in Section~\ref{sec.pr.pro.hyp.Al}.

\paragraph{Other examples.} Alternative linear estimators are classical in the statistical or learning literature: 
\begin{itemize}
\item Pinsker filters \cite{Pin:1980,Efr_Pin:1996}, that is, $A_{w,\alpha}=\diag\paren{ (1 - (k^{\alpha}/w))_+ \, , \, k=1\dots n} $ for some parameters $w,\alpha>0$. This example matches case (i) in Proposition~\ref{pro.hyp.Al}. 
\item Linear spectral methods for statistical inverse problems \cite{Bis_Hoh_Mun_Ruy:2007,LoG_Ros_Odo_DeV_Ver:2008}, such as spectral cut-off (or principal components regression) and $\ell_2$-boosting.
\item Symmetrized $k$-nearest neighbors \cite{Yan:1981}.
\end{itemize}
More examples and references can be found in \cite[Chapter~5]{Was:2006} and \cite[Chapter~3]{Tsy:2009}, for instance.

\section{Linear estimator selection}  \label{sec.heur}
In this section, we first describe the statistical framework of linear estimator selection, then introduce the notion of minimal penalty. Finally, we briefly review the related work on linear estimator selection. 
\subsection{Unbiased risk estimation heuristics} \label{sec.heur.ideal-sel}
Usually, several estimators of the form $\Fh = AY$ can be used. 
The problem that we consider in this paper is then to select one of them, that is, to choose a matrix $A\,$.
Let us assume that a family of matrices $(\Al)_{\lamm}$ is given (examples are shown in Section~\ref{sec.heur.ex}), hence a family of estimators $\sparen{\Fhl}_{\lamm}$ can be used, with $\Fhl \egaldef \Al Y\,$.
The goal is to choose \emph{from data} some $\lh \in \Lambda\,$, so that the quadratic risk of $\Fh_{\lh}$ is as small as possible. 

The best choice would be the \emph{oracle}:
\[ \lo \in \arg\min_{\lamm} \set{n^{-1} \snorm{\Fhl - F}_2^2 } \enspace , \]
which cannot be used since it depends on the unknown signal $F\,$. 
Therefore, the goal is to define a data-driven $\lh$ satisfying an {\em oracle inequality} of the form
\begin{equation} \label{eq.oracle}
n^{-1} \snorm{\Fh_{\lh} - F}_2^2 \leq C_n \inf_{\lamm} \set{ n^{-1} \snorm{\Fhl - F}_2^2 } + R_n \enspace ,
\end{equation}
with large probability, where the leading constant $C_n$ should be close to~1 (at least for large $n$) and the remainder term $R_n$ should be negligible compared to the risk of the oracle.
Many classical selection methods are built upon the ``unbiased risk estimation'' heuristics: If $\lh$ minimizes a criterion $\crit(\lambda)$ such that 
\[ \forall \lamm \, , \qquad \E\croch{ \crit(\lambda) } \approx \E\croch{ n^{-1} \snorm{\Fhl -F}_2^2 } \enspace , \]
then $\lh$ satisfies an oracle inequality such as in Eq.~\eqref{eq.oracle} with large probability.
For instance, %
cross-validation \cite{All:1974,Sto:1974} and generalized cross-validation (GCV)  \cite{Cra_Wah:1979}
are built upon this heuristics.

One way of implementing this heuristics is penalization, which consists in minimizing the sum of the empirical risk and a penalty term, i.e., using a criterion of the form:
\[ \crit(\lambda) = n^{-1} \snorm{ \Fhl - Y }_2^2 + \pen(\lambda) \enspace . \]

The unbiased risk estimation heuristics, also called Mallows' heuristics, then leads to the {\em optimal (deterministic) penalty}
  \[ \penid(\lambda) \egaldef \E\croch{ n^{-1} \snorm{ \Fhl - F }_2^2}  - \E \croch{ n^{-1} \snorm{ \Fhl - Y }_2^2 } \enspace . \] 
  When $\Fhl = \Al Y\,$, we have:
  \begin{align} \label{eq.riskFhl}
  \|\Fhl - F\|_2^2 &= \norm{ (\Al - \Id_n) F}_2^2 + \norm{ \Al \varepsilon}_2^2 + 2 \prodscal{\Al \varepsilon}{ (\Al - \Id_n) F} \enspace , \\ \label{eq.riskempFhl}
 \|\Fhl - Y\|_2^2 &= \| \Fhl - F\|_2^2 + \norm{\varepsilon}_2^2 - 2 \prodscal{\varepsilon}{\Al \varepsilon} + 2 \prodscal{\varepsilon}{(\Id_n - \Al) F} \enspace ,
  \end{align} 
  where 
$\varepsilon = Y - F \in \R^n$ and $\forall t,u \in \R^n\,$, $\prodscal{t}{u} = \sum_{i=1}^n t_i u_i\,$.
  Since $\varepsilon$ is centered with covariance matrix $\sigma^2 \Id_n\,$, 
  Eq.~\eqref{eq.riskFhl} and Eq.~\eqref{eq.riskempFhl} imply that 
  \begin{equation} \label{eq.penid}
  \penid(\lambda) 
  = \frac{ 2 \sigma^2 \tr(\Al)} {n} \enspace , 
  \end{equation}
  up to the term $-\E\scroch{n^{-1} \snorm{\varepsilon}_2^2}\! = - \sigma^2\,$, which can be dropped off since it does not vary with $\lambda\,$.

Note that $\dfl = \tr(\Al)$ is called the {\em effective dimensionality} or {\em degrees of freedom}~\cite{Zha:2005}, so that the optimal penalty in Eq.~\eqref{eq.penid} is proportional to the dimensionality associated with the matrix $\Al$---for projection matrices, we get back the dimension of the subspace, which is classical in model selection.

The expression  of the optimal penalty in Eq.~\eqref{eq.penid} led to several selection procedures, in particular Mallows' $C_L$ (called $C_p$ in the case of projection estimators) \cite{Mal:1973}, where $\sigma^2$ is replaced by some estimator $\widehat{\sigma^2}\,$.
The estimator of $\sigma^2$ usually used with $C_L$ is based upon the value of the empirical risk at some~$\lambda_0$ with $\df(\lambda_0)$ large; it has the drawback of overestimating the risk, in a way which depends on~$\lambda_0$ and $F$ \cite{Efr:1986}.
GCV, which implicitly estimates $\sigma^2\,$, has the drawback of overfitting if the family $(\Al)_{\lamm}$ contains a matrix too close to $\Id_n$ \cite{Cao_Gol:2006}, so that examples have been given where GCV is not asymptotically optimal \cite{KCLi:1986}; GCV also overestimates the risk even more than $C_L$ for most $\Al$ (see (7.9) and Table~4 in \cite{Efr:1986}). 

In this paper, we define an estimator of $\sigma^2$ directly related to the selection task which does not have similar drawbacks. 
Our estimator relies on the concept of minimal penalty, introduced by Birg\'e and Massart \cite{Bir_Mas:2006} and further studied in~\cite{Arl_Mas:2009:pente}.

\subsection{Minimal and optimal penalties} \label{sec.heur.mini-opt}
We deduce from Eq.~\eqref{eq.riskFhl} the {\em bias-variance decomposition} of the risk:
  \begin{align} \label{eq.EriskFhl}
  \E\croch{ n^{-1} \|\Fhl - F\|_2^2 } &= n^{-1} \norm{ (\Al - \Id_n) F}_2^2 + \frac{ \tr (\Al^\top \Al) \sigma^2} {n} = \mbox{bias} + \mbox{variance} \enspace ,
  \end{align}
  and from Eq.~\eqref{eq.riskempFhl} the expectation of the empirical risk:
  \begin{align}
  \label{eq.EriskempFhl}
  \E\croch{ n^{-1} \|{\Fhl - Y}\|_2^2 - \norm{\varepsilon}_2^2} &= n^{-1} \norm{ (\Al - \Id_n) F}_2^2 - \frac{ \paren{ 2 \tr(\Al) - \tr (\Al^\top \Al) } \sigma^2 } {n}   \enspace . 
  \end{align}

Note that the variance term in Eq.~\eqref{eq.EriskFhl} is not proportional to the effective dimensionality $\dfl = \tr (\Al)$ but to $\tr(\Al^\top \Al)\,$. 
Although several papers argue these terms are of the same order (for instance, they are equal when $\Al$ is a projection matrix), this may not hold in general. 
If Eq.~\eqref{eq.hyp.Al.0} holds for all matrices $A\in \set{\Al}_{\lamm}\,$, we only have
\begin{equation} \label{eq.comp-df1-df2} 
0 \leq \frac{\tr (\Al^\top \Al)}{2-\cAl} \leq \tr(\Al) \leq \frac{2 \tr(\Al) - \tr (\Al^\top \Al)}{\cAl} \leq \frac{2 \tr(\Al)}{\cAl}  \enspace . \end{equation}

\begin{figure}
\begin{center}
\includegraphics[scale=.45]{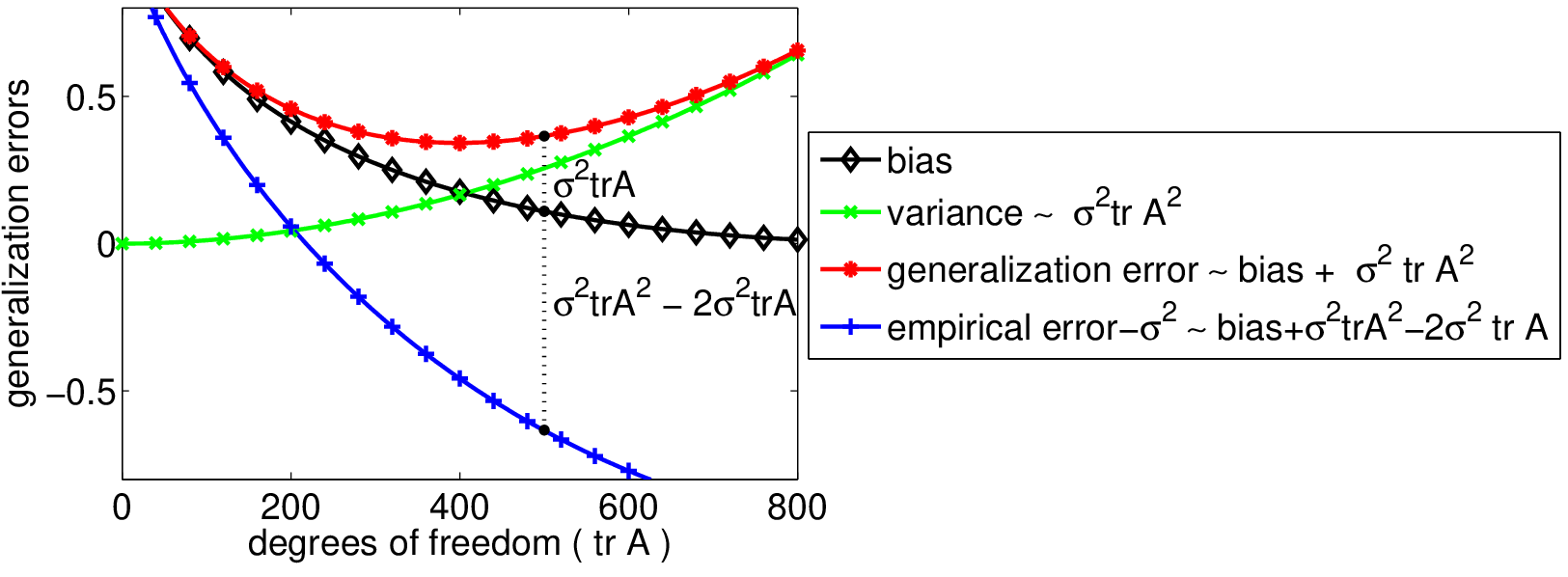}
\end{center}
\vspace*{-.5cm}

\caption{Bias-variance decomposition of the generalization error, and minimal/optimal penalties.}
\label{fig:biasvariance}
\end{figure}

In order to give a first intuitive interpretation of Eq.~\eqref{eq.EriskFhl} and Eq.~\eqref{eq.EriskempFhl}, let us consider the kernel ridge regression example, where $A = K ( K + \lambda \idm)^{-1}$,  and assume that the risk and the empirical risk behave as their expectations in Eq.~\eqref{eq.EriskFhl} and Eq.~\eqref{eq.EriskempFhl}; see also Fig.~\ref{fig:biasvariance}. 
Completely rigorous arguments based upon concentration inequalities are developed in the Appendix and summarized in Section~\ref{sec.algo}, leading to the same conclusions as the present informal reasoning.

First, as proved by Lemma~\ref{le.ridge.monotone} in Section~\ref{A.sec.conc.quad-Gauss}, the bias $n^{-1} \norm{ (\Al - \Id_n) F}_2^2$ is a non-increasing function of the dimensionality $\dfl = \tr(\Al)\,$, and the variance $ \tr (\Al^\top \Al) \sigma^2 n^{-1}$ is an increasing function of $\dfl\,$, as well as $2 \tr(\Al) - \tr (\Al^\top \Al)\,$.
Therefore, Eq.~\eqref{eq.EriskFhl} shows that the optimal~$\lambda$ realizes the best trade-off between bias (which decreases with $\dfl$) and variance (which increases with $\dfl$), which is a classical fact in model selection (see Figure~\ref{fig:biasvariance}). 

Second, the expectation of the empirical risk in Eq.~\eqref{eq.EriskempFhl} can be decomposed into the bias and a negative variance term which is the opposite of 
\begin{equation} \label{eq.penmin}
\penmin(\lambda) \egaldef n^{-1} \paren{ 2 \tr(\Al) - \tr (\Al^\top \Al) } \sigma^2 \enspace .
\end{equation}
As suggested by the notation $\penmin\,$, we will show it is a {\em minimal penalty} in the following sense. 
If 
\[ 
\forall D \geq 0, \qquad \lhmin(D) \in \arg\min_{\lamm} \set{ n^{-1} \|{\Fhl - Y}\|_2^2 + D \penmin(\lambda) } \enspace ,
\]
then, up to concentration inequalities that are detailed in Section~\ref{sec.algo.oracle}, $\lhmin(D)$ behaves like a minimizer of 
\[ g_D(\lambda) = \E\croch{ n^{-1} \snorm{\Fhl - Y}_2^2 + D \penmin(\lambda) } - \sigma^2 = n^{-1} \norm{ (\Al - \Id_n) F}_2^2 + (D-1) \penmin(\lambda) \enspace . \]
Therefore, two main cases can be distinguished:
\begin{itemize}
\item if $D<1\,$, then $g_D(\lambda)$ decreases with $\dfl$ so that $\df(\lhmin(D))$ is huge: $\lhmin(D)$ overfits.
\item if $D>1\,$, then $g_D(\lambda)$ increases with $\dfl$ when $\dfl$ is large enough, so that $\df(\lhmin(D))$ is much smaller than when $D<1\,$.
\end{itemize}
As a conclusion, $\penmin(\lambda)$ is the minimal amount of penalization needed so that a minimizer $\lh$ of a penalized criterion is not clearly overfitting.

Following an idea first proposed in \cite{Bir_Mas:2006} and further analyzed or used in several other papers such as \cite{Leb:2005,Arl_Mas:2009:pente,Mau_Mic:2008:slope}, we now propose to use that $\penmin(\lambda)$ is a minimal penalty for estimating $\sigma^2$ and plug this estimator into Eq.~\eqref{eq.penid}. Indeed, if we penalize the empirical risk $n^{-1} \snorm{\Fhl - Y}_2^2$ by $C \frac{\penmin(\lambda)}{\sigma^2}$ (which does not depend on $\sigma^2$), then the argument above suggests that around the value $D = C \sigma^{-2}=1$ (i.e., around $C=\sigma^2$), we have a jump in the selected degrees of freedom.
This leads to the algorithm described in Section~\ref{sec.algo.def}.

Note that the minimal penalty given by Eq.~\eqref{eq.penmin} is new; it generalizes previous results \cite{Bir_Mas:2006,Arl_Mas:2009:pente} where $ \penmin(\Al) = n^{-1} \tr(\Al) \sigma^2$ because all $\Al$ were assumed to be projection matrices, i.e., $\Al^\top \Al = \Al\,$.
Furthermore, our results generalize the slope heuristics $\penid \approx 2 \penmin$ (only valid for projection estimators \cite{Bir_Mas:2006,Arl_Mas:2009:pente}) to general linear estimators for which $\penid/\penmin \in (1,2]\,$.

\subsection{Related work} \label{sec.related}
Several procedures have been proposed in the literature for selecting among linear estimators. 
The most classical ones are Mallows' $C_L$ \cite{Mal:1973} and GCV \cite{Cra_Wah:1979} (which have already been introduced) and cross-validation (see \cite{Arl_Cel:2010:surveyCV} for references). 

Recently, Baraud, Giraud and Huet \cite{Bar_Gir_Hue:2010} proposed an estimator selection procedure via penalization, that applies in particular to linear estimator selection; a possible drawback of their procedure is that it strongly assumes the noise is Gaussian, since the Gaussian distribution explicitly appears in the definition of their penalty.
Two penalized maximum likelihood criteria have also been proposed for selecting the ridge regression parameter \cite{Tra:2009}, but they are only supported by simulation experiments. 

Finally, let us mention here an aggregation procedure recently proposed by Dalalyan and Salmon \cite{Dal_Sal:2011} for affine estimators. Their goal is different from ours (aggregating instead of selecting), but still related since they prove some oracle inequalities for their final estimator.

\section{Main results} \label{sec.algo}
In this section, we first describe our algorithm and then present our theoretical results. 

\subsection{Algorithm} \label{sec.algo.def}
The following algorithm first computes an estimator of $\Ch$ of $\sigma^2$ using the minimal penalty in Eq.~\eqref{eq.penmin}, then considers the optimal penalty in Eq.~\eqref{eq.penid} for selecting $\lambda\,$.
\begin{enumerate}
\item[\textbf{Input:}] $Y \in \R^n$ and $\set{\Al}_{\lL}$ a collection of matrices
\item $\forall C>0\,$, compute $\lhminb(C) \in \arg\min_{\lamm} \sset{ \snorm{\Fhl - Y}_2^2 + C  \paren{ 2 \tr(\Al) - \tr (\Al^\top \Al) } } \,$.
\item Find $\Ch$ such that $\df(\lhminb(\Ch)) \in \croch{ n/10 , n/3} \, $.
\item Select $ \lh \in \arg\min_{\lamm} \sset{ \snorm{\Fhl - Y}_2^2 + 2 \Ch  \tr(\Al)  } \, $.
\end{enumerate}
In the steps 1 and 2 of the above algorithm, in practice, a grid in log-scale is used, and our theoretical results from the next section suggest to use a step-size of order $n^{-1/2}\,$. 
Step~1 can be solved efficiently (at least when $\Lambda$ is finite or $\Lambda$ is embedded with a total order) thanks to Algorithm~2 in \cite{Arl_Mas:2009:pente}.

Note that it may not be possible in all cases to find a $C$ such that $\df(\lhminb(C)) \in \scroch{ n/10 , n/3}\,$; therefore, our condition in step 2, could be relaxed to finding a $\Ch$ such that for all $C > \Ch (1 + \delta)\,$, $\df(\lhminb(C)) < n/10$ and for all
$C < \Ch / (1 + \delta) \,$, $\df(\lhminb(C)) > n/ 10\,$, with $\delta \propto \sqrt{\ln(n)/n}\,$.

Alternatively, using the same grid in log-scale, we can select $\Ch$ with maximal jump between successive values
of $\df(\lhminb(C))$---note that our theoretical result then does not entirely hold, as we show the presence of a jump around $\sigma^2\,$, but do not show the absence of similar jumps elsewhere. See examples in Section~\ref{sec.simu}.

\subsection{Assumptions} \label{sec.algo.hyp}
Before stating our main results, let us state the main assumptions we make on $(\Al)_{\lamm}$ and on the distribution of the noise. We essentially consider a set of linear estimators which corresponds to a union of a discrete set and a union of matrices obtained from ridge regression (which are themselves parameterized by a single continuous parameter).
\begin{itemize}
\item
 Assumption on the matrices $\Al$: 
some constants $\cAl \in (0,2)$ and $\majnormAl \geq 1$ exists such that 
\begin{equation}
\tag{\ensuremath{\mathbf{H\Al}}}
\label{hyp.Al}
\left. 
\begin{split}
&\forall \lambda \in \Lambda \, , \quad 
\Al \in \mathcal{M}_n(\R) \mbox{ is deterministic,} \quad 
\normmat{\Al} \leq \majnormAl \quad \\
&\tr\paren{\Al} \leq n \quad \mbox{and} \quad 
\tr(\Al^{\top} \Al) \leq (2-\cAl) \tr(\Al) \enspace .
\end{split}
\right\}
\end{equation}
\item
Assumption on $\Lambda$: $\Id \in \set{\Al}_{\lL}$ and 
\begin{equation}
\tag{\ensuremath{\mathbf{H\Lambda}}}
\label{hyp.Lam}
\left.
\begin{aligned} 
&\Lambda \subset \Lambda_0 \cup \bigcup_{j=1}^{\NLamfusrid} \paren{ \set{j} \times [0,+\infty] }
\quad 
\mbox{with} \quad \card(\Lambda_0) \leq \CLamfusdis n^{\alLamfusdis} 
\qquad \NLamfusrid \leq \CLamfusrid n^{\alLamfusrid} \\
&\mbox{and} \quad \forall j \in \set{1, \ldots, \NLamfusrid} \, , \quad 
\forall x \in (0,+\infty) \, , \quad 
A_{(j,x)} = K_j (K_j + n x \Id_n)^{-1} \\
&\qquad \mbox{with} \quad K_j \in \M_n(\R) \backslash \set{\zeromat{n}}
\quad \mbox{symmetric positive semi-definite}
\\ &\qquad \mbox{and} \quad A_{(j,0)} = \Id \quad A_{(j,+\infty)}=\zeromat{n} 
\end{aligned}
\right\}
\end{equation}
\item
 Assumption on the noise: 
\begin{align}
\tag{\ensuremath{\mathbf{H\mathcal{N}\sigma^2}}}
\label{hyp.eps-Gauss-hom} 
\varepsilon_1 , \ldots, \varepsilon_n \mbox{ are i.i.d. } \sim \mathcal{N}(0,\sigma^2)
\end{align}
\item Assumption on the bias:
\begin{gather} 
\label{hyp.thm.mini.2} \tag{\ensuremath{\mathbf{Abias}}}
\exists \lhypbias \in \Lambda, \quad \df(\lhypbias) \leq \sqrt{n} \quad \mbox{and} \quad b(\lhypbias) \leq \sigma^2 \sqrt{n \ln(n)} \enspace . 
\end{gather}
\end{itemize}
The above assumption set is discussed in details in Section~\ref{sec.discuss-assumptions}. 
\begin{remark}
By Proposition~\ref{pro.hyp.Al}, under assumption \eqref{hyp.Lam}, assumption \eqref{hyp.Al} holds with $\cAl=1$ and $\majnormAl=1$ for all $\lambda \in \bigcup_j \set{j} \times [0,+\infty]\,$. 
Furthermore, if all matrices $(\Al)_{\lL_0}$ are among the examples of Proposition~\ref{pro.hyp.Al}, assumption \eqref{hyp.Al} holds with $\cAl=1$ and $\majnormAl=\sup_{\lL_0} \norm{\Al}\,$, which happens to be close to 1 with large probability in all the examples we considered in our simulation experiments.
\end{remark}


\subsection{Minimal penalty} \label{sec.algo.penmin}

\begin{theorem} \label{thm.mini}
Let $\lhminb$ be defined by 
\begin{equation} \label{eq.lhmin} 
\forall C \geq 0, \qquad \lhminb(C) \in \arg\min_{\lamm} \set{ \norm{\Fhl - Y}_2^2 + C \paren{ 2 \tr(\Al) - \tr(\Al^\top \Al)} } \enspace . \end{equation}
Assume \eqref{hyp.Lam}, \eqref{hyp.Al}, \eqref{hyp.eps-Gauss-hom} and \eqref{hyp.thm.mini.2} hold true.
Let 
\[ \betaThmMiniBelow = \cteA \majnormAl \quad \betaThmMiniAbove = \frac{\cteB \majnormAl}{\cAl} \qquad \alLamMix = \max\set{\alLamfusdis \, , \, \alLamfusrid + 2} \quad \CLamMix = 6 \CLamfusdis + \CLamfusrid \enspace . \]
Then, for every $\delta \geq 2\,$, a constant $\nminThmMini(\cAl,\alLamMix+\delta,\majnormAl)$ exists such that for every $n \geq \nminThmMini \,$, 
\begin{align} \label{A.eq.thm.lhmin.below}
\forall 0 \leq C < \paren{ 1 - \betaThmMiniBelow (\alLamMix+\delta) \sqrt{\frac{\ln(n)}{n}}
} \sigma^2 , \quad \df(\lhminb(C)) \geq \frac{n}{3} 
\\
 \label{A.eq.thm.lhmin.above}
 \mbox{and} \quad 
\forall C > \paren{ 1 +  \betaThmMiniAbove (\alLamMix+\delta) \sqrt{\frac{\ln(n)}{n}} } \sigma^2  \, , \quad \df(\lhminb(C)) \leq \frac{n}{10}  
\end{align}
hold with probability at least $1 - \CLamMix n^{-\delta}\,$.
\end{theorem}
Theorem~\ref{thm.mini} is proved in Section~\ref{A.sec.thm.mini}.

The first important consequence of Theorem~\ref{thm.mini} is that under mild assumptions (see Section~\ref{sec.discuss-assumptions}), with a large probability, the constant $\Ch$ obtained by the algorithm of Section~\ref{sec.algo.def} satisfies 
\[ \paren{ 1 - \betaThmMiniBelow (\alLamMix+\delta) \sqrt{\frac{\ln(n)}{n}}
} \sigma^2 \leq \Ch \leq \paren{ 1 +  \betaThmMiniAbove (\alLamMix+\delta) \sqrt{\frac{\ln(n)}{n}} } \sigma^2  \enspace . \]
In particular, $\Ch$ is a consistent estimator of $\sigma^2$ , but Theorem~\ref{thm.mini} actually provides precise non-asymptotic \emph{multiplicative} bounds for $\Ch$ (with a convergence rate of order $\sqrt{\ln(n)/n}$), that are crucial for deriving a non-asymptotic oracle inequality for the algorithm of Section~\ref{sec.algo.def}, as emphasized in Remark~\ref{rk.Thm.opt.sigma-known} below. 

Compared to classical estimators of $\sigma^2\,$, such as the one usually used with Mallows' $C_L$, $\Ch$ does not depend on the choice of some model assumed to have almost no bias, which can lead to overestimating $\sigma^2$ by an unknown amount \cite{Efr:1986}.

\begin{remark}
On the same event and under the same assumptions as Theorem~\ref{thm.mini}, 
\[ \forall C > \paren{ 1 +  \frac{\betaThmMiniAbove}{10} (\alLamMix+\delta) \frac{\sqrt{\ln(n)}}{n^{1/4}} } \sigma^2  \, , \quad \df(\lhminb(C)) \leq n^{3/4},   \]
by taking $\majdimapres = n^{3/4}$ in Proposition~\ref{A.pro.mini}. Therefore, a clearer jump can be observed by looking at $\df(\lhminb(C))$ at a slightly less precise scale, which results in a possible loss in the estimation of $\sigma^2$. 
\end{remark}
\begin{remark}
The precise values $n/3$ and $n/10$ in Eq.~\eqref{A.eq.thm.lhmin.below}--\eqref{A.eq.thm.lhmin.above} have no particular meaning\textup{:} $(n/3,n/10)$ could be replaced by $(n/\kappa,n/\kappa^{\prime})$ for any $\kappa^{\prime}>\kappa>2\,$. If all matrices $\Al$ correspond to $k$-NN or OLS estimators, we can even take any $\kappa^{\prime}>\kappa>1\,$, for instance, $(9n/10,n/10)\,$.
\end{remark}


\subsection{Oracle inequality}  \label{sec.algo.oracle}

Define 
\begin{equation} \label{def.lhopt}
\forall C \geq 0 \, , \qquad \lhopt(C) \in \arg\min_{\lamm} \set{ \norm{\Fhl - Y}_2^2 + 2 C \tr(\Al) } \enspace . \end{equation} 
We have the following general theorem. 
\begin{theorem} \label{thm.oracle.lhopt}
Let $\lhopt(C)$ be defined by Eq.~\eqref{def.lhopt} and assume that \eqref{hyp.Lam}, \eqref{hyp.Al} and \eqref{hyp.eps-Gauss-hom} hold true. 
Let 
\begin{gather*} 
\betaThmOracleD = 
\max\set{ 32 \majnormAl^2 \, , \, 24 \majnormAl^2 + \cteF } 
\qquad \mbox{and} \qquad
\betaThmOracleE = 
\cteE \enspace .
\end{gather*} 
Then, for every $\delta \geq 2\,$, if $n \geq \nminThmOracle$ some absolute constant, 
with probability at least $1 - \CLamMix n^{-\delta}\,$,
for every $C>0$ and every $\eta \in (0,2)\,$,
\begin{equation} \label{eq.thm.lhopt}
\begin{split}
 n^{-1} \norm{\Fh_{\lhopt(C)} - F}_2^2 
 &\leq 
\paren{1+\eta} \inf_{\lamm} \set{ n^{-1}  \norm{\Fh_{\lambda} - F}_2^2 + \frac{2 (C-\sigma^2)_+ \tr(\Al)}{n} } 
\\
& \quad 
+ \frac{32}{\eta} \paren{C \sigma^{-2} - 1}^2 \sigma^2 \un_{C \leq \sigma^2} 
+ \paren{\betaThmOracleD + \frac{\betaThmOracleE}{\eta}} \frac{\ln(n) (\delta+\alLamMix) \sigma^2}{n}
\end{split}
\end{equation}
and
\begin{equation} \label{eq.thm.lhopt-alt}
\begin{split}
 n^{-1} \norm{\Fh_{\lhopt(C)} - F}_2^2 
 &\leq 
\paren{1+\eta} \inf_{\lamm} \set{ n^{-1}  \norm{\Fh_{\lambda} - F}_2^2 } 
+ \frac{32}{\eta} \paren{C \sigma^{-2} - 1}^2 \sigma^2 
\\
& \quad + \paren{\betaThmOracleD + \frac{\betaThmOracleE}{\eta}} \frac{\ln(n) (\delta+\alLamMix) \sigma^2}{n}
\enspace .
\end{split}
\end{equation}
\end{theorem}
Theorem~\ref{thm.oracle.lhopt} is proved in Section~\ref{A.sec.thm.oracle}. Before applying it to the proposed algorithm, let us make a few remarks. 
\begin{remark}
Note that the two inequalities in Eq.~(\ref{eq.thm.lhopt}) and Eq.~(\ref{eq.thm.lhopt-alt}) differ from their treatment of underestimation and overestimation of $C$ (compared to $\sigma^2$). As shown \emph{explicitly} in  Eq.~(\ref{eq.thm.lhopt}), overestimation of $C$ leads to a increase of generalization cost which grows as 
$\frac{2 \tr(\Al)}{n} $, which is small, while underestimation leads to a constant factor, which is more problematic, and requires greater care when estimating $C$ from data, which we do in this paper.
\end{remark}

\begin{remark}
When $C$ is deterministic, we can deduce from Eq.~\eqref{eq.thm.lhopt}--\eqref{eq.thm.lhopt-alt} an oracle inequality in expectation; we refer to the proof of Theorem~\ref{thm.algo} for details.
\end{remark}

\begin{remark}
If $\eta = (\ln(n))^{-1}$ and $\absj{C \sigma^{-2} - 1 } = O\sparen{\sqrt{\ln(n)/n}}\,$, the remainder term in Eq.~\eqref{eq.thm.lhopt}--\eqref{eq.thm.lhopt-alt} is of order $ (\ln(n))^3 \sigma^2 n^{-1}\,$, which is negligible in front of the risk of the oracle provided that $\vdeux(\lo)$ grows with $n$ faster than $(\ln(n))^3\,$, since the risk of $\Fhlo$ is at least of order $\vdeux(\lo) n^{-1}\,$. 
This usually holds when the bias is not exactly zero for some $\lamm$ with $\tr(\Al^\top \Al)$ too small, as often assumed in the model selection literature for proving asymptotic optimality results.
\end{remark}

\begin{remark}
The term $\paren{C \sigma^{-2} - 1}^2 / \theta $ could be lowered with some assumption on $\sup_{\lL} \frac{\tr(\Al)}{\tr(\Al^{\top}\Al)}$. Eq.~\eqref{eq.thm.lhopt}--\eqref{eq.thm.lhopt-alt} actually correspond to the worst-case situation, where no assumption is made. For instance, if all estimators are least-squares or $k$-nearest neighbours estimators, this term can be replaced by $\absj{C \sigma^{-2} - 1}$.
\end{remark}

\begin{remark} \label{rk.Thm.opt.sigma-known}
When the noise-level $\sigma^2$ is known, taking $C=\sigma^2$ in Theorem~\ref{thm.oracle.lhopt} \textup{(}that is, considering Mallows' $C_L$ penalty\textup{)} yields the following oracle inequality instead of Eq.~\eqref{eq.thm.lhopt-alt}:
\begin{equation*} 
n^{-1} \norm{\Fh_{\lhopt(C)} - F}_2^2 
 \leq 
\paren{1+\eta} \inf_{\lamm} \set{ n^{-1}  \norm{\Fh_{\lambda} - F}_2^2 } 
+ \paren{\betaThmOracleD + \frac{\betaThmOracleE}{\eta}} \frac{\ln(n) (\delta+\alLamMix) \sigma^2}{n}
\enspace .
\end{equation*}
As underlined in previous works \cite{Bar_Gir_Hue:2007,Bar_Gir_Hue:2010}, dealing with the case of {\em unknown variance} is more challenging, even in the model selection case. 
In Theorem~\ref{thm.oracle.lhopt} above, the remainder term $\paren{C \sigma^{-2} - 1}^2 \sigma^2 $ underlines how important it is to have $C$ close to $\sigma^2$. 
We conjecture this term is essentially unimprovable in the case where $\tr(\Al^{\top} \Al)$ is close to its lower bound $n^{-1} \tr(\Al)^2$ for $\lambda$ ``close'' to the oracle, since the penalty $2 \sigma^2 n^{-1} \tr(\Al)$ can then be an order of magnitude higher than the risk, which is the sum of the bias and of $\sigma^2 n^{-1} \tr(\Al^{\top} \Al)$. 
Therefore, estimating $\sigma^2$ with a precision as high as the one guaranteed in Theorem~\ref{thm.mini} is crucial to derive an oracle inequality valid for all linear estimators. 
\end{remark}


\subsection{Combined result} \label{sec.algo.thm-combined}
As a corollary of Theorems~\ref{thm.mini} and~\ref{thm.oracle.lhopt}, we get the following non-asymptotic oracle inequality (with leading constant arbitrarily close to 1) for the algorithm proposed in Section~\ref{sec.algo.def}.
\begin{theorem} \label{thm.algo}
Let $\Ch$ and $\lh$ be defined as in the algorithm of Section~\ref{sec.algo.def}.
Assume that \eqref{hyp.thm.mini.2}, \eqref{hyp.Lam}, \eqref{hyp.Al} and \eqref{hyp.eps-Gauss-hom} hold true.
Then, for every $\delta \geq 2\,$, a numerical constant $\betaThmMixA$ and a constant $\nminThmMix(\cAl,\delta+\alLamMix,\majnormAl)$ exist such that if $n \geq \nminThmMix\,$, 
with probability at least $1 - \CLamMix n^{-\delta}\,$, for every $\eta \in (0,2)\,$, 
\begin{equation} \label{eq.thm.oracle}
\begin{split}
\frac{1}{n} \norm{\Fh_{\lh} - F}_2^2 \leq 
\paren{ 1 + \eta } 
\inf_{\lamm} \set{  \frac{1}{n} \norm{\Fhl - F}_2^2 }
+  \betaThmMixA \majnormAl^2 \paren{ \alLamMix + \delta }^2  \frac{\ln(n) \sigma^2}{\eta n} 
\enspace .
\end{split}
\end{equation}
As a consequence, if $n \geq \nminThmMix(\cAl,\alLamMix+2,\majnormAl)\,$, for every $\eta \in (0,2)\,$, 
\begin{equation} \label{eq.thm.oracle.E}
\begin{split}
\E\croch{ \frac{1}{n} \norm{\Fh_{\lh} - F}_2^2 }
\leq \paren{ 1 + \eta } 
\E\croch{ \inf_{\lamm} \set{  \frac{1}{n} \norm{\Fhl - F}_2^2 } } 
+  \betaThmMixA \majnormAl^2 \paren{ \alLamMix + 2 }^2  \frac{\ln(n) \sigma^2}{\eta n} 
\\
+ \frac { 2 \sqrt{\CLamMix} \paren{ 2 \sigma^2 \majnormAl^2 + \paren{1 + \majnormAl}^2 n^{-1} \norm{F}^2 } } { n } 
\enspace .
\end{split}
\end{equation}
\end{theorem}
Theorem~\ref{thm.algo} is proved in Section~\ref{A.sec.pr.thm.algo}. Its main consequences are detailed in Section~\ref{sec.thm.algo.consequences}. 

Note that taking $\eta = \frac{\betaThmMixA \majnormAl^2 \paren{ \alLamMix + 2 }^2  \ln(n) \sigma^2}{\E\croch{ \inf_{\lamm} \set{  \norm{\Fhl - F}_2^2 }}} $ in Eq.~\eqref{eq.thm.oracle.E} yields a non-asymptotic oracle inequality with leading constant one, that directly implies an asymptotic optimality result when $\frac{\ln(n) \sigma^2}{n}$ is negligible in front of the risk of the oracle. 

\subsection{Discussion of the assumptions} \label{sec.discuss-assumptions}
\paragraph{Assumption \eqref{hyp.Al}.}  
It holds true with $\cAl=1$ in all the main examples detailed in Section~\ref{sec.heur.ex}, as proved by Proposition~\ref{pro.hyp.Al}.
Note that $\normmat{\Al} \leq \majnormAl$ barely is an assumption.
With $\majnormAl=1$, it means that $\Al$ actually shrinks $Y\,$, which holds true for all examples except nearest-neighbors and Nadaraya-Watson. 
In general, since $\normmat{\Al}$ is observable, one only has to check that $\majnormAl \egaldef \sup_{\lL} \normmat{\Al}$ is not much larger than 1, as we observed in all our simulation experiments. 
Note that $\normmat{\Al} \leq \majnormAl < +\infty $ is assumed in \cite{KCLi:1987} for proving asymptotic optimality results when selecting among nearest-neighbors or Nadaraya-Watson estimators.

\paragraph{Assumption~\eqref{hyp.thm.mini.2}.}
It holds if $b(\lambda) \leq n \sigma^2 c \dfl^{-d}$ for some $c \geq 0$ and $d \geq 1$ (and some $\lL$ exists with $\df(\lambda) \leq \sqrt{n}$ close to $\sqrt{n}$), a standard assumption in the context of model selection.
Besides, \eqref{hyp.thm.mini.2} is much less restrictive and can even be relaxed, see Appendix~\ref{A.sec.thm.mini}. For instance, it is sufficient to have $b(\lambda) \leq n \sigma^2 c \dfl^{-d}$ for one among several families of estimators, without having to know which one it is.

\paragraph{Gaussian noise and \eqref{hyp.Lam}.}   
For proving Theorems~\ref{thm.mini}, \ref{thm.oracle.lhopt} and \ref{thm.algo}, a key ingredient is a uniform concentration inequality for four functions of $\lambda$ and $\varepsilon\,$, that is, a lower bound on the probability of the event $\Omega_{(\alLamMix + \delta) \ln(n)}(\Lambda)$ defined in Section~\ref{A.sec.conc.event}. 
In particular, assumptions \eqref{hyp.Lam} and \eqref{hyp.eps-Gauss-hom} are only used for proving this event has a probability at least $1-\CLamMix n^{-\delta}\,$. 
So, any alternative assumption set under which a similar uniform concentration inequality could be proved could be used instead of \eqref{hyp.Lam} and \eqref{hyp.eps-Gauss-hom}, leading to the same results (except for the precise values of the constants). We refer to the proofs for details. 

For instance, kernel ridge regression could be replaced in assumption \eqref{hyp.Lam} by Pinsker filters (see Section~\ref{sec.heur.ex}) with the one-dimensional parameter set $w \in (0,+\infty)$. 

When $\varepsilon$ is sub-Gaussian, the key concentration results (Lemma~\ref{le.Omega.finite}) can certainly be proved for $\xi = \sigma^{-1} \varepsilon$ (possibly with additional small deviation terms that do not change the core of the proof) at the price of additional technicalities (at least when $\NLamfusrid=0$, i.e., $\Lambda$ is finite), which implies that Theorems~\ref{thm.mini}, \ref{thm.oracle.lhopt} and~\ref{thm.algo} would still be valid.
Note that assuming the noise is Gaussian is classical when proving non-asymptotic oracle inequalities \cite{Cao_Gol:2006,Bar_Gir_Hue:2010}; only asymptotic results exist about linear estimators selection with moment conditions on the noise \cite{KCLi:1986,KCLi:1987}. Considering heavy-tailed noise is of clear interest but beyond the scope of this paper.

\subsection{Main consequences of Theorem~\ref{thm.algo} and comparison with previous results} \label{sec.thm.algo.consequences}
\paragraph{Oracle inequality.}   
The algorithm of Section~\ref{sec.algo.def} satisfies a non-asymptotic  oracle inequality with high probability, as shown by Eq.~\eqref{eq.thm.oracle}: 
The risk of the selected estimator $\Fh_{\lh}$ is close to the risk of the oracle, up to a remainder term which is negligible when the dimensionality $\df(\lo)$ of the oracle grows with $n$ faster than $\ln(n)\,$, a typical situation when the bias is never equal to zero, for instance in kernel ridge regression.

Eq.~\eqref{eq.thm.oracle} is {\em non-asymptotic}, meaning that it holds for every fixed $n$ as soon as the assumptions explicitly made in Theorem~\ref{thm.algo} are satisfied.
Most results (all but the more recent ones for linear estimator selection \cite{Cao_Gol:2006,Bar_Gir_Hue:2010}) are asymptotic, meaning that $n$ is implicitly assumed to be larged compared to each parameter of the problem. 
This assumption can be problematic for several learning problems, for instance in multiple kernel learning when the number $p$ of kernels  may grow with $n\,$.

Another important feature of Eq.~\eqref{eq.thm.oracle} is that it holds with high probability, which is stronger than most results only true in expectation, that is, similar to Eq.~\eqref{eq.thm.oracle.E}, or even weaker. As emphasized by \cite{KCLi:1986}, the difference is significant, since some examples exist where a procedure (GCV) is asymptotically optimal in expectation (like Eq.~\eqref{eq.thm.oracle.E}) but not a.s. (like Eq.~\eqref{eq.thm.oracle}).

Several oracle inequalities have been proved in the statistical literature for Mallows' $C_L$ with $\sigma^2$ known, either asymptotic \cite{KCLi:1986,KCLi:1987} or non-asymptotic \cite{Cao_Gol:2006}. 
When $\sigma^2$ is unknown (and replaced by a consistent estimator), up to the best of our knowledge, guarantees for $C_L$ are only available for the model selection problem (see \cite{Bir_Mas:2006,Arl_Mas:2009:pente} and references therein).

\paragraph{Comparison with other procedures.}  
Oracle inequalities or asymptotic optimality results have been proved for several other linear estimator selection procedures. 
\subparagraph{Generalized Cross Validation (GCV) \cite{Cra_Wah:1979}} was mostly studied for model selection and (kernel) ridge regression. GCV is asymptotically optimal in expectation under mild assumptions \cite{Cra_Wah:1979}, but additional restrictions are needed for its almost sure asymptotic optimality \cite{KCLi:1986}: some ridge regression example exists where GCV is not asymptotically optimal whereas $C_L$ is \cite{KCLi:1986}. 
Up to the best of our knowkedge, except for model selection, non-asymptotic oracle inequalities for GCV only exist for (kernel) ridge regression \cite{Cao_Gol:2006}; their result requires a prior upper bound $\tr(\Al) \leq n/5$ for all $\lL$, and it is only valid in expectation, so it is weaker than Eq.~\eqref{eq.thm.oracle}. Moreover, compared to \cite{Cra_Wah:1979}, our results are applicable to all linear estimators and identify key assumptions regarding the bias and variance of our collection of linear estimators (i.e., Assumption~\eqref{hyp.thm.mini.2}).

Asymptotic optimality results also exist for GCV for $k$-nearest neighbors regression \cite{KCLi:1987}, that require the same assumption $\norm{\Al}\leq \majnormAl$ than we have. 
\subparagraph{Cross-validation methods \cite{Arl_Cel:2010:surveyCV}} also satisfy some asymptotic optimality results in the nearest-neighbor regression case \cite{KCLi:1987}. Compared to GCV or to the algorithm of Section~\ref{sec.algo.def}, cross-validation methods also suffer from a large computational cost, that can only be lowered by considering $V$-fold cross-validation with $V$ rather small, at the price of an increased risk. See also the simulation study of Section~\ref{sec.simu} for a numerical comparison.

\subparagraph{Baraud, Giraud and Huet's penalization procedure \cite{Bar_Gir_Hue:2010}}
also satisfies a non-asymptotic oracle inequality under mild assumptions (with a leading constant $C>1$ and a remainder term that can be large) assuming the noise is Gaussian. 
Compared to our minimal penalty algorithm, their procedure is more general and can deal with arbitrarily large collections $\Lambda$ (putting aside computational complexity). 
Nevertheless, their algorithm is slightly more complex, and probably more dependent on the Gaussian assumption on the noise, since the Gaussian distribution explicitly appears in the definition of their penalty. 

\paragraph{No overfitting.} 
Finally, let us mention a special feature of the minimal penalty algorithm that may not appear in the comparison of theoretical results, in particular because it plays an important role only at second order and when $n$ is small. 
By construction, the algorithm of Section~\ref{sec.algo.def} selects $\lh$ with an effective dimensionality larger than $\lhminb(\Ch)$ at which the jump occurs. Therefore, our algorithm {\em never overfits too much}, in addition to the theoretical risk bounds we have proved. 
This is a quite interesting property compared for instance to GCV, which is likely to overfit if it is not corrected, because GCV minimizes a criterion proportional to the empirical risk.

\section{Simulations} \label{sec.simu}
In this section, we report simulations experiments with examples of fixed design regression to illustrate our theoretical results. For simulations on random design regression, see~\cite{Arl_Bac:2009:minikernel_nips}.

We consider the following example: we take $n=200$  and use $n$ points uniformly spaced in $[0,1]$, i.e., $x_i = (i-1)/(n-1) \in [0,1]$. We then consider $Y_i = \sin( 25 \pi X_i) + \varepsilon_i$ or $Y_i = \sin( 25 \pi X_i^3) + \varepsilon_i$, where $\varepsilon_i$ are independent standard Gaussian random variables.

We performed experiments with kernel ridge regression with kernel  $k(x,y) \!=\! \prod_{i=1}^d e^{-|x_i-y_i| }\,$ (where the goal is to learn the regularization parameter),  with nearest-neighbor regression (where the goal is to learn the number of neighbors), and with locally-weighted regression (where the goal is to learn the bandwidth of the kernel).

\begin{figure}
%
%
\begin{center}
\hspace*{-1cm}
\includegraphics[scale=.485]{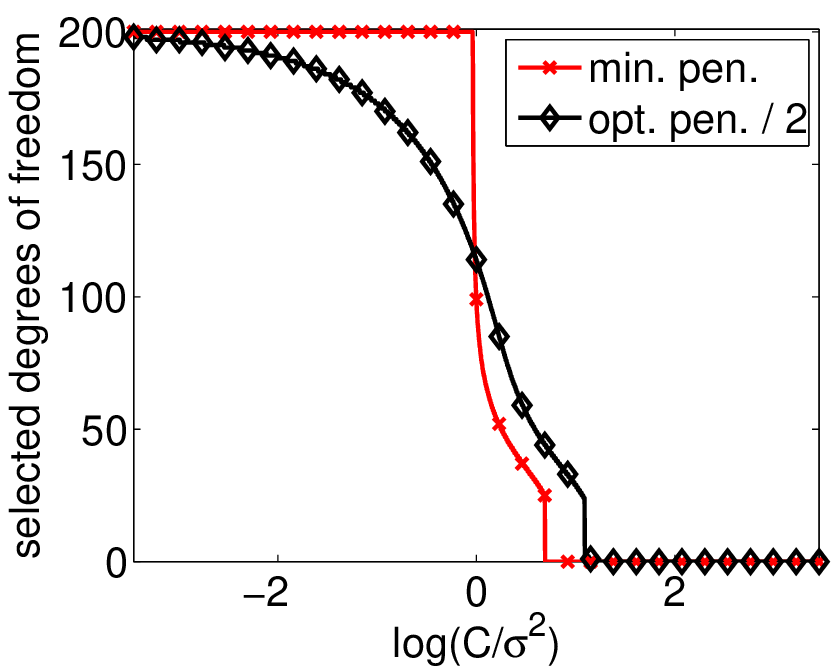} \hspace*{.36cm}
\includegraphics[scale=.485]{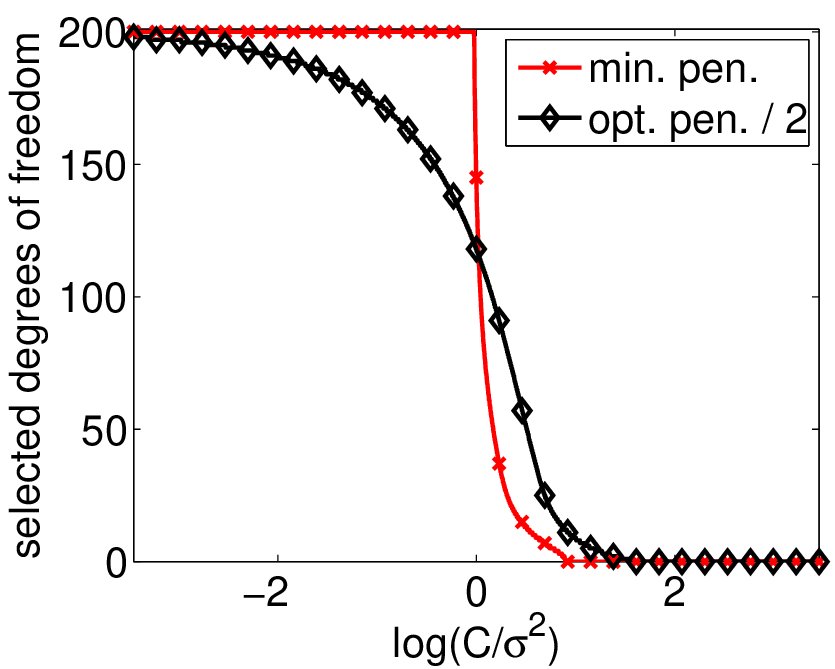} \hspace*{-1cm}

\hspace*{-1cm}\includegraphics[scale=.485]{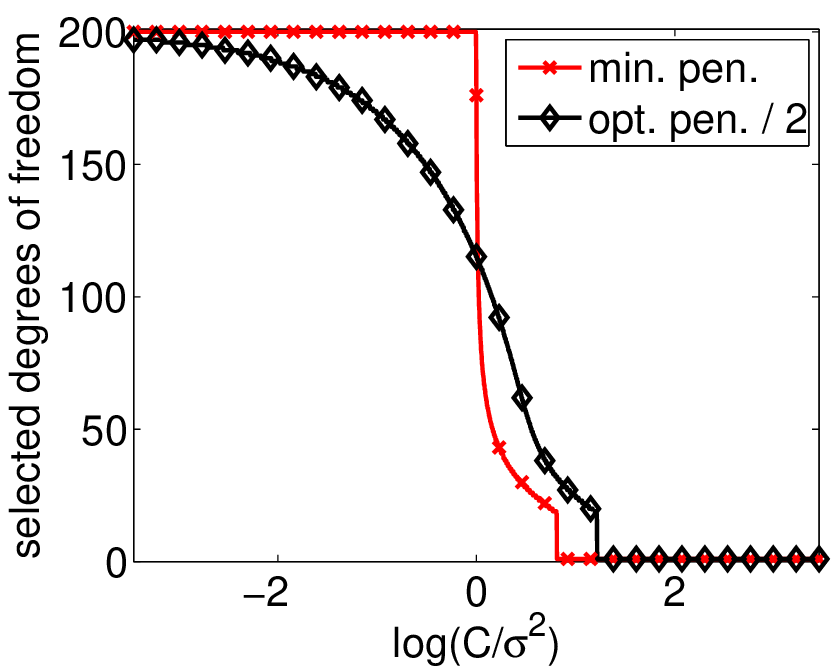} \hspace*{.36cm}
\includegraphics[scale=.485]{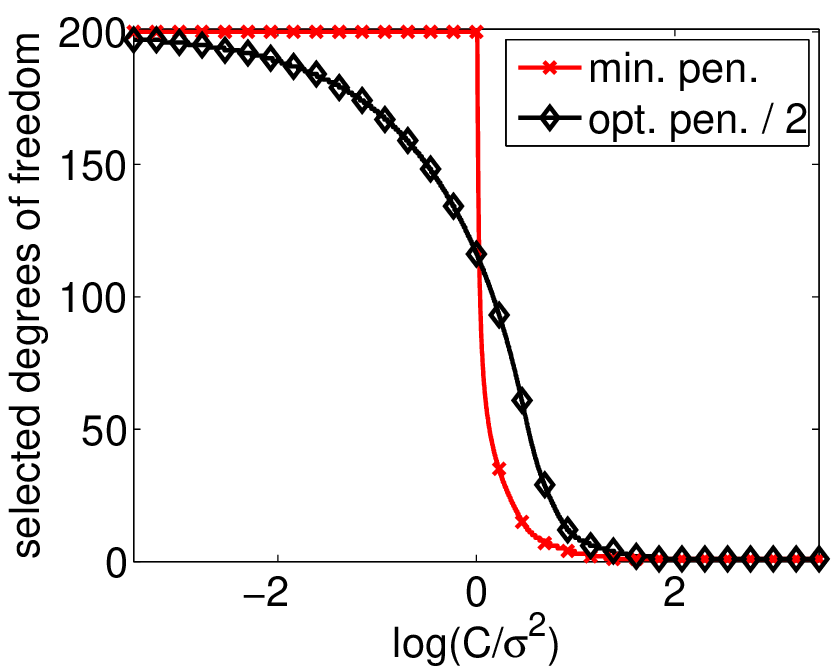}\hspace*{-1cm}

\hspace*{-1cm}\includegraphics[scale=.485]{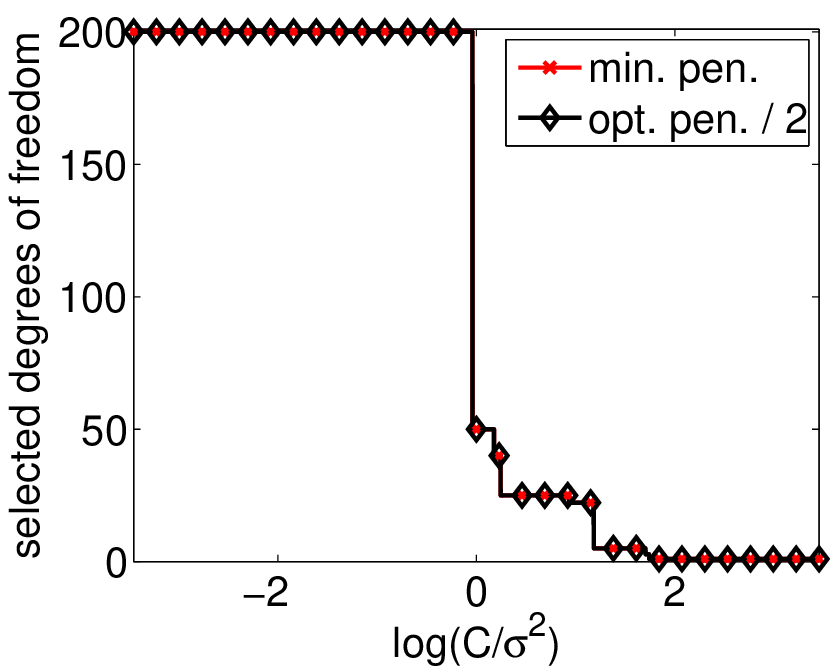} \hspace*{.36cm}
\includegraphics[scale=.485]{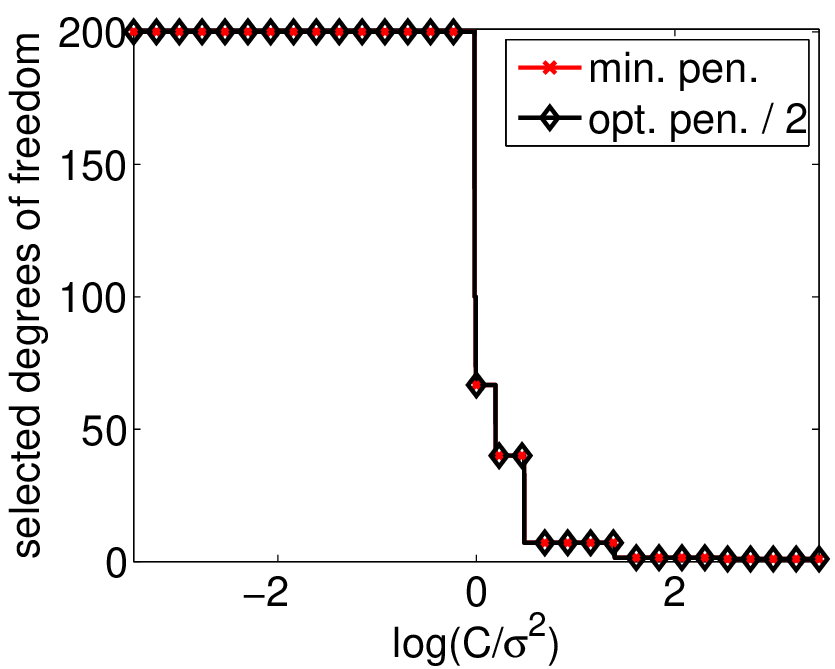}
\hspace*{-1cm}
\end{center}
%

\caption{Selected degrees of freedom vs. penalty strength $\log(C/\sigma^2)\,$, for a a fixed design problem and kernel ridge regression (top), Nadaraya-Watson estimator (middle) and $K$-nearest neighbor regression: note that when penalizing by the minimal penalty, there is a strong jump at $C = \sigma^2\,$, while when using half the optimal penalty, this is not the case. Left: $Y_i = \sin (25\pi X_i) + \varepsilon_i$, Right: $Y_i = \sin (25\pi X_i^3) + \varepsilon_i$. }
\label{fig:jump:2}
\end{figure}

\paragraph{Jump.}   
In Figure~\ref{fig:jump:2} we study the size of the jump for kernel ridge regression. With half the optimal penalty (which is used in traditional variable selection for linear regression), we do not get any jump, while with the minimal penalty we always do. Note that on the left plots, we may get sharp (but of lower amplitude) jumps away from $C = \sigma^2$. In this particular situation, this is due to the periodicity of the sine function.

\begin{figure}
%
%
\begin{center}
\hspace*{-1cm}
\includegraphics[scale=.485]{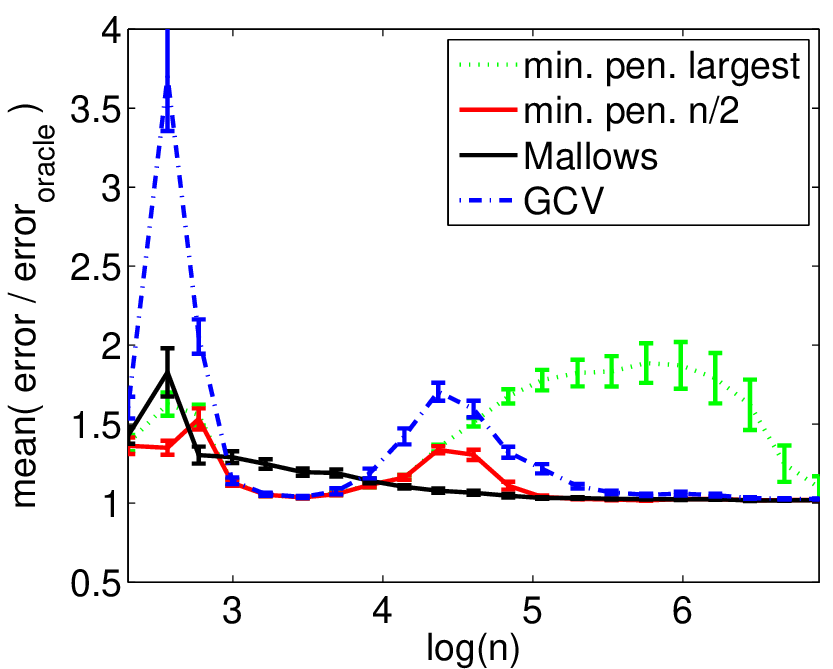} \hspace*{.36cm}
\includegraphics[scale=.485]{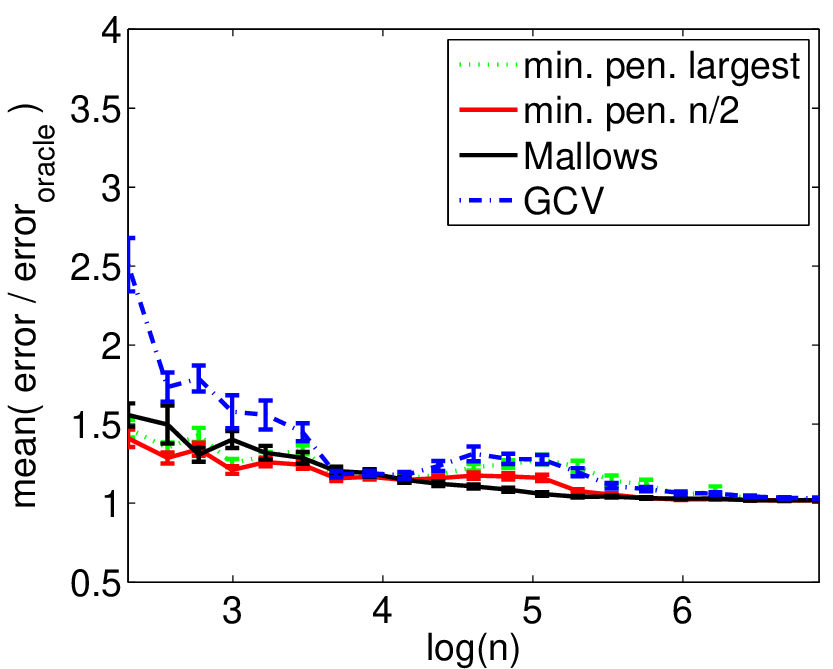} \hspace*{-1cm}

\hspace*{-1cm}\includegraphics[scale=.485]{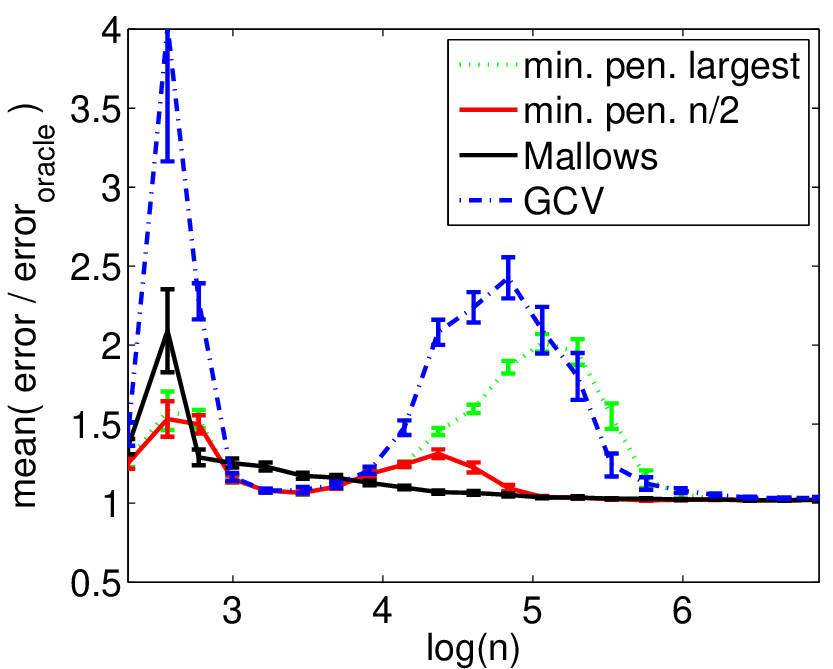} \hspace*{.36cm}
\includegraphics[scale=.485]{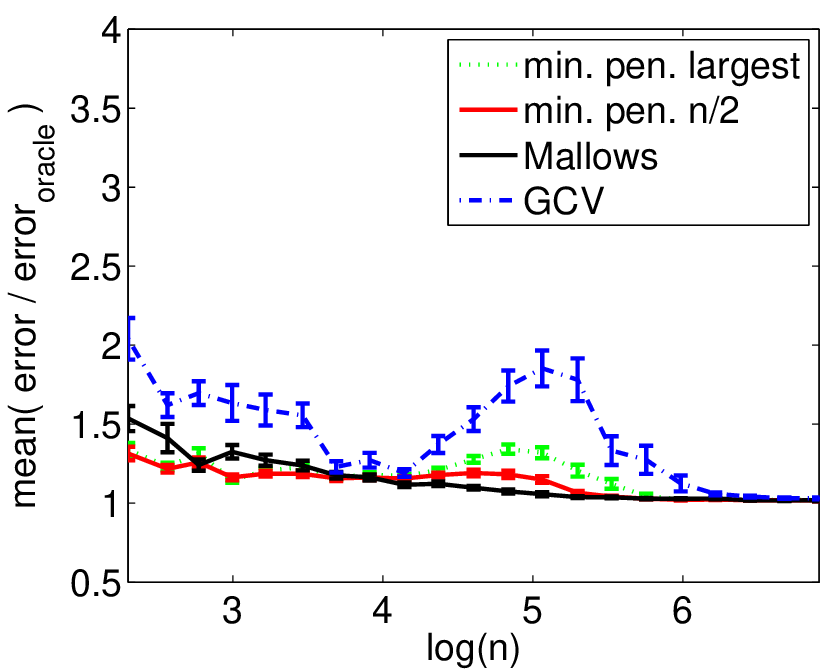} \hspace*{-1cm}

\hspace*{-1cm}\includegraphics[scale=.485]{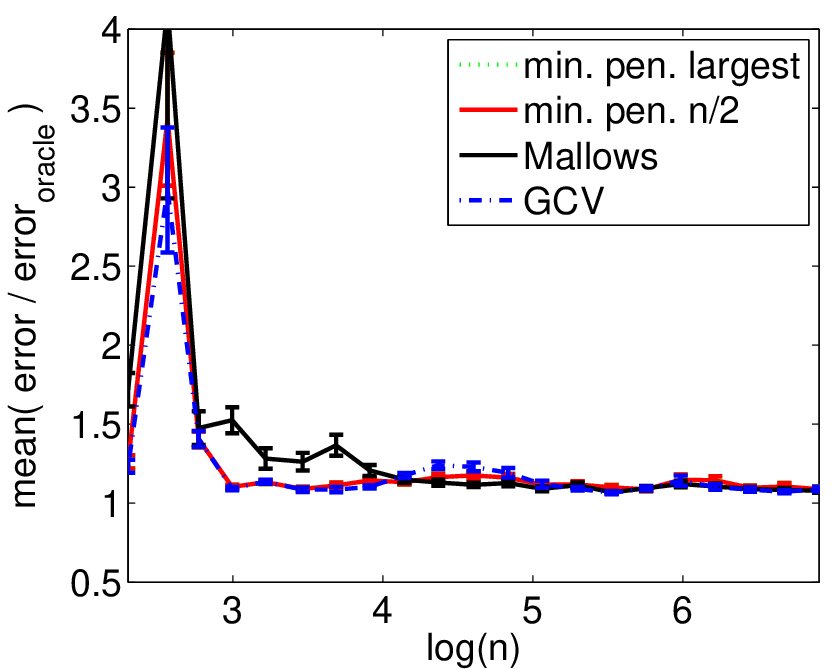} \hspace*{.36cm}
\includegraphics[scale=.485]{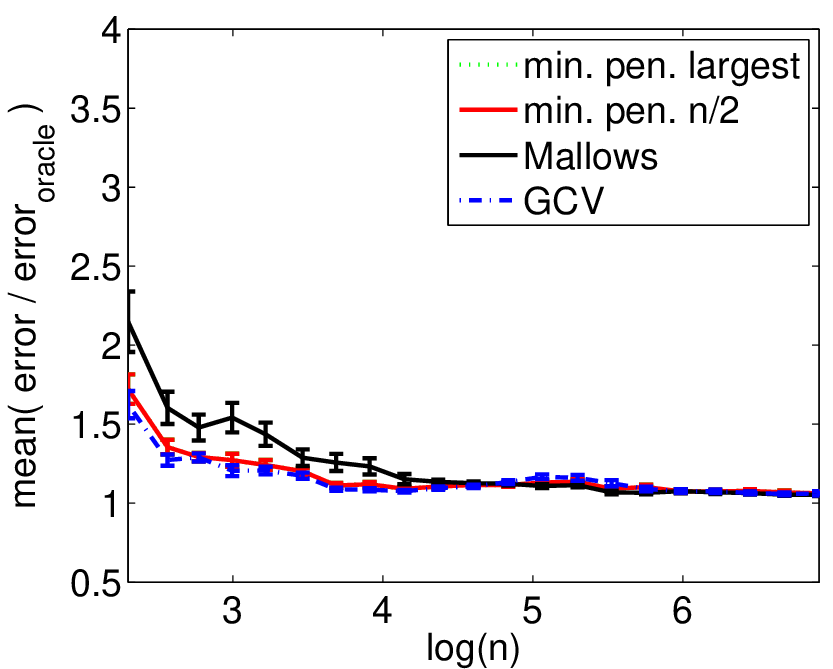}  
\hspace*{-1cm}
\end{center}
%

\caption{Comparison of various smoothing parameter selection (minimal with two types of jump selection, GCV, Mallows) for various values of numbers of observations. Left: $Y_i = \sin (25\pi X_i) + \varepsilon_i$, Right: $Y_i = \sin (25\pi X_i^3) + \varepsilon_i$. }
\label{fig:jump:3}
\end{figure}

\paragraph{Comparison of estimator selection methods.}   
In Figure~\ref{fig:jump:3}, we plot model selection results for 20 replications of data, comparing GCV \cite{Cra_Wah:1979} and Mallows $C_L$ penalty (which assumes the knowledge of $\sigma^2$). For the minimal penalty, we consider two strategies for finding the largest jump. We select the set of parameters so that the degrees of freedom are integers in $[0,n]$.
In the first strategy, we select the value $C$ with the largest jump, while in the second strategy we select the first $C$ so that the selected degrees of freedom goes below $n/2$.

We compare to the oracle (which can be computed because we can enumerate $\Lambda$). We see in Figure~\ref{fig:jump:3} that (a) the largest jump is not a good heuristic as spurrious jumps may occur, (b) the minimal penalty technique outperforms GCV, and (c) that in some cases, Mallows (which assumes more knowledge) is outperformed by our minimal penalty strategy. 
Indeed, when the signal-to-noise ratio is small, overpenalizing a bit (i.e., multiplying Mallows' penalty by a factor $\kappa>1$) is often useful, which turns out to be done automatically with the minimal penalty.

\section{Conclusion} \label{sec.conclu}
\paragraph{A new light on the slope heuristics.}   
Theorems~\ref{thm.mini}, \ref{thm.oracle.lhopt} and \ref{thm.algo} generalize some results first proved in \cite{Bir_Mas:2006} where all $\Al$ are assumed to be projection matrices.
To this extent, Birg\'e and Massart's slope heuristics has been modified in a way that sheds a new light on the ``magical'' factor~2 between the minimal and the optimal penalty, as proved in \cite{Bir_Mas:2006,Arl_Mas:2009:pente}.
Indeed, Theorems~\ref{thm.mini}--\ref{thm.oracle.lhopt} show that for general linear estimators,
\begin{equation} \label{eq.slope.gal} \frac{ \penid(\lambda) } {\penmin(\lambda)} = \frac{ 2 \tr(\Al) } { 2 \tr(\Al) - \tr(\Al^\top \Al)} \enspace , \end{equation}
which can take any value in $(1,2]$ in general (assuming $\cAl=1$ as in all the major examples we have in mind); this ratio is only equal to~2 when $\tr(\Al) \approx \tr(\Al^\top \Al)\,$, hence mostly when $\Al$ is a projection matrix or a $k$-NN matrix.

\paragraph{Covariance matrix estimation.} 
A natural extension of the present work appears in the multitask regression example, when $p$ regression problems (the tasks) are solved simultaneously. Then, a key quantity is the $p \times p$ covariance matrix of the tasks, which has to be estimated for an optimal selection of regularization parameters (for instance). As shown in \cite{Sol_Arl_Bac:2011}, the minimal penalty strategy can be used successfully for estimating a covariance matrix under mild assumptions, by applying the algorithm of Section~\ref{sec.algo.def} to $p(p+1)/2$ well-chosen one-dimensional regression problems. 

\paragraph{Future directions.}   
The good empirical performances of elbow heuristics based algorithms (i.e., based on the sharp variation of a certain quantity around good hyperparameter values \cite{Han_OLe:1993,Gro_Wol:2009,Rez_Hos:2009}) suggest that Theorem~\ref{thm.algo} can be generalized to many learning frameworks (and potentially to non-linear estimators), probably with small modifications in the algorithm, but always relying on the concept of minimal penalty.

\medskip

In the case of projection estimators, the slope heuristics still holds when the design is random and data are heteroscedastic \cite{Arl_Mas:2009:pente}; 
we conjecture a generalization of Eq.~\eqref{eq.slope.gal} is still valid for heteroscedastic data with some (but not all) linear estimators.

Another interesting open problem would be to extend the results of Section~\ref{sec.algo} to more general continuous sets $\Lambda$, such as the ones appearing naturally in multiple kernel learning.
We conjecture that Theorem~\ref{thm.algo} is valid without modification for a ``small'' continuous $\Lambda\,$, that is, of ``small'' dimension. 
On the contrary, in applications such as the Lasso with $p \gg n$ variables, the natural set $\Lambda$ cannot be well covered by a grid of cardinality $n^{\alpha}$ with $\alpha$ small, and our minimal penalty algorithm and Theorem~\ref{thm.algo} certainly have to be modified.


\appendix
\section*{Appendix}

\section{Notation and first computations} \label{A.sec.nota}
Recall that 
\[ Y = F + \varepsilon \]
where $F = (f(x_i))_{1 \leq i \leq n} \in \R^n$ is deterministic, $\varepsilon = (\varepsilon_i)_{1 \leq i \leq n}\in \R^n$ is centered with covariance matrix $\sigma^2 \Id_n$ and $\Id_n$ is the $n \times n$ identity matrix.

For every $x \in \R\,$, $x_+ = \max\set{x \, , \, 0}$ denotes the positive part of $x\,$.

In the proofs, we use repeatedly that 
\begin{equation} 
\label{eq.maj-theta}
\forall a,b \geq 0 \, , \, \forall \theta >0 \, , \quad 2 \sqrt{ab} \leq \theta a + \theta^{-1} b \enspace ,
\end{equation}
with equality for $\theta = \sqrt{b/a}$.

\subsection{General framework}

For every $\lamm$, $\Fhl = \Al Y$ for some $n \times n$ real-valued matrix $\Al\,$, so that  
\begin{align} \label{A.eq.riskFhl}
\norm{ \Fhl - F}_2^2 &= \norm{ (\Al - \Id_n) F}_2^2 + \norm{ \Al \varepsilon}_2^2 + 2 \prodscal{\Al \varepsilon}{ (\Al - \Id_n) F} \enspace , \\ \label{A.eq.riskempFhl}
\norm{\Fhl - Y}_2^2 &= \norm{ \Fhl - F}_2^2 + \norm{\varepsilon}_2^2 - 2 \prodscal{\varepsilon}{\Al \varepsilon} + 2 \prodscal{\varepsilon}{(\Id_n - \Al) F} \enspace , 
\end{align} 
where $\forall t,u \in \R^n$, $\prodscal{t}{u} = \sum_{i=1}^n t_i u_i$ and $\norm{t}_2^2 = \prodscal{t}{t}\,$.

Now, define, for every $\lamm\,$,
\begin{gather*}
b(\lambda) = \norm{ (\Al - \Id_n) F}_2^2 
\qquad 
\vun(\lambda) = \tr(\Al) \sigma^2
\qquad 
\vdeux(\lambda) = \tr(\Al^\top \Al) \sigma^2 
\\
\delun(\lambda) = \prodscal{\varepsilon}{\Al \varepsilon} - \tr(\Al) \sigma^2 
\qquad
\deldeux(\lambda) = \norm{ \Al \varepsilon}_2^2 - \tr(\Al^\top \Al) \sigma^2
\\
\delta_3(\lambda) = 2 \prodscal{\Al \varepsilon}{ (\Al - \Id_n) F} 
\qquad 
\delta_4(\lambda) =  2 \prodscal{\varepsilon}{(\Id_n - \Al) F} \enspace ,
\end{gather*}
so that Eq.~\eqref{A.eq.riskFhl} and~\eqref{A.eq.riskempFhl} can be rewritten
\begin{align} \label{A.eq.decomp.riskFhl}
\norm{ \Fhl - F}_2^2 &= b(\lambda) + \vdeux(\lambda) + \deldeux(\lambda) + \delta_3(\lambda) \\
\label{A.eq.decomp.riskempFhl}
\norm{ \Fhl - Y}_2^2 &= \norm{ \Fhl - F}_2^2 - 2 \vun(\lambda) - 2 \delun(\lambda) + \delta_4(\lambda) + \norm{\varepsilon}_2^2  \enspace .
\end{align}

Note that $b (\lambda)$, $\vun(\lambda)$ and $\vdeux(\lambda)$ are deterministic, and  for all $\lamm$ and $i=1\ldots 4\,$, $\delta_i(\lambda)$ is random with zero mean. In particular, we deduce the following expressions of the risk and the empirical risk of $\Fhl\,$:
\begin{align} \label{A.eq.EriskFhl}
\E\croch{ n^{-1} \norm{\Fhl - F}_2^2 } &= n^{-1} \norm{ (\Al - \Id_n) F}_2^2 + \frac{ \tr (\Al^\top \Al) \sigma^2} {n} \enspace , \\
\label{A.eq.EriskempFhl}
\E\croch{ n^{-1} \norm{\Fhl - Y}_2^2 } - \sigma^2 &= n^{-1} \norm{ (\Al - \Id_n) F}_2^2 - \frac{ \paren{ 2 \tr(\Al) - \tr (\Al^\top \Al) } \sigma^2 } {n}   \enspace . 
\end{align}

\subsection{The event $\Omega_x$} \label{A.sec.conc.event}

In this section, we define the large probability event on which all our main results hold. 
Let $\Comega \in [0,+\infty)^6$ be fixed. 
Then, for any $x \geq 0\,$, we define the event 
\[ \Omega_x = \Omega_x(\Lambda) = \Omega_x \paren{ \Lambda, \Comega} \] 
on which, for every $\lamm$ and every $\theta_1, \theta_2, \theta_3, \theta_4 \in(0,1]\,$,
\begin{align}
\label{A.eq.conc.del1}
\absj{\delun(\lambda)} &\leq \theta_1 \sigma^2 \tr(\Al^{\top} \Al) + \paren{\ComegaA + \ComegaB \theta_1^{-1} } x \sigma^2  \\
\label{A.eq.conc.del2}
\absj{\deldeux(\lambda)} &\leq \theta_2 \sigma^2 \tr(\Al^{\top} \Al) + \paren{\ComegaC + \ComegaD \theta_2^{-1} } x \sigma^2 \\
\label{A.eq.conc.del3}
\absj{\delta_3(\lambda)} &\leq \theta_3 \norm{ (\Id_n - \Al) F}_2^2 + \ComegaE \theta_3^{-1} x \sigma^2 \\
\label{A.eq.conc.del4}
\absj{\delta_4(\lambda)} &\leq \theta_4 \norm{ (\Id_n - \Al) F}_2^2 + \ComegaF \theta_4^{-1} x \sigma^2 \enspace . 
\end{align}

\subsection{Splitting assumption \eqref{hyp.Lam}}

We actually deal with assumption \eqref{hyp.Lam} by considering separately the case of discrete $\Lambda$ (assumption~\eqref{hyp.Lam-dis}) and the case of ridge regression with a continuous one-dimensional parameter (assumption~\eqref{hyp.ridge}):
\begin{itemize}
\item The set $\Lambda$ is finite (with a polynomial size w.r.t. $n$): 
\begin{align}
\tag{\ensuremath{\mathbf{H\Lambda{}dis}}}
\label{hyp.Lam-dis} 
\card(\Lambda) \leq \CLamdis n^{\alLamdis}
\end{align}
\item
In the kernel ridge regression example:
\begin{equation}
\tag{\ensuremath{\mathbf{Hridge}}}
\label{hyp.ridge} 
\left.
\begin{aligned}
\Lambda = [0,+\infty]
 \qquad 
A_0 = \Id_n 
\qquad A_{\infty} = \zeromat{n} 
\qquad \\
\exists K \in \M_n(\R) \backslash \set{ \zeromat{n} } \mbox{ symmetric positive semi-definite}  \\
\mbox{such that} \quad \forall \lambda \in (0,+\infty) \, ,  \quad \Al = K (K+n\lambda \Id_n)^{-1} 
\end{aligned}
\right\}
\end{equation}
\end{itemize}
Once the concentration results will be proved under each assumption, an union bound will yield the desired results under the composite assumption \eqref{hyp.Lam}.

\subsection{Kernel ridge regression} \label{A.sec.ridge}
Let us now prove a few useful elementary results that are specific to kernel ridge regression, that is, when assumption \eqref{hyp.ridge} holds true.

Under assumption~\eqref{hyp.ridge}, a key remark is the following.
Since $K$ is symmetric positive, non-negative numbers $\mu_1, \ldots, \mu_n$ and some orthogonal matrix $P \in \M_n(\R)$ exists such that $K=P^{\top} \diag(\mu_1, \ldots , \mu_n) P\,$. 
Then, 
\begin{equation}
\label{eq.Omega.ridge.key} \forall \lamm \, , \quad \Al = P^{\top} \Dl P \quad \mbox{where} \quad \Dl = \diag\paren{ \paren{ \frac{\mu_j}{\mu_j + n \lambda} }_{1 \leq j \leq n} } \enspace , 
\end{equation}
with the convention $D_0 = \Id_n$ and $D_{\infty} = \zeromat{n}\,$.
We also define 
\[ (f_j)_{1 \leq j \leq n} = PF \enspace .\]

%
\begin{lemma} \label{le.ridge.monotone}
If assumption~\eqref{hyp.ridge} holds true, then, 
\[ \tr(\Al) \enspace , \quad \tr(\Al^{\top} \Al) \quad \mbox{and} \quad 2 \tr(\Al) - \tr(\Al^{\top} \Al) \]
are decreasing continuous functions of $\lambda$ over $[0,+\infty]\,$, all equal to $n$ for $\lambda=0$ and to $0$ for $\lambda = + \infty\,$\textup{;}
$b(\lambda)$ is a nondecreasing continuous function of $\lambda$, with $b(0)=0$ and $b(+\infty)=\norm{F}^2\,$.
\end{lemma}

\begin{proof}[Proof of Lemma~\ref{le.ridge.monotone}]
According to Eq.~\eqref{eq.Omega.ridge.key}, for every $\lambda \in [0,+\infty)$,
\begin{align*}
\tr(\Al) &= \dfl = \sum_{j=1}^n \paren{ \frac{ \mu_j } { \mu_j + n \lambda} } \\
\tr(\Al^\top \Al)  &= \sum_{j=1}^n \paren{ \frac{ \mu_j } { \mu_j + n \lambda} }^2 \\
2 \tr(\Al) - \tr(\Al^\top \Al) 
&=  \sum_{j=1}^n \croch{ \frac{ 2 \mu_j } { \mu_j + n \lambda} - \paren{ \frac{ \mu_j } { \mu_j + n \lambda} }^2 } 
= \sum_{j=1}^n \croch{ \frac{ \mu_j (\mu_j + 2 n \lambda) } { \paren{ \mu_j + n \lambda}^2 } } \\
b(\lambda) &= \norm{ (\Al - \Id_n) F}_2^2 = \sum_{j=1}^n \paren{ 1 - \frac{ \mu_j } { \mu_j + n \lambda} }^2 f_j^2 \enspace , 
\end{align*}
The result follows since $\mu_j \geq 0$ for every $j\,$. 
\end{proof}

\subsection{A useful lemma} \label{A.sec}

\begin{lemma} \label{le.retire.condition.A3}
Let $n \geq 1$ an integer. 
Then, for any matrix $A \in \M_n(\R)$, 
\begin{equation} \label{eq.le.retire.condition.A3.base}
\tr(A) \leq \sqrt{ n \tr(A^{\top}A) } \enspace ,
\end{equation}
from which we get
\begin{align} \label{eq.le.retire.condition.A3.v1}
&\forall A \in \M_n(\R) \, , \, x \geq 0 \, , \, \theta>0 \, , \quad 
x \tr(A) \leq \theta \tr(A^{\top} A) + \frac{x^2 n }{4 \theta} 
\enspace . 
\end{align}
\end{lemma}
\begin{proof}[Proof of Lemma~\ref{le.retire.condition.A3}]
First, since $(A,B) \mapsto \tr(A^{\top} B)$ is a scalar product on $\M_n(\R) \,$, by Cauchy-Schwarz inequality, for every $A \in \M_n(\R)$, 
\[ \paren{\tr(A)}^2 = \paren{\tr(\Id^{\top} A)}^2 \leq \tr\paren{A^{\top} A} \tr\paren{\Id^{\top} \Id} =  n \tr\paren{A^{\top} A} \enspace . \]
Therefore, for every $x \geq 0$, $\theta>0$, 
\[ 
x \tr(A) \leq 2 \sqrt{ \tr(A^{\top} A) \frac{x^2 n}{4}} \leq \theta \tr(A^{\top} A) + \frac{x^2 n}{4 \theta}  \qquad \mbox{by Eq.~\eqref{eq.maj-theta}.}
\]
\end{proof}

\section{Key concentration inequalities} \label{A.sec.conc}

In this section, we state the concentration inequalities showing the event $\Omega_x$ has a large probability. 
Our main contributions are the following. Proposition~\ref{A.pro.conc.quad.gal.Gauss.momexp.improved} slightly extends a result by Laurent and Massart \cite{Lau_Mas:2000} for the concentration of quadratic forms of a Gaussian vector. 
In the kernel ridge regression example, we prove uniform concentration inequalities (over a continuous set) for some linear forms of a Gaussian vector (Proposition~\ref{pro.Omega.ridge.del34.Gauss} and Lemma~\ref{le.conc-linear.Gauss.general}), and for some quadratic forms of a random vector (Proposition~\ref{pro.Omega.ridge.del12}). 
All results are proved in Section~\ref{A.sec.proof-conc}. 

\subsection{Linear functions of $\varepsilon$} \label{A.sec.conc.lin-Gauss}
In the Gaussian case, concentration inequalities for $\delta_3(\lambda)$ and $\delta_4(\lambda)$ come from the following standard result.
\begin{proposition} \label{A.pro.conc.linear.gal.Gauss}
Let $\xi$ be a standard Gaussian vector in $\R^n\,$, $\alpha\in \R^n$ and 
$ Z = \prodscal{ \xi}{\alpha} = \sum_{j=1}^n \alpha_j \xi_j $. 
Then, for every $x \geq 0\,$,
\begin{equation*} 
\Prob\paren{ \absj{Z} \leq \sqrt{ 2x } \norm{\alpha}_2 } \geq 1 - e^{-x} \enspace . 
\end{equation*}
\end{proposition} 
\begin{proof}[Proof of Proposition~\ref{A.pro.conc.linear.gal.Gauss}]
The result is clear when $\alpha=0$. Otherwise, $\prodscal{ \xi}{\alpha} \sim \mathcal{N} \paren{ 0 , \norm{\alpha}^2}$ which gives the result.  
\end{proof}

In the kernel ridge regression case, we get a similar result (with larger constants) uniformly over $\Lambda = [0,+\infty]$. Up to the best of our knowledge, Proposition~\ref{pro.Omega.ridge.del34.Gauss} is a new result that could be useful for studying kernel ridge regression in a more general framework. 

\begin{proposition} \label{pro.Omega.ridge.del34.Gauss}
Assume that \eqref{hyp.ridge} and \eqref{hyp.eps-Gauss-hom} hold true. 
Then, 
for every $x \geq 0$, an event $\Omegaridge_x$ of probability at least $1 - \exp(-x + \betaLRUc + \ln(n)) $ exists on which 
for every $\lamm\,$, 
\begin{align}
\label{eq.conc.del3.ridge-unif.Gauss}
\absj{\delta_3(\lambda)} &\leq \cteLinRidge \sigma \norm{(\Id_n-\Al) F} \sqrt{x}, \\
\label{eq.conc.del4.ridge-unif.Gauss}
\absj{\delta_4(\lambda)} &\leq \cteLinRidge \sigma \norm{(\Id_n-\Al) F} \sqrt{x} 
\enspace . 
\end{align}
\end{proposition}
Proposition~\ref{pro.Omega.ridge.del34.Gauss} is proved in Section~\ref{A.sec.proof.le.Omega.ridge.del34.Gauss}.

\subsection{Quadratic functions of $\varepsilon$} \label{A.sec.conc.quad-Gauss}

In the Gaussian case, concentration inequalities for $\delun(\lambda)$ and $\deldeux(\lambda)$ come from the following result.
\begin{proposition} \label{A.pro.conc.quad.gal.Gauss.momexp.improved}
Let $\xi$ be a standard Gaussian vector in $\R^n\,$, $M$ a real-valued $n \times n$ matrix and 
$ Z \egaldef \prodscal{\xi}{M \xi} - \tr(M) = \xi^{\top} M \xi - \tr(M)\,$.
Then, $\E\croch{Z}=0$ and for every $x \geq 0\,$,
\begin{equation} \label{A.eq.pro.conc.quad.gal.Gauss.momexp.maj.precis.improved}
\Prob\paren{ Z \leq \sqrt{ 2 x \paren{ \tr(M^2) + \tr(M^{\top}M) }} + 2 \normmat{M} x } 
\geq 1 - \exp\paren{ - x } \enspace . 
\end{equation}
\end{proposition}
Since $\tr(M^2) \leq \tr(M^{\top}M)$ (by Lemma~\ref{le.trM^2<trMTM}), we get that for every $x \geq 0\,$,
\begin{align} \label{A.eq.pro.conc.quad.gal.Gauss.momexp.maj.improved}
\Prob\paren{  \prodscal{\xi}{M \xi} \leq \tr(M) + 2 \sqrt{ x \tr(M^{\top}M)} + 2 \normmat{M} x } 
&\geq 1 - \exp\paren{ - x } 
\enspace .  
\end{align}
Proposition~\ref{A.pro.conc.quad.gal.Gauss.momexp.improved} extends \cite[Lemma~1]{Lau_Mas:2000}; it is proved in Section~\ref{A.sec.proof.pro.conc.quad.gal.Gauss.momexp.improved}. 
The main deviation term in Eq.~\eqref{A.eq.pro.conc.quad.gal.Gauss.momexp.maj.precis.improved} is optimal since 
$ \var(Z) =  \tr(M^2) + \tr(M^{\top} M) $ as shown in Section~\ref{sec.supmat.var-quad-form}. 

\medskip

In the kernel ridge regression case, we get a similar result uniformly over $\Lambda = [0,+\infty]$. Up to the best of our knowledge, Proposition~\ref{pro.Omega.ridge.del12} is a new result that could be useful for studying kernel ridge regression in a more general framework.

\begin{proposition} \label{pro.Omega.ridge.del12}
Assume that \eqref{hyp.ridge} holds true. 
Then, assumption~\eqref{hyp.Al} is satisfied with $\majnormAl=\cAl=1$, and some $\Lambda_1 \subset \Lambda$ exists such that 
$\card(\Lambda_1) \leq 2 n$ and 
for every $\Comega \in [0,+\infty)^6\,$ and $x \geq 0\,$, 
on $\Omega_x \paren{ \Lambda_1, \Comega }\,$, 
for every $\theta_1,\theta_2 \in (0,1]\,$, for every $\lamm\,$, 
\begin{align}
\label{eq.conc.del1.ridge-unif}
\absj{\delun(\lambda)} &\leq \theta_1  \sigma^2 \tr(\Al^{\top}\Al) + \paren{ \ComegaA + \ComegaB \theta_1^{-1}} x \sigma^2 + 2 \sigma^2 \\
\label{eq.conc.del2.ridge-unif}
\absj{\deldeux(\lambda)} &\leq \theta_2  \sigma^2 \tr(\Al^{\top}\Al) + \paren{ \ComegaC + \ComegaD \theta_2^{-1}} x \sigma^2 + 2 \sigma^2 
\enspace . 
\end{align}
\end{proposition}
\begin{remark}
Proposition~\ref{pro.Omega.ridge.del12} does not rely on any assumption on the distribution of the noise $\varepsilon\,$. Therefore, it can be used in the Gaussian case \textup{(}combined with Lemma~\ref{le.Omega.finite}\textup{)}, but also under any alternative assumption for which concentration inequalities for $\delun(\lambda)$ and $\deldeux(\lambda)$ can be proved. 
\end{remark}
Proposition~\ref{pro.Omega.ridge.del12} is proved in Section~\ref{A.sec.proof.le.Omega.ridge.del12}.

\subsection{Lower bounds on the probability of $\Omega_x(\Lambda)$} \label{A.sec.conc.finite}
We are now in position to state lower bounds on the probability of $\Omega_x(\Lambda,\Comega)$ provided $\Comega$ is well-chosen. First, we consider the case of a finite $\Lambda$. 

\begin{lemma} \label{le.Omega.finite}
Under assumptions \eqref{hyp.Lam-dis}, \eqref{hyp.Al} and \eqref{hyp.eps-Gauss-hom}, $\Prob(\Omega_x(\Lambda,\Comega)) \geq 1 - 6 \card(\Lambda) e^{-x}$ 
\begin{gather*}
\mbox{with} \quad \ComegaA = 2 \majnormAl \quad \ComegaB = 1 \quad \ComegaC = 2 \majnormAl^2 \quad \ComegaD = \majnormAl^2  \quad \ComegaE = 2 \majnormAl^2 \quad \mbox{and} \quad  \ComegaF = 2 \enspace . 
\end{gather*}
\end{lemma}
Lemma~\ref{le.Omega.finite}, a consequence of Propositions~\ref{A.pro.conc.linear.gal.Gauss} and~\ref{A.pro.conc.quad.gal.Gauss.momexp.improved}, is proved in Section~\ref{A.sec.proof.le.Omega.finite}. 

\medskip

In the kernel ridge regression example, we prove a similar concentration result (with slightly larger constants) uniformly over the continuous set $\Lambda=[0,+\infty]$.
\begin{lemma} \label{le.Omega.cont-ridge}
Under assumptions \eqref{hyp.ridge} and \eqref{hyp.eps-Gauss-hom}, 
$\Prob(\Omega_x(\Lambda,\Comega)) \geq 1 - \exp\paren{\betaLRUd+ \ln(n)} e^{-x}$
\begin{gather*}
\mbox{if} \quad 
\ComegaA = 2 \quad \ComegaB = 1 \quad \ComegaC = 2 \quad \ComegaD =  1 \quad \ComegaE = 
\cteC \quad \mbox{and} \quad  \ComegaF = \cteC \enspace . 
\end{gather*}
\end{lemma}
Lemma~\ref{le.Omega.cont-ridge}, a consequence of Propositions~\ref{pro.Omega.ridge.del34.Gauss}, \ref{A.pro.conc.quad.gal.Gauss.momexp.improved} and~\ref{pro.Omega.ridge.del12}, is proved in Section~\ref{A.sec.proof.le.Omega.cont-ridge}.

\section{Proof of Proposition~\ref{pro.hyp.Al}} \label{sec.pr.pro.hyp.Al}
We consider separately the examples to which Proposition~\ref{pro.hyp.Al} applies. 
\subsection{Case (i)}
Since $A$ is symmetric, it can be diagonalized in an orthonormal basis, with eigenvalues $a_1, \ldots , a_n \in [0,1]$ (by assumption), so that 
\[ \tr(A^{\top}A) = \sum_{i=1}^n a_i^2 \leq \sum_{i=1}^n a_i = \tr(A) \leq n \enspace . \]
In particular, in example (ia), $A$ is an orthogonal projection matrix, so $A$ is symmetric with $\Sp(A) \subset [0,1]$ and $A^{\top}A = A$ implies $\tr(A^{\top} A) = \tr(A)\,$.
For example (ib), by Eq.~\eqref{eq.Omega.ridge.key}, $A$ is symmetric and $\Sp(A) \subset [0,1]\,$. 

\subsection{Case (ii)}
\[
\tr(A^{\top} A) 
= \sum_{1 \leq i,j \leq n} A_{i,j}^2 
\leq \sum_{i=1}^n \sum_{j=1}^n A_{i,j} A_{i,i} = \sum_{i=1}^n A_{i,i} = \tr(A) 
\leq \sum_{i=1}^n \sum_{k=1}^n A_{i,k} = n \enspace .
\]
In particular, if \eqref{eq.kNN} holds true, then $\tr(A) = n/k$ 
and \[
\tr(A^{\top} A) 
= \sum_{1 \leq i,j \leq n} A_{i,j}^2 
= \frac{1}{k}  \sum_{1 \leq i,j \leq n} A_{i,j} = \frac{n}{k} \enspace .
\qed \]

\section{Minimal penalty (proof of Theorem~\ref{thm.mini})} \label{A.sec.thm.mini}
Let us recall the definition \eqref{eq.lhmin} of $\lhminb(C)$:
\begin{equation*} 
\forall C \geq 0, \qquad \lhminb(C) \in \arg\min_{\lamm} \set{ \norm{\Fhl - Y}_2^2 + C \paren{ 2 \tr(\Al) - \tr(\Al^\top \Al)} } \enspace . \end{equation*}
We prove in this section Theorem~\ref{thm.mini}, which is actually a corollary of the following proposition: 
\begin{proposition} \label{A.pro.mini}
Let $\lhminb$ be defined by Eq.~\eqref{eq.lhmin}. 
Assume that $\Id \in \set{\Al}_{\lL}$, \eqref{hyp.Al} holds true, and let $\majdimbiais \in [0,n]$ and  $\majbiais \geq \sqrt{n \ln(n)}$ be such that 
\begin{gather} 
\label{A.hyp.thm.mini.2.bis} \tag{$\mathbf{Abiais^{\prime}}$}
\exists \lhypbias \in \Lambda, \quad \df(\lhypbias) \leq  \majdimbiais \quad \mbox{and} \quad b(\lhypbias) \leq \sigma^2 \majbiais
\enspace . 
\end{gather}
Let $\Comega \in [0,+\infty)^6$ and define 
\[ \KthetaA \egaldef  2 \ComegaA  + \ComegaC + 3 \ComegaE + 3 \ComegaF \qquad \mbox{and} \qquad \KthetaB \egaldef 6 \ComegaB + 3 \ComegaD \enspace . \]
Then, for every $\gamma \geq 2$, $\ctepenminB>0$, a constant $\nminProMini>0$ only depending on $\KthetaA,\KthetaB,\gamma,\ctepenminB$ exists such that if $n \geq \nminProMini \,$, for every  $\mindimavant \in [ 0 , \frac{n}{2} )$ and $\majdimapres \in \scroch{ \frac{10}{\cAl} \max\sset{ \KthetaB \gamma \ln(n) , 2 \majdimbiais} , n }$
\begin{align} 
\label{A.eq.thm.lhmin.below.bis} 
\forall 0 \leq C < \croch{ 1 -  \paren{ 2 \sqrt{2 \KthetaB} + \frac{1}{\ctepenminB}  } \frac{  \gamma \sqrt{n \ln(n)}} {n-2\mindimavant}}  \sigma^2 , \quad \df(\lhminb(C)) &\geq \mindimavant
\\
\label{A.eq.thm.lhmin.above.bis}
\forall C > \croch{ 1 + \paren{ \frac{25}{27} + \frac{1}{\ctepenminB} + \frac{40}{9} \sqrt{\KthetaB}} \frac{\gamma \majbiais } {\cAl \majdimapres} } \sigma^2 , \quad \df(\lhminb(C)) &\leq \majdimapres  \enspace .
\end{align}
hold on the event $\Omega_{\gamma \ln(n)}(\Lambda,\Comega)\,$.

In particular, under assumptions \eqref{hyp.Lam} and \eqref{hyp.eps-Gauss-hom}, if $n \geq \nminProMiniB(\gamma,\majnormAl)$, with probability at least $1 - \paren{ 6 \card(\Lambda_0) + \NLamfusrid \exp\paren{\betaLRUd + \ln(n)}} n^{-\gamma}\,$, 
\begin{align} 
\label{A.eq.thm.lhmin.below.bis.cas1} 
\forall \mindimavant \in \left[ 0 , \frac{n}{2} \right) \, , \, 
\forall 0 \leq C < \croch{ 1 -  \frac{ \cteAb \majnormAl \gamma \sqrt{n \ln(n)}}{n - 2 \mindimavant}  }  \sigma^2 , \quad \df(\lhminb(C)) &\geq \mindimavant
\\
\label{A.eq.thm.lhmin.above.bis.cas1}
\forall \majdimapres \in \croch{ \frac{10}{\cAl} \max\set{ \KthetaB \gamma \ln(n) , 2 \majdimbiais} , n } \, , \, 
\forall C > \croch{ 1 + \frac{\cteBa \gamma \majbiais \majnormAl} {\cAl \majdimapres} } \sigma^2 , \quad \df(\lhminb(C)) &\leq \majdimapres  
\enspace .
\end{align}
\end{proposition}
\begin{proof}[Proof of Theorem~\ref{thm.mini}]
Taking $\majdimbiais = \sqrt{n}$ and $\majbiais = \sqrt{n \ln(n)}\,$, assumption \eqref{A.hyp.thm.mini.2.bis} becomes \eqref{hyp.thm.mini.2}. 
Let $\mindimavant=n/3$, $\majdimapres=n/10$ (which is possible in Proposition~\ref{A.pro.mini} as soon as $n/\ln(n) \geq 100 \KthetaB \gamma$ and $n \geq (200/\cAl)^2$). 
Then, choosing $\gamma = \alLamMix + \delta\,$, 
Eq.~\eqref{A.eq.thm.lhmin.below.bis} becomes Eq.~\eqref{A.eq.thm.lhmin.below}, 
Eq.~\eqref{A.eq.thm.lhmin.above.bis} becomes Eq.~\eqref{A.eq.thm.lhmin.above}, 
and they hold with probability at least 

\begin{align*}
&\qquad 1 - \paren{ 6 \card(\Lambda_0) + \NLamfusrid \exp\paren{\betaLRUd + \ln(n)}} n^{-\alLamMix - \delta}
\\
&\geq 
\detail{ \detailtag  1 - \paren{ 6 \CLamfusdis n^{\alLamfusdis} + \CLamfusrid \exp\paren{\betaLRUd + (\alLamfusrid+1) \ln(n)}} n^{-\alLamMix - \delta}  
   \\ &=   }
1 - \paren{ 6 \CLamfusdis n^{\alLamfusdis -\alLamMix} + \CLamfusrid \exp\paren{\betaLRUd +  ( \alLamfusrid + 1 -\alLamMix )  \ln(n)}} n^{ - \delta}
\\
&\geq 
\detail{   \detailtag 1 - \paren{ 6 \CLamfusdis + \CLamfusrid \exp\paren{\betaLRUd -  \ln(n)}} n^{ - \delta}   
\\ &\geq }
1 - \paren{ 6 \CLamfusdis  + \CLamfusrid} n^{ - \delta}
\end{align*}
as soon as $\ln(n) \geq \betaLRUd\,$, since $\alLamMix \geq 2 + \alLamfusrid$ and $\alLamMix \geq \alLamfusdis\,$.
\end{proof}

\begin{remark} \label{rk.mini.blowup}
On the event $\Omega_{\gamma \ln(n)}(\Lambda)$ where Eq.~\eqref{A.eq.thm.lhmin.below.bis} and~\eqref{A.eq.thm.lhmin.above.bis} hold and under the same assumptions, 
we can derive from Eq.~\eqref{A.eq.decomp.riskFhl}, \eqref{A.eq.conc.del2} with $\theta_2=1/2$, \eqref{A.eq.conc.del3} with $\theta_3=1$, that  
\[ \forall \lamm \, , \quad n^{-1} \norm{ \Fhl - F }_2^2 \geq \frac{\tr(\Al^{\top}\Al) \sigma^2}{2 n} - \frac{(\ComegaC+2\ComegaD+\ComegaE) \gamma \ln(n) \sigma^2 }{n} \enspace . \]
Since $\tr(\Al^{\top}\Al) \geq n^{-1} \paren{\dfl}^2\,$ we deduce that  for all $\lambda \in \Lambda$:
\[  \dfl \geq \frac{n}{\ln(n)} \ \Rightarrow \  n^{-1} \norm{ \Fhl - F }_2^2 \geq \sigma^2 \paren{ \frac{1}{2 \paren{\ln(n)}^2} - \frac{(\ComegaC+2\ComegaD+\ComegaE) \gamma \ln(n)}{n}}  \enspace . \]
Hence, the blow up of $\df(\lhminb(C))$ holding when the penalty is below the minimal penalty also implies a blow up of the risk $n^{-1} \norm{ \Fh_{\lhminb(C)} - F }_2^2\,$.
\end{remark}

\begin{remark}
Using assumption \eqref{hyp.Al} and Lemma~\ref{le.retire.condition.A3}, we deduce that for every $\lamm\,$, 
\[ \frac{\paren{\df(\lambda)}^2}{n} \leq \tr(\Al^{\top}\Al) \leq 2 \tr(\Al) = 2 \df(\lambda) \enspace . \]
Therefore, on the event defined by Proposition~\ref{A.pro.mini}, a jump in $\tr \paren{ A_{\lhminb(C)}^{\top} A_{\lhminb(C)} }$ also occurs when $C$ goes through $\sigma^2\,$. 
Indeed, 
\[ \df(\lhminb(C)) \geq \mindimavant \quad \mbox{implies that} \quad \tr \paren{ A_{\lhminb(C)}^{\top} A_{\lhminb(C)} } \geq \mindimavant^2 \]
and this lower bound becomes $n/9$ if $\mindimavant=n/3\,$.
Furthermore, 
\[ \df(\lhminb(C)) \leq \majdimapres \quad \mbox{implies that} \quad \tr \paren{ A_{\lhminb(C)}^{\top} A_{\lhminb(C)} } \leq 2 \majdimapres \]
which is equal to $2 n^{3/4}$ when $\majdimapres = n^{ 3/4}\,$, for instance.
\end{remark}

Let us now prove Proposition~\ref{A.pro.mini}. The proof is organized is as follows: 
\begin{enumerate}
\item Section~\ref{A.sec.thm.mini.start} makes use of the definition of the event $\Omega_{\gamma \ln(n)}(\Lambda,\Comega)$ for controlling uniformly over $C$ and $\lL$ the criterion $\crit_C$ minimized by $\lhminb(C)$. 
\item Section~\ref{A.sec.thm.mini.below} considers the case $C<\sigma^2$, showing that if $A_{\lidentite}=\Id$ and $\lL$ satisfies $\df(\lambda) \leq \mindimavant$, then $\crit_C(\lidentite) \leq \crit_C(\lambda)$, hence $\lhminb(C) \geq \mindimavant$. 
\item Section~\ref{A.sec.thm.mini.above} considers the case $C>\sigma^2$, showing that if $\lL$ satisfies $\df(\lambda) \geq \majdimapres$, then $\crit_C(\lhypbias) < \crit_C(\lambda)$, hence $\lhminb(C) \leq  \majdimapres$. 
\end{enumerate}

\subsection{General starting point} \label{A.sec.thm.mini.start}
Combining Eq.~\eqref{eq.lhmin} with Eq.~\eqref{A.eq.decomp.riskFhl} and \eqref{A.eq.decomp.riskempFhl}, for every $C \geq 0\,$,
$\lhminb(C)$ also minimizes over $\lamm$
\begin{align*}
\crit_C(\lambda) &\egaldef \norm{\Fhl - Y}_2^2 - \norm{\varepsilon}_2^2 + C \paren{ 2 \tr(\Al) - \tr(\Al^\top \Al)} \\
&= b(\lambda) + (\sigma^{-2} C-1) \paren{ 2 \vun(\lambda) - \vdeux(\lambda)} - 2 \delun(\lambda) + \deldeux(\lambda) + \delta_3(\lambda) + \delta_4(\lambda) \enspace .
\end{align*}

On the event $\Omega_{\gamma \ln(n)}(\Lambda,\Comega)$, using Eq.~\eqref{A.eq.conc.del1}, \eqref{A.eq.conc.del2}, \eqref{A.eq.conc.del3} and~\eqref{A.eq.conc.del4}, we get for every $\lamm$
%
 %
 \begin{align*}
 \crit_C(\lambda) &\geq ( 1 - \theta_3 - \theta_4)b(\lambda) + (\sigma^{-2}C-1) (2 \vun(\lambda) - \vdeux(\lambda)) 
         - (2 \theta_1 + \theta_2) \vdeux(\lambda) \\
 & \quad - \paren{ 2 \ComegaA + 2 \ComegaB \theta_1^{-1} + \ComegaC + \ComegaD \theta_2^{-1} 
 		+ \ComegaE \theta_3^{-1} + \ComegaF \theta_4^{-1} } \gamma \ln(n) \sigma^2
\enspace ,
 \end{align*}
 and
 \begin{align*}
 \crit_C(\lambda) &\leq ( 1 + \theta_3 + \theta_4)b(\lambda) + (\sigma^{-2}C-1) (2 \vun(\lambda) - \vdeux(\lambda)) 
         + (2 \theta_1 + \theta_2) \vdeux(\lambda) \\
 & \quad + \paren{ 2 \ComegaA + 2 \ComegaB \theta_1^{-1} + \ComegaC + \ComegaD \theta_2^{-1} 
 		+ \ComegaE \theta_3^{-1} + \ComegaF \theta_4^{-1} } \gamma \ln(n) \sigma^2
\enspace .
 \end{align*}
 Taking $\theta_1 = \theta_2 = \theta/3 $ for some $\theta \in (0,3]\,$, and $\theta_3 = \theta_4 = 1/3\,$, this implies:
\begin{align} 
\label{A.eq.penmin.minor-crit}
\crit_C(\lambda) &\geq \frac{b(\lambda)} {3} + (\sigma^{-2} C-1) (2 \vun(\lambda) - \vdeux(\lambda)) - \theta \vdeux(\lambda) - \paren{ \KthetaA + \KthetaB\theta^{-1}} \gamma \ln(n) \sigma^2 \\
\label{A.eq.penmin.major-crit}
\crit_C(\lambda) &\leq \frac{5b(\lambda)} {3} + (\sigma^{-2} C-1) (2 \vun(\lambda) - \vdeux(\lambda)) + \theta \vdeux(\lambda) + \paren{ \KthetaA + \KthetaB\theta^{-1}} \gamma \ln(n) \sigma^2
\enspace .  		
\end{align}

\subsection{Below the minimal penalty} \label{A.sec.thm.mini.below}
We assume in this subsection  that $C \in [0,\sigma^2)$ and $\mindimavant \in [ 0, n/2 ) \,$. 
Let $\lamm\,$. Three cases can be distinguished:
\begin{enumerate}
\item If $\KthetaB \gamma \ln(n) / 2 \leq \dfl \leq \mindimavant\,$, then Eq.~\eqref{A.eq.penmin.minor-crit} yields 
\begin{align} \notag 
\crit_C(\lambda) 
&\geq 2 (C- \sigma^2) \dfl  -  2 \theta \dfl \sigma^2 - \theta^{-1} \KthetaB \gamma \ln(n) \sigma^2 - \KthetaA \gamma \ln(n) \sigma^2  
\detail{   \\ \detailtag   
 &\geq 2 (C- \sigma^2) \dfl  -  2 \sqrt{2 \KthetaB \gamma \ln(n) \dfl} \sigma^2 - \KthetaA \gamma \ln(n) \sigma^2    }
\\ 
&\geq 
2 (C- \sigma^2) \mindimavant  -  2 \sqrt{2 \KthetaB \gamma \ln(n) \mindimavant} \sigma^2 - \KthetaA \gamma \ln(n) \sigma^2 
\enspace , 
\label{A.eq.penmin.minor-crit.below.2} 
\end{align}
by taking $\theta = \sqrt{\KthetaB \gamma \ln(n) / (2 \dfl)} \leq 1\,$.
\item If $\dfl \leq \KthetaB \gamma \ln(n) / 2 \,$, taking $\theta=1$ in Eq.~\eqref{A.eq.penmin.minor-crit} yields 
\begin{align} 
\notag 
\crit_C(\lambda) 
&\geq 
2 (C- \sigma^2) \dfl  -  2 \sigma^2 \dfl - (\KthetaA + \KthetaB) \gamma \ln(n) \sigma^2  
\detail{   \\ \detailtag   
   &\geq  - 4 \dfl \sigma^2 - \paren{  \KthetaB + \KthetaA}  \gamma \ln(n) \sigma^2   }
\\ 
&\geq 
- \paren{ 3 \KthetaB + \KthetaA}  \gamma \ln(n) \sigma^2
\label{A.eq.penmin.minor-crit.below.2b} 
\enspace .
\end{align}
\item For $\lidentite  \in \Lambda$ such that $A_{\lidentite} = \Id$, $\df(\lidentite) = n$ and $b(\lidentite)=0$. So, Eq.~\eqref{A.eq.penmin.major-crit} yields 
\begin{align} \notag
\crit_C(\lidentite) 
&\leq (C- \sigma^2 + \theta \sigma^2) n + \paren{ \KthetaA + \KthetaB\theta^{-1}} \gamma \ln(n) \sigma^2
\\ 
&=  (C- \sigma^2) n + 2 \sqrt{n \KthetaB \gamma \ln(n) } \sigma^2 + \KthetaA \gamma \ln(n) \sigma^2
 \enspace , 
\label{A.eq.penmin.major-crit.below.3} 
\end{align}
by taking $\theta = \sqrt{\KthetaB \gamma \ln(n) / n} \leq 1\,$, assuming $n / \ln(n) \geq \KthetaB \gamma \,$.
\end{enumerate}

\paragraph{First condition on $C$: case 1 vs. case 3.}
Comparing Eq.~\eqref{A.eq.penmin.minor-crit.below.2}  and Eq.~\eqref{A.eq.penmin.major-crit.below.3}, we get that 
\begin{equation} \label{A.eq.penmin.comp-crit.below.1} 
\crit_C(\lidentite) < \inf_{\lamm, \, \KthetaB \gamma \ln(n) / 2 \leq \dfl \leq \mindimavant } \set{ 
\crit_C(\lambda) } 
\end{equation}
if 
\detail{
\begin{align*} \detailtag
& \qquad 2 (C- \sigma^2) \mindimavant  -  2 \sqrt{2 \KthetaB \gamma \ln(n) \mindimavant} \sigma^2 - \KthetaA \gamma \ln(n) \sigma^2 
\\ \detailtag &> (C- \sigma^2) n  + 2 \sqrt{n \KthetaB \gamma \ln(n) } \sigma^2 + \KthetaA \gamma \ln(n) \sigma^2
\end{align*}
which holds if 
}
\[ 
 (1 - \sigma^{-2} C ) (n - 2 \mindimavant )
>
2 \paren{ \sqrt{n} + \sqrt{2  \mindimavant}} \sqrt{\KthetaB \gamma \ln(n)  }  + 2 \KthetaA \gamma \ln(n) 
\enspace .
\]
The right-hand side of the above inequality is smaller than 
\begin{align*}
\gamma \sqrt{n \ln(n)} \croch{ 2 \sqrt{2 \KthetaB} + \frac{2 \KthetaA}{\ctepenminA} }
\end{align*}
where we used that $ \mindimavant < n/2$ and $\gamma \geq 2\,$, 
and we assumed that $\sqrt{n / \ln(n)} \geq \ctepenminA$ for some $\ctepenminA>0$ to be chosen later.
Hence, Eq.~\eqref{A.eq.penmin.comp-crit.below.1} holds as soon as $\sqrt{n / \ln(n)} \geq \ctepenminA$ and 
\begin{equation} \label{A.eq.penmin.cond-C.below.1}
1 - \sigma^{-2} C 
>
2 \croch{ \sqrt{2 \KthetaB} + \frac{\KthetaA}{ \ctepenminA }  } \frac{\gamma \sqrt{n \ln(n)}} {n - 2 \mindimavant }
\enspace .
\end{equation}

\paragraph{Second condition on $C$: case 2 vs. case 3.}
Comparing Eq.~\eqref{A.eq.penmin.minor-crit.below.2b}  and Eq.~\eqref{A.eq.penmin.major-crit.below.3}, we get that 
\begin{equation} \label{A.eq.penmin.comp-crit.below.2} 
 \crit_C(\lidentite) < \inf_{\lamm, \, \dfl \leq \KthetaB \gamma \ln(n) / 2 } \set{ 
\crit_C(\lambda) }  
\end{equation}
if 
\detail{
\begin{align*}
\detailtag
- \paren{ 3 \KthetaB + \KthetaA}  \gamma \ln(n) \sigma^2
>   (C- \sigma^2) n  + 2 \sqrt{n \KthetaB \gamma \ln(n) } \sigma^2 + \KthetaA \gamma \ln(n) \sigma^2
\end{align*}
which holds if 
}
\[ 
(1 - \sigma^{-2} C ) n
>
2 \sqrt{n \KthetaB \gamma \ln(n)  }  + \paren{ 3 \KthetaB + 2 \KthetaA} \gamma \ln(n) 
\enspace .
\]
The right-hand side of the above inequality is smaller than 
\begin{align*}
\paren{ \sqrt{2 \KthetaB} + \frac{ 3 \KthetaB + 2 \KthetaA }{\ctepenminA} } \gamma \sqrt{n \ln(n)}
\end{align*}
where we used that $\gamma \geq 2$ and we assumed that $\sqrt{n / \ln(n)} \geq \ctepenminA$.
Hence, Eq.~\eqref{A.eq.penmin.comp-crit.below.2} holds as soon as $\sqrt{n / \ln(n)} \geq \ctepenminA$ and 
\begin{equation} \label{A.eq.penmin.cond-C.below.2}
1 - \sigma^{-2} C 
>
\paren{ \sqrt{2 \KthetaB} + \frac{ 3 \KthetaB + 2 \KthetaA }{\ctepenminA} } \gamma \sqrt{\frac{\ln(n)}{n}}
\enspace .
\end{equation}

\paragraph{Combining the two conditions.}
Finally, we have proved that $\df(\lhminb(C)) > \mindimavant $ if $\sqrt{n / \ln(n)} \geq \ctepenminA$ and 
conditions \eqref{A.eq.penmin.cond-C.below.1} and \eqref{A.eq.penmin.cond-C.below.2} are both satisfied, 
hence if 
\[ 1 - \sigma^{-2} C  >  
\paren{ 2  \sqrt{2 \KthetaB} + \frac{1}{\ctepenminB} } \frac{\gamma \sqrt{ n \ln(n) } }{ n - 2 \mindimavant }   \]
and 
$\sqrt{n / \ln(n)} \geq \ctepenminA = \paren{ 2 \KthetaA + 3 \KthetaB } \ctepenminB \,$.

\subsection{Above the minimal penalty} \label{A.sec.thm.mini.above}
We assume in this subsection  that $C > \sigma^2\,$. 
Let $\lamm\,$. 
As in Section~\ref{A.sec.thm.mini.below}, we consider three cases.
\begin{enumerate}
\item If $\KthetaB \gamma \ln(n) / 2  \leq \dfl \leq \majdimbiais\,$, then, using assumption \eqref{hyp.Al} in Eq.~\eqref{A.eq.penmin.major-crit} yields 
\begin{align} \notag
\crit_C(\lambda) 
&\leq \frac{5b(\lambda)} {3} + (C- \sigma^2) (2 \dfl -\tr(\Al^{\top} \Al))  + \theta \tr(\Al^{\top} \Al) \sigma^2+ \paren{ \KthetaA + \KthetaB\theta^{-1}} \gamma \ln(n) \sigma^2
\\ \notag
&\leq 
\frac{5b(\lambda)} {3} + 2 (C- \sigma^2) \dfl + 2 \theta \dfl \sigma^2 + \paren{ \KthetaA + \KthetaB\theta^{-1}} \gamma \ln(n) \sigma^2
\\ \notag 
&= 
\frac{5b(\lambda)} {3} + 2 (C- \sigma^2) \dfl + 2 \sqrt{2 \dfl \KthetaB \gamma \ln(n) } \sigma^2 + \KthetaA \gamma \ln(n) \sigma^2
\\
&\leq 
\frac{5 b(\lambda)}{3} + 2 (C- \sigma^2) \majdimbiais  + 2 \sqrt{2 \majdimbiais \KthetaB \gamma \ln(n) } \sigma^2 + \KthetaA \gamma \ln(n) \sigma^2 \enspace , 
\label{A.eq.penmin.major-crit.above.2} 
\end{align}
by taking $\theta = \sqrt{\KthetaB \gamma \ln(n) / (2 \dfl)} \leq 1\,$.
\item If $\dfl \leq \KthetaB \gamma \ln(n) / 2\,$, then, using assumption \eqref{hyp.Al} in Eq.~\eqref{A.eq.penmin.major-crit} with $\theta=1$ yields 
\begin{align} \notag
\crit_C(\lambda) 
\detail{
&\leq 
\frac{5b(\lambda)} {3} + 2 (C- \sigma^2) \dfl + (2 \sigma^2 - C) \tr(\Al^{\top} \Al) + \paren{ \KthetaA + \KthetaB} \gamma \ln(n) \sigma^2 \detailtag
\\ 
}
&\leq 
\frac{5b(\lambda)} {3} + 2 (C- \sigma^2) \dfl + 2 \dfl \sigma^2 + \paren{ \KthetaA + \KthetaB} \gamma \ln(n) \sigma^2
\notag 
\\ 
&\leq 
\frac{5b(\lambda)} {3} + (C- \sigma^2) \KthetaB \gamma \ln(n) + \paren{ \KthetaA + 2 \KthetaB} \gamma \ln(n) \sigma^2
\label{A.eq.penmin.major-crit.above.2b} 
\enspace .
\end{align}
\item If $\dfl \geq \majdimapres\,$, then, using assumption \eqref{hyp.Al}  in Eq.~\eqref{A.eq.penmin.minor-crit} yields 
\begin{align} \notag 
\crit_C(\lambda) 
&\geq (C-\sigma^2) (2 \dfl - \tr(\Al^{\top}\Al)) -  \theta \tr(\Al^{\top}\Al) \sigma^2 - \paren{ \KthetaA + \theta^{-1} \KthetaB } \gamma \ln(n) \sigma^2 
\\ \notag 
&\geq \cAl (C-\sigma^2) \dfl  - 2 \theta \dfl \sigma^2 - \paren{ \KthetaA + \theta^{-1} \KthetaB } \gamma \ln(n) \sigma^2 
\\ \notag 
&= 
\cAl (C-\sigma^2) \dfl  -  2 \sqrt{2 \KthetaB \gamma \ln(n) \dfl} \sigma^2 - \KthetaA \gamma \ln(n) \sigma^2 
\\
&\geq 
\cAl (C-\sigma^2) \majdimapres  -  2 \sqrt{2 \KthetaB \gamma \ln(n) n} \sigma^2 - \KthetaA \gamma \ln(n) \sigma^2 
\enspace , 
\label{A.eq.penmin.minor-crit.above.2} 
\end{align}
by taking $\theta = \sqrt{\KthetaB \gamma \ln(n) / (2 \dfl)} \leq 1$ since $\majdimapres \geq \KthetaB \gamma \ln(n) /2\,$, and using that \eqref{hyp.Al} implies $\df(\lambda) \leq n\,$.
\end{enumerate}

Eq.~\eqref{A.eq.penmin.minor-crit.above.2} implies 
\begin{equation} \label{A.eq.penmin.minor-crit.above.3} 
\inf_{\lamm, \, \dfl \geq \majdimapres} \set{ \crit_C (\lambda) } \geq 
\cAl (C-\sigma^2)\majdimapres  -  2 \sqrt{2 \KthetaB \gamma \ln(n) n} \sigma^2 - \KthetaA \gamma \ln(n) \sigma^2 
\enspace . \end{equation}

Let $\lambda = \lhypbias$ given by assumption \eqref{A.hyp.thm.mini.2.bis}. Two cases can occur:

\paragraph{If $\lhypbias$ matches case~1:}
If $\majdimbiais \geq \df(\lhypbias)\geq \KthetaB \gamma \ln(n)/2\,$, 
taking $\lambda = \lhypbias$ in Eq.~\eqref{A.eq.penmin.major-crit.above.2} implies
\begin{equation} \label{A.eq.penmin.major-crit.above.3} 
\crit_C(\lhypbias) \leq 
\frac{5 \sigma^2 \majbiais}{3} + 2 (C- \sigma^2) \majdimbiais  + 2 \sqrt{2 \majdimbiais \KthetaB \gamma \ln(n) } \sigma^2 + \KthetaA \gamma \ln(n) \sigma^2 \enspace . 
\end{equation}

Comparing Eq.~\eqref{A.eq.penmin.major-crit.above.3} and Eq.~\eqref{A.eq.penmin.minor-crit.above.3}, we get that 
\begin{equation} \label{A.eq.penmin.comp-crit.above.1} 
 \crit_C(\lhypbias) < \inf_{\lamm, \, \dfl \geq \majdimapres} \set{ \crit_C(\lambda) }  
\end{equation}
hence $\df(\lhminb(C)) < \majdimapres$ if 
\[ 
(C-\sigma^2) \paren{ \cAl  \majdimapres - 2 \majdimbiais }
>
 \frac{5 \sigma^2 \majbiais}{3} 
 +  2  \paren{ n^{1/2} + \sqrt{\majdimbiais}} \sqrt{ 2 \KthetaB \gamma \ln(n)} \sigma^2 
 + 2\KthetaA \gamma \ln(n) \sigma^2  
 \enspace ,
\]
which holds if
\[
(\sigma^{-2} C- 1) \paren{ \cAl \majdimapres - 2 \majdimbiais }
>
 \frac{5 \majbiais}{3} 
 +  4  \sqrt{ 2 \KthetaB \gamma n \ln(n)} 
 + 2\KthetaA \gamma \ln(n) 
\]
since $\majdimbiais \leq n$.
Using in addition that $\gamma \geq 2$ and $\majbiais \geq \sqrt{n \ln(n)}\,$, the right-hand side of the above equation is smaller than 
\begin{align*}
\frac{5 \majbiais}{3} 
 +  4 \sqrt{\KthetaB} \gamma \sqrt{ n \ln(n)} 
 + 2 \KthetaA \gamma \ln(n) 
&\leq 
\detail{ \detailtag
\frac{5 \majbiais}{3} 
 +  \paren{ 4 \sqrt{\KthetaB}  + 2 \KthetaA \ctepenminC^{-1}} \gamma \sqrt{ n \ln(n)} 
\\ &\leq 
}
\croch{ \frac{5}{6}  +  4 \sqrt{\KthetaB}  + 2 \KthetaA \ctepenminC^{-1} } \gamma  \majbiais \enspace , 
\end{align*}
assuming that $\sqrt{n / \ln(n)} \geq \ctepenminC$ for some $\ctepenminC >0$ to be chosen later.
Now, since $\majdimapres \geq 20 \majdimbiais / \cAl\,$,  
$\cAl  \majdimapres - 2 \majdimbiais  \geq 9 \cAl  \majdimapres / 10 $ and Eq.~\eqref{A.eq.penmin.comp-crit.above.1}  holds as soon as 
\begin{equation} \label{A.eq.penmin.cond.above.1} 
\sigma^{-2} C- 1 
> 
 \frac{5}{9} \croch{ \frac{5}{3}  +  8 \sqrt{\KthetaB}  + \frac{4 \KthetaA}{ \ctepenminC } } \frac{\gamma  \majbiais}{\cAl \majdimapres} 
\enspace . 
\end{equation}


\paragraph{If $\lhypbias$ matches case~2:}
If $\df(\lhypbias)\leq \KthetaB \gamma \ln(n)/2\,$, taking $\lambda = \lhypbias$ in Eq.~\eqref{A.eq.penmin.major-crit.above.2b} implies 
\begin{equation} \label{A.eq.penmin.major-crit.above.2b-prime}
\crit_C(\lhypbias) \leq \frac{5 \majbiais \sigma^2} {3} + (C- \sigma^2) \KthetaB \gamma \ln(n) + \paren{ \KthetaA + 2 \KthetaB} \gamma \ln(n) \sigma^2 \enspace .
\end{equation}
Comparing Eq.~\eqref{A.eq.penmin.major-crit.above.2b-prime} and Eq.~\eqref{A.eq.penmin.minor-crit.above.3}, we get that 
\begin{equation} \label{A.eq.penmin.comp-crit.above.2} 
 \crit_C(\lhypbias) < \inf_{\lamm, \, \dfl \geq \majdimapres} \set{ \crit_C(\lambda) }  
\end{equation}
hence $\df(\lhminb(C)) < \majdimapres$ if 
\detail{
\begin{equation} \detailtag \begin{split}
\frac{5 \majbiais \sigma^2} {3} + (C- \sigma^2) \KthetaB \gamma \ln(n) + \paren{ \KthetaA + 2 \KthetaB} \gamma \ln(n) \sigma^2
\\ <
\cAl (C-\sigma^2) \majdimapres  -  2 \sqrt{2 \KthetaB \gamma \ln(n) n} \sigma^2 - \KthetaA \gamma \ln(n) \sigma^2 
 \enspace ,
 \end{split}
\end{equation}
}
which holds if 
\[ 
(\sigma^{-2} C- 1) \paren{ \cAl \majdimapres - \KthetaB \gamma \ln(n)}  
>
\frac{5 \majbiais} {3}  + \paren{ 2 \KthetaA + 2 \KthetaB} \gamma \ln(n) 
+ 2 \sqrt{2 \KthetaB \gamma \ln(n) n} 
 \enspace .
\]
The right-hand side of the above equation is smaller than 
\begin{align*}
\croch{ \frac{5 } {6} + 2 \sqrt{\KthetaB } + \frac{ 2 \paren{ \KthetaA +  \KthetaB}}{\ctepenminC}
 } \gamma \majbiais 
\end{align*}
since $\gamma \geq 2\,$, $\majbiais \geq \sqrt{n \ln(n)}$ and $\sqrt{n/\ln(n)} \geq \ctepenminC>0\,$.
Now, $\cAl \majdimapres - \KthetaB \gamma \ln(n) \geq 9 \cAl \majdimapres / 10$ since $ \majdimapres \geq 10 \KthetaB \gamma \ln(n) / \cAl\,$, so that 
Eq.~\eqref{A.eq.penmin.comp-crit.above.2}  holds as soon as 
\begin{equation} \label{A.eq.penmin.cond.above.2} 
\sigma^{-2} C- 1
>
\croch{ \frac{5 } {3} + 4 \sqrt{\KthetaB } + \frac{ 4 (\KthetaA + \KthetaB)}{\ctepenminC}  } \frac{5}{9} \frac {\gamma \majbiais} {\cAl \majdimapres}
 \enspace .
\end{equation}

\paragraph{Combining the two conditions.}
Finally, we have proved that whatever the value of $\df(\lhypbias)\,$, $\df(\lhminb(C)) < \majdimapres $ holds if conditions~\eqref{A.eq.penmin.cond.above.1} and~\eqref{A.eq.penmin.cond.above.2} are both satisfied and if 
$\sqrt{n/\ln(n)} \geq \ctepenminC \,$, $\majdimapres \geq \frac{10}{\cAl} \max\sset{ \KthetaB \gamma \ln(n) , 2 \majdimbiais} \,$.
Eq.~\eqref{A.eq.thm.lhmin.above.bis} follows by choosing $\ctepenminC = 3 (\KthetaA + \KthetaB) \ctepenminB \,$.
Merging all assumptions on $n$ made in the proof, the constant $\nminProMini$ is defined by 
\[ \inf_{n \geq n_3} \set{ \frac{n}{\ln(n)} } \geq \max\set{ \KthetaB \gamma , 16 \paren{ \KthetaA + \KthetaB }^2 \ctepenminB^2 } \enspace . \]

\subsection{Second part of Proposition~\ref{A.pro.mini}}
We apply Lemmas~\ref{le.Omega.finite} and~\ref{le.Omega.cont-ridge} and the union bound, taking $ \ctepenminB = 3 $ and noticing that 
\begin{align*}
\KthetaA \detail{ \detailtag &= 4 \majnormAl + 2 \majnormAl^2 + 3 \paren{ \max\sset{ 2 \majnormAl^2 , \frac{\cteLinRidge^2}{4} } + \frac{\cteLinRidge^2}{4} }  \\ }
&\leq 6 \majnormAl^2 + \max\sset{ 6 \majnormAl^2 , \cteG } + \cteG 
\quad \mbox{and} \quad
\KthetaB = 6 + 3 \majnormAl^2 \leq \cteAa \majnormAl^2 \enspace . \qed
\end{align*}

\section{Oracle inequality (proof of Theorem~\ref{thm.oracle.lhopt})} \label{A.sec.thm.oracle}
Recall the definition \eqref{def.lhopt} of $\lhopt(C)\,$:
\begin{equation*} 
\forall C \geq 0 \, , \qquad \lhopt(C) \in \arg\min_{\lamm} \set{ \norm{\Fhl - Y}_2^2 + 2 C \tr(\Al) } \enspace . 
\end{equation*} 
We prove in this section Theorem~\ref{thm.oracle.lhopt}, as a corollary of the following proposition.
\begin{proposition} \label{A.pro.oracle}
Let $\lhopt(C)$ be defined by Eq.~\eqref{def.lhopt}.
Assume that \eqref{hyp.Al} holds true. 
Let $\Comega \in [0,+\infty)^6$ and define 
\[ 
\KthetaC \egaldef 4  \ComegaA
\qquad 
\KthetaD \egaldef 8 \paren{  \ComegaC + 2 \ComegaD + 2 \ComegaE } 
\quad \mbox{and} \quad 
\KthetaE \egaldef 2 \paren{4 \ComegaB + \ComegaF}
 \enspace . \]

Then, for every $\gamma \geq 2\,$, 
on the event $\Omega_{\gamma \ln(n)}(\Lambda,\Comega)\,$, 
for every $\theta \in (0,1/4)\,$,
\begin{equation} \label{A.eq.thm.lhopt}
\begin{split}
n^{-1} \norm{\Fh_{\lhopt(C)} - F}_2^2 
\leq \frac{1 + 2 \theta}{1 -  4\theta } \inf_{\lamm} \set{ n^{-1} \norm{\Fh_{\lambda} - F}_2^2 + \frac{2 (C-\sigma^2)_+ \tr(\Al)}{n} } 
\\
+ \frac{\paren{C \sigma^{-2} - 1}^2 \un_{C \leq \sigma^2} \sigma^2}{\theta \paren{ 1-4\theta }} 
+ \frac{ \paren{ \KthetaC  + \frac{3}{4} \KthetaD \theta  +  \KthetaE \theta^{-1}} }{ 1 - 4\theta } \frac{\ln(n) \gamma \sigma^2}{n}
\enspace .
\end{split} 
\end{equation}
and
\begin{equation} \label{A.eq.thm.lhopt.bis}
\begin{split}
n^{-1} \norm{\Fh_{\lhopt(C)} - F}_2^2 
\leq \frac{1 + 4 \theta}{1 -  4\theta } \inf_{\lamm} \set{ n^{-1} \norm{\Fh_{\lambda} - F}_2^2  } 
+ \frac{\paren{C \sigma^{-2} - 1}^2 \sigma^2}{\theta \paren{ 1-4\theta }} 
\\
+ \frac{ \paren{ \KthetaC  + \KthetaD \theta  +  \KthetaE \theta^{-1}} }{ 1 - 4\theta } \frac{\ln(n) \gamma \sigma^2}{n}
\enspace .
\end{split} 
\end{equation}
\end{proposition}

\begin{corollary}  \label{A.cor.oracle}
%
%
Under assumptions \eqref{hyp.Al}, \eqref{hyp.Lam} and \eqref{hyp.eps-Gauss-hom}, 
for every $\gamma \geq 2\,$, 
with probability at least $1 - \paren{ 6 \card(\Lambda_0) + \NLamfusrid \exp\paren{\betaLRUd + \ln(n)}} n^{-\gamma}\,$, 
for every $C>0$, for every $\theta \in (0,1/4)\,$, Eq.~\eqref{A.eq.thm.lhopt} and~\eqref{A.eq.thm.lhopt.bis} hold true, with $(\KthetaC,\KthetaD,\KthetaE)$ replaced by $(\betaThmOracleA,\betaThmOracleB,\betaThmOracleC)$ where 
\[ 
\betaThmOracleA = 8 \majnormAl
 \qquad 
\betaThmOracleB = 8 \paren{ 4 \majnormAl^2 + \max\set{ 4 \majnormAl^2 , \cteDa } } 
 \quad \mbox{and} \quad  
\betaThmOracleC 
= \cteD
 \enspace .\] 
\end{corollary}

We can now prove Theorem~\ref{thm.oracle.lhopt}:
\begin{proof}[Proof of Theorem~\ref{thm.oracle.lhopt}]
The proof is similar to the one of Theorem~\ref{thm.mini}, starting from Corollary~\ref{A.cor.oracle} instead of Proposition~\ref{A.pro.mini}.
In addition, we remark that $\forall \theta \in (0,1/8)\,$, 
\[ \frac{1+2\theta}{1-4\theta} \leq \frac{1+4\theta}{1-4\theta} \leq 1 + 16 \theta \] and we take $\theta = \eta / 16\,$. 
Choosing $\gamma = \alLamMix + \delta \geq 2\,$ in Corollary~\ref{A.cor.oracle}, Eq.~\eqref{A.eq.thm.lhopt} and~\eqref{A.eq.thm.lhopt.bis} become Eq.~\eqref{eq.thm.lhopt} and~\eqref{eq.thm.lhopt-alt}, 
which holds with probability at least $1 - \paren{ 6 \CLamfusrid  + \CLamfusdis} n^{ - \delta}\,$, assuming that $n$ is larger than some numerical constant $\nminThmOracle = \exp(\betaLRUd)\,$. 
\end{proof}

Let us now prove Proposition~\ref{A.pro.oracle} and Corollary~\ref{A.cor.oracle}.
The proof is organized is as follows: 
\begin{enumerate}
\item In Section~\ref{A.sec.thm.oracle.start}, standard algebraic manipulations reduce the problem to bounding $\widehat{\Delta}(\lambda) \egaldef - 2 \delun(\lambda) + \delta_4(\lambda)$ uniformly over $\lL$. 
\item Section~\ref{A.sec.thm.oracle.conc} make use of the definition of the event $\Omega_{\gamma \ln(n)}(\Lambda,\Comega)$ for bounding $\widehat{\Delta}(\lambda)$.
\item In Section~\ref{A.sec.thm.oracle.small-terms}, remainder terms proportional to $\tr(\Al)$ are upper bounded in terms of $\tr(\Al^{\top} \Al)$, so they can be compared to the risk of $\Fhl$. 
\end{enumerate}

\subsection{General starting point} \label{A.sec.thm.oracle.start}
Combining Eq.~\eqref{A.eq.decomp.riskempFhl} and~\eqref{def.lhopt}, we obtain that for every $C>0$ and every $\lamm\,$, 
\begin{equation} \label{A.eq.debut-oracle}
\begin{split}
\norm{\Fh_{\lhopt(C)} - F}_2^2 + 2 (C-\sigma^2) \tr(A_{\lhopt(C)}) + \widehat{\Delta}(\lhopt(C)) 
\\ 
\leq \inf_{\lamm} \set{ \norm{\Fh_{\lambda} - F}_2^2 + 2 (C-\sigma^2) \tr(\Al) + \widehat{\Delta}(\lambda) } 
\enspace .
\end{split} 
\end{equation}
where 
\[ \forall \lamm \, , \quad \widehat{\Delta}(\lambda) \egaldef - 2 \delun(\lambda) + \delta_4(\lambda) \enspace . \]
Inequality \eqref{A.eq.debut-oracle} implies an oracle inequality as soon as $\widehat{\Delta}(\lambda)$ is small compared to $\snorm{\Fh_{\lambda} - F}_2^2$ and $C-\sigma^2$ is small enough.

\subsection{With concentration inequalities} \label{A.sec.thm.oracle.conc}

On the event $\Omega_{\gamma \ln(n)}$, using Eq.~\eqref{A.eq.conc.del1} and~\eqref{A.eq.conc.del4} with $2\theta_1=\theta_4=\theta \in(0,1]\,$, we get that for every $\theta \in (0,1]$ and $\lamm\,$, 
\begin{equation} 
\label{A.eq.debut-oracle.conc1}
\absj{\widehat{\Delta}(\lambda)} \leq 
\theta \croch{ b(\lambda) + \vdeux(\lambda)}
+ \paren{\frac{\KthetaC}{2} + \frac{\KthetaE}{2\theta} } \gamma \ln(n) \sigma^2
 \end{equation}
Using also Eq.~\eqref{A.eq.decomp.riskFhl}, \eqref{A.eq.conc.del2} and~\eqref{A.eq.conc.del3} with $\theta_2=\theta_3=1/2\,$, we get that for every $\lamm\,$, 
\begin{align} 
\norm{\Fhl - F}^2 
\detail{ \detailtag
&= b(\lambda) + \vdeux(\lambda) + \delta_2(\lambda) + \delta_3(\lambda) 
        \\  
}
&\geq 
\frac{1}{2} \paren{ b(\lambda) + \vdeux(\lambda)}
- \frac{\KthetaD}{8}  \gamma \ln(n) \sigma^2 
\label{A.eq.debut-oracle.conc2a}
 \end{align}
so that 
\begin{align} 
b(\lambda) + \vdeux(\lambda)
\leq
2 \norm{\Fhl - F}^2 + \frac{\KthetaD}{4} \gamma \ln(n) \sigma^2 
\enspace . 
\label{A.eq.debut-oracle.conc2}
 \end{align}

Combining Eq.~\eqref{A.eq.debut-oracle.conc1} and~\eqref{A.eq.debut-oracle.conc2}, we get 
\begin{align*}
\absj{\widehat{\Delta}(\lambda)} 
&\leq 
2 \theta \norm{\Fhl - F}^2 
+ \paren{\frac{\KthetaC}{2} + \frac{\theta \KthetaD}{4}  + \frac{\KthetaE}{2\theta} } \gamma \ln(n) \sigma^2
\end{align*}
so that Eq.~\eqref{A.eq.debut-oracle} implies that  
\begin{equation} \label{A.eq.debut-oracle.2}
\begin{split}
&\qquad (1 - 2\theta) \norm{\Fh_{\lhopt(C)} - F}_2^2 + 2 (C-\sigma^2) \tr(A_{\lhopt(C)}) 
\\
&\leq \inf_{\lamm} \set{ (1 + 2 \theta) \norm{\Fh_{\lambda} - F}_2^2 + 2 (C-\sigma^2) \tr(\Al) } 
+ \paren{ \KthetaC  + \frac{\theta \KthetaD}{2}  +  \frac{\KthetaE}{\theta} } \gamma \ln(n) \sigma^2
\enspace .
\end{split} 
\end{equation}

\subsection{Handling the small terms proportional to $\tr(\Al)$} \label{A.sec.thm.oracle.small-terms}
We will now make use of Lemma~\ref{le.retire.condition.A3} for handling the terms $(C-\sigma^2) \tr(\Al)$ that appear in Eq.~\eqref{A.eq.debut-oracle.2}. 
By Eq.~\eqref{eq.le.retire.condition.A3.v1} with $x=2 \absj{ C\sigma^{-2}-1 } \un_{C \leq \sigma^2}$ and Eq.~\eqref{A.eq.debut-oracle.conc2}, for every $\lL$, 
\begin{align} \notag 
2 \absj{C-\sigma^2} \tr(\Al) \un_{C \leq \sigma^{2}} 
&\leq \theta \vdeux \paren{\lambda} + \frac{n \paren{ C\sigma^{-2}-1 }^2 \un_{C \leq \sigma^{2}} \sigma^2}{\theta}  
\\
&\leq 2 \theta \norm{\Fhl - F}_2^2 + \frac{ \theta \KthetaD }{4}   \gamma \ln(n) \sigma^2 
+ \frac{n \paren{ C\sigma^{-2}-1 }^2 \sigma^2}{\theta}  \un_{C \leq \sigma^{2}} \enspace  . \label{A.eq.gestion-tr(Al)}
\end{align}
Applying Eq.~\eqref{A.eq.gestion-tr(Al)} with $\lambda=\lhopt(C)$, Eq.~\eqref{A.eq.debut-oracle.2} implies 
\begin{equation} \label{A.eq.debut-oracle.3}
\begin{split}
(1 - 4\theta) \norm{\Fh_{\lhopt(C)} - F}_2^2 
\leq \inf_{\lamm} \set{ (1 + 2 \theta) \norm{\Fh_{\lambda} - F}_2^2 + 2 (C-\sigma^2) \tr(\Al) } 
\\
+ \frac{n \paren{ C\sigma^{-2}-1 }^2 \sigma^2}{\theta}  \un_{C \leq \sigma^{2}}
+ \paren{ \KthetaC  + \frac{3 \theta \KthetaD}{4}  +  \frac{\KthetaE}{\theta} } \gamma \ln(n) \sigma^2
\enspace ,
\end{split} 
\end{equation}
which proves Eq.~\eqref{A.eq.thm.lhopt} since $\theta \in (0,1/4)\,$.
Applying again  Eq.~\eqref{A.eq.debut-oracle.conc2} and Eq.~\eqref{eq.le.retire.condition.A3.v1} with $x=2 \absj{ C\sigma^{-2}-1 } \un_{C \geq \sigma^2}$, Eq.~\eqref{A.eq.debut-oracle.3} implies
\begin{equation} \notag 
\begin{split}
(1 - 4\theta) \norm{\Fh_{\lhopt(C)} - F}_2^2 
\leq \inf_{\lamm} \set{ (1 + 4 \theta) \norm{\Fh_{\lambda} - F}_2^2  } 
\\
+ \frac{n \paren{ C\sigma^{-2}-1 }^2 \sigma^2}{\theta}  
+ \paren{ \KthetaC  + \theta \KthetaD  +  \frac{\KthetaE}{\theta} } \gamma \ln(n) \sigma^2
\enspace ,
\end{split} 
\end{equation}
which proves Eq.~\eqref{A.eq.thm.lhopt.bis} since $\theta \in (0,1/4)\,$. \qed

\subsection{Proof of Corollary~\ref{A.cor.oracle}}
The reasoning is the same as for proving the second part of Proposition~\ref{A.pro.mini}: we take $\Comega$ according to Lemmas~\ref{le.Omega.finite} and~\ref{le.Omega.cont-ridge}, so that the probability of $\Omega(\Lambda,\Comega)$ can be lower-bounded by the  union bound. \qed
%

\section{Proof of Theorem~\ref{thm.algo}} \label{A.sec.pr.thm.algo}
Theorem~\ref{thm.algo} is a straightforward consequence of Theorems~\ref{thm.mini} and~\ref{thm.oracle.lhopt}. 
Indeed, let us first remark the events defined by Theorems~\ref{thm.mini} and~\ref{thm.oracle.lhopt} are the same, namely $\Omega_{(\alLamMix+\delta) \ln(n)} (\Lambda,\Comega)$ for some well-chosen $\Comega\,$. 
So, on $\Omega_{(\alLamMix+\delta) \ln(n)} (\Lambda,\Comega)\,$, by Eq.~\eqref{A.eq.thm.lhmin.below} and~\eqref{A.eq.thm.lhmin.above}, 
\[ \absj{C \sigma^{-2} - 1 } \leq \max\set{ \betaThmMiniBelow \, , \, \betaThmMiniAbove }  (\alLamMix+\delta) \sqrt{\frac{\ln(n)}{n}}  \enspace . \]
Since Eq.~\eqref{eq.thm.lhopt-alt} also holds on $\Omega_{(\alLamMix+\delta) \ln(n)} (\Lambda,\Comega)\,$, 
we get Eq.~\eqref{eq.thm.oracle} with some numerical constant 
\[ \betaThmMixA \geq \majnormAl^{-2} \paren{ 32 \max\set{ \betaThmMiniBelow , \betaThmMiniAbove}^2 + \betaThmOracleD + \frac{\betaThmOracleE}{2} }  \enspace . \]
 
\medskip

We deduce from Eq.~\eqref{eq.thm.oracle} an oracle inequality in expectation by noting that 
if $n^{-1} \norm{\Fh_{\lh} - F}_2^2 \leq R_{n,\delta}$ on $\Omega_{(\alLamMix+\delta) \ln(n)}\,$, then 
\begin{align}
\E\croch{ \frac{1}{n} \norm{\Fh_{\lh} - F}_2^2 } 
&=  \E\croch{ \frac{\1_{\Omega_{(\alLamMix+\delta) \ln(n)}}}{n} \norm{\Fh_{\lh} - F}_2^2 } 
	+ \E\croch{ \frac{\1_{\Omega_{(\alLamMix+\delta) \ln(n)}^c}}{n} \norm{\Fh_{\lh} - F}_2^2 } 
\notag \\ 
\detail{
&\leq \E\croch{ R_{n,\delta} } 
	+  \frac{1}{n} \sqrt{ \Prob(\Omega_{(\alLamMix+\delta) \ln(n)}^c) } \sqrt{ \E\croch{ \norm{\Fh_{\lh} - F}_2^4 } } 
\detailtag \\ 
}
&\leq \E\croch{ R_{n,\delta} } 
	+  \frac{1}{n} \sqrt{ \CLamMix n^{-\delta} } \sqrt{ \E\croch{ \norm{\Fh_{\lh} - F}_2^4 } } 
\label{eq.thm.oracle.E.pr.1}
\end{align}
by Cauchy-Schwarz inequality. 
Now, remark that for every $\lamm\,$, 
\begin{align*} 
      \norm{\Fhl - F}_2^2 
\leq 2 \norm{\Al \varepsilon}^2 + 2 \norm{(I-\Al) F}^2 
\leq 2 \majnormAl^2 \norm{\varepsilon}^2 + 2 \paren{1 + \majnormAl}^2 \norm{F}^2
\end{align*}
where we used that $\normmat{\Al} \leq \majnormAl$ by assumption \eqref{hyp.Al}.
So,
\begin{align*}
 \E\croch{ \norm{\Fh_{\lh} - F}_2^4 } 
&\leq \E\croch{ \sup_{\lamm} \norm{\Fhl - F}_2^4 } 
\detail{
\\ \detailtag
&\leq 4 \E\croch{ \paren{ \majnormAl^2 \norm{\varepsilon}^2 + \paren{1 + \majnormAl}^2 \norm{F}^2 }^2 } \\ \detailtag
&= 4 \paren{ \majnormAl^4 \E\croch{ \norm{\varepsilon}^4} + 2 \paren{1 + \majnormAl}^2 \majnormAl^2 \norm{F}^2 \E\croch{\norm{\varepsilon}^2} + \paren{1 + \majnormAl}^4 \norm{F}^4 } \\ \detailtag
&= 4 \paren{ \majnormAl^4 (n^2 + 2n) \sigma^4 + 2 \paren{1 + \majnormAl}^2 \majnormAl^2 \norm{F}^2 n \sigma^2 + \paren{1 + \majnormAl}^4 \norm{F}^4 } \\
&
}
\leq 4 \paren{ (n+1) \sigma^2 \majnormAl^2 + \paren{1 + \majnormAl}^2 \norm{F}^2 }^2 \enspace .
\end{align*}
Using also Eq.~\eqref{eq.thm.oracle.E.pr.1} and~\eqref{eq.thm.oracle},  Eq.~\eqref{eq.thm.oracle.E} follows, taking $\delta=2\,$. \qed

\section{Proof of the concentration inequalities} \label{A.sec.proof-conc}

\subsection{Proof of Proposition~\ref{pro.Omega.ridge.del34.Gauss}} \label{A.sec.proof.le.Omega.ridge.del34.Gauss}
The proof of Proposition~\ref{pro.Omega.ridge.del34.Gauss} relies on Corollary~\ref{cor.le.conc-linear.Gauss.general}, which is itself a consequence of a general concentration result (Lemma~\ref{le.conc-linear.Gauss.general}). 
Both results, which are proved at the end of the subsection, rely on a rather classical argument: a concentration inequality for the supremum of a Gaussian process, and an upper bound for its expectation in terms of entropy. Then, the proof reduces to bounding the length of a $\mathcal{C}^1$ path inside the unit Euclidean ball, which is done in Eq.~\eqref{eq.proof.le.conc-linear.Gauss.general.1} under the assumptions of Lemma~\ref{le.conc-linear.Gauss.general}. 
\begin{lemma} \label{le.conc-linear.Gauss.general}
Let $n \in \N \backslash \set{0}$, $Z \geq 0\,$. 
Let $\xi$ be a standard Gaussian vector in $\R^n$, and $u : (0,+\infty) \mapsto \R^n$ be some function such that 
\begin{equation} \label{eq.hyp.u.le-conc-lin-Gauss-gal}  \tag{\bf Hp}
\left.
\begin{aligned}
\forall t \in (0,+\infty) \, , \quad \norm{u(t)}^2 = \sum_{j=1}^n u_j(t)^2 \leq 1 
\qquad 
\forall j \in \set{1, \ldots, n} \, , \quad u_j \in \mathcal{C}^1 ((0,+\infty)) \\
\mbox{and either} \quad 
u_j^{\prime} \equiv 0 \mbox{ or } \card\set{ t \in (0,+\infty) \telque u_j^{\prime}(t) = 0 } \leq Z 
\end{aligned}
\right\}
\end{equation}
Then, $u$ admits a continuous extension to $[0,+\infty]$, and for every $x \geq 0$,
\begin{equation} \label{eq.le-conc-lin-Gauss-gal}
\Prob\paren{ \sup_{t \in [0,+\infty]} \absj{ \prodscal{\xi}{u(t)} } \leq \sqrt{2 x} + 12 \sqrt{\ln( 2 + 4 n (Z+1))} + 6 \sqrt{\pi} } \geq 1-  e^{-x} \enspace . 
\end{equation}
\end{lemma}

\begin{corollary} \label{cor.le.conc-linear.Gauss.general}
Let $D_R,D_S \in \N\,$, and for every $j \in \set{1, \ldots, n}\,$, let $R_j \in \R[X]$ be a polynomial of degree at most $D_R\,$, and $S_j \in \R[X]$ be a polynomial of degree at most $D_S\,$ with no positive root. If the polynomials $(R_j)_{1 \leq j \leq n}$ have no common positive root, then 
\[ \forall t \in (0,+\infty) \, , \quad x(t) = \paren{ \frac{R_j(t)}{S_j(t)} }_{1 \leq j \leq n} \] is well-defined and non-zero, so that $u : t \to \norm{x(t)}^{-1} x(t)$ is well-defined on $(0,+\infty)\,$. 
Moreover, $u$ satisfies the assumption \eqref{eq.hyp.u.le-conc-lin-Gauss-gal} with $Z = (3 D_R + (3n-2) D_S - 1)_+\,$, so that Eq.~\eqref{eq.le-conc-lin-Gauss-gal} holds true. 
\end{corollary}

We can now prove Proposition~\ref{pro.Omega.ridge.del34.Gauss}.
\begin{proof}[Proof of Proposition~\ref{pro.Omega.ridge.del34.Gauss}]
Let $P$ be an orthogonal matrix such that Eq.~\eqref{eq.Omega.ridge.key} holds true (that is, $P$ is an orthogonal matrix which diagonalizes $K$), $\xi = \sigma^{-1} P \varepsilon$ and $PF=(f_j)_{1 \leq j \leq n}\,$. 
Then, $\xi$ is a standard Gaussian vector in $\R^n$ and 
for every $\lambda\in(0,+\infty)\,$, 
\begin{gather*} 
\absj{\delta_3(\lambda)} = 2 \absj{\prodscal{\varepsilon}{\Al^{\top} (\Al - \Id_n) F}} 
= 2 \sigma  \absj{\prodscal{\xi}{\Dl (\Dl - \Id_n) PF}} = 2 \sigma \lambda \absj{\prodscal{\xi}{x(\lambda)}} 
\\
\mbox{with} \qquad 
x(\lambda) = \lambda^{-1} \paren{ \frac{\mu_j}{\mu_j+n\lambda} \paren{ \frac{\mu_j}{\mu_j+n\lambda} - 1} f_j }_{1 \leq j \leq n} 
= \paren{  \frac{- n \mu_j f_j}{(\mu_j+n\lambda)^2}  }_{1 \leq j \leq n} 
\enspace . 
\end{gather*}
Therefore, Corollary~\ref{cor.le.conc-linear.Gauss.general} can be applied with $D_R=0$ and $D_S=2\,$, hence $Z=6n-5\,$: for every $y \geq 0$, an event of probability at least $1 - e^{-y}$ exists on which $\forall \lambda \in [ 0,+\infty] \,$,
\begin{align*}
\absj{\delta_3(\lambda)} 
&\leq 2 \sigma \norm{\Al^{\top} (\Al - \Id_n) F} \paren{ \sqrt{2y} + 12 \sqrt{\ln( 2 + 4 n (6n-4))} + 6 \sqrt{\pi} } 
\\
&\leq 2 \sigma \norm{\Al} \norm{ (\Al - \Id_n) F} \paren{ \sqrt{2} \sqrt{y}  + 12 \sqrt{ 2 \ln(n)}  + \paren{ 6 \sqrt{\pi} + 12 \sqrt{\ln( 24)} } }
\detail{
\\ &\leq  2 \sigma \norm{ (\Al - \Id_n) F} \sqrt{ 2 + 288 + 1 } \sqrt{ y + \ln(n) + \paren{ 6 \sqrt{\pi} + 12 \sqrt{\ln( 24)} }^2 } \detailtag
}
\\
&\leq \cteLinRidge \sigma \norm{ (\Al - \Id_n) F} \sqrt{ y + \ln(n) + \betaLRUa }
\end{align*}
by Cauchy-Schwarz inequality. 
So, for every $x \geq 0\,$, taking $y = x -  \ln(n) - \betaLRUa\,$, an event of probability at least 
$1 - \exp(\ln(n) + \betaLRUa) \exp(-x)$ exists on which Eq.~\eqref{eq.conc.del3.ridge-unif.Gauss} holds true for every $\lambda\in[0,+\infty]\,$.

\medskip

Similarly, for every $\lambda\in [0,+\infty]\,$, 
\begin{gather*} 
\absj{\delta_4(\lambda)} = 2 \absj{\prodscal{\varepsilon}{(\Id_n - \Al) F}} 
= 2 \sigma \absj{\prodscal{\xi}{(\Id_n - \Dl) PF}} = 2 \sigma \lambda \absj{\prodscal{\xi}{x(\lambda)}} 
\\
\mbox{with} \qquad 
x(\lambda) = \lambda^{-1} \paren{ \paren{ 1 - \frac{\mu_j}{\mu_j+n\lambda} } f_j }_{1 \leq j \leq n} 
= \paren{  \frac{n f_j}{\mu_j+n\lambda}  }_{1 \leq j \leq n} 
\enspace ,
\end{gather*}
so that Corollary~\ref{cor.le.conc-linear.Gauss.general} can be applied with $D_R=0$ and $D_S=1\,$, hence $Z=3n-3\,$: for every $y \geq 0$, an event of probability at least $1 - e^{-y}$ on which $\forall \lambda \in [0,+\infty]\,$,
\begin{align*}
\absj{\delta_4(\lambda)} 
&\leq 2 \sigma \norm{(\Id_n - \Al) F} \paren{ \sqrt{2y} + 12 \sqrt{\ln( 2+4n(3n-2))} + 6 \sqrt{\pi} } 
\\
&\leq 2 \sigma  \norm{ (\Al - \Id_n) F} \paren{ \sqrt{2} \sqrt{y}  + 12 \sqrt{ 2 \ln(n)}  + \paren{ 6 \sqrt{\pi} + 12 \sqrt{\ln( 12)} } }
\detail{
\\ \detailtag &\leq  2 \sigma \norm{ (\Al - \Id_n) F} \sqrt{ 2 + 288 + 1} \sqrt{ y + \ln(n) + \paren{ 6 \sqrt{\pi} + 12 \sqrt{\ln( 12)} }^2 }
}
\\
&\leq \cteLinRidge \sigma \norm{ (\Al - \Id_n) F} \sqrt{ y + \ln(n) + \betaLRUb }
\enspace .
\end{align*}
So, for every $x \geq 0\,$, taking $y = x -  \ln(n) - \betaLRUb\,$, an event of probability at least 
$1 - \exp(\ln(n) + \betaLRUb) \exp(-x)$ exists on which Eq.~\eqref{eq.conc.del4.ridge-unif.Gauss} holds true for every $\lambda\in [0,+\infty] \,$.

The result follows by taking an union bound\nodetail{.}
\detail{ since 
\begin{equation} \detailtag \exp(\betaLRUa + \ln(1 + \exp(\betaLRUb-\betaLRUa) ) ) 
\leq \exp( \betaLRUa + e^{\betaLRUb-\betaLRUa} ) \leq e^{\betaLRUc} \enspace . \end{equation} }
\end{proof}

\begin{proof}[Proof of Lemma~\ref{le.conc-linear.Gauss.general}]
Let $T = \set{u(t) \telque t \in (0,+\infty)}\,$. 
Since $T$ is a $\mathcal{C}^1$ path, its length $L(T)$ is well-defined (possibly infinite) and satisfies
\begin{align}  
L(T) &= \int_0^{+\infty} \sqrt{ \sum_{j=1}^n \paren{ u_j^{\prime \, 2}(t)}} dt 
\leq \sum_{j=1}^n  \int_0^{+\infty} \absj{  u_j^{\prime}(t) } dt 
\leq 2 n (Z+1) \enspace .
\label{eq.proof.le.conc-linear.Gauss.general.1}
\end{align}
Indeed, for every $j = 1 , \ldots , n\,$, either $u_j^{\prime}\equiv 0$ so that $\int_0^{+\infty} \absj{  u_j^{\prime}(t) } dt=0$, or $u_j^{\prime}$ has at most $Z$ zeros, so $\int_0^{+\infty} \absj{  u_j^{\prime}(t) } dt$ is the sum of the amplitudes of variation of $u_j$ over at most $(Z+1)$ intervals, each term being smaller or equal to 2 since $\absj{u_j(t)} \leq 1$ for every $t>0$ by assumption.

Eq.~\eqref{eq.proof.le.conc-linear.Gauss.general.1} also implies $u_j(t)$ has finite limits when $t \to 0$ and when $t \to +\infty$ for every $j$, so $u(0)$ and $u(+\infty)$ can be defined by continuity and $L(\overline{T}) = L(T)$ with $\overline{T} \egaldef \set{ u(t) \telque t \in [0,+\infty] } \,$.

For any set $S$, let $H(\delta,S) = \ln(N(\delta,S))$ be the metric entropy of $S$ w.r.t. $\norm{\cdot}$. 
Since $\overline{T}$ is a continuous path, 
\begin{equation} \label{eq.maj-entropie-path} 
 \forall \delta > 0 \, , \quad N(\delta,\overline{T} \cup (-\overline{T}))  \leq 2 N(\delta,\overline{T})  \leq 2 \left\lceil \frac{L(\overline{T})}{\delta} \right\rceil \leq 2 + \frac{2 L(\overline{T})}{\delta} \enspace .
\end{equation}
In particular, $\sqrt{H(\cdot, \overline{T} \cup (-\overline{T}))}$ is integrable at 0, so Theorem~3.18 in \cite{Mas:2003:St-Flour} yields
\begin{equation} \label{eq.Esup-proc-Gauss} 
\E\croch{ \sup_{z \in \overline{T}} \absj{ \prodscal{\xi}{z} }} 
= \E\croch{ \sup_{z \in \overline{T} \cup (-\overline{T})} \set{ \prodscal{\xi}{z} }} 
\leq 12 \int_0^{1} \sqrt{H(\delta, \overline{T} \cup (-\overline{T}))} d\delta \enspace .
\end{equation}
Combining Eq.~\eqref{eq.proof.le.conc-linear.Gauss.general.1}, \eqref{eq.maj-entropie-path} and~\eqref{eq.Esup-proc-Gauss}, we get that 
\begin{align*}
\E\croch{ \sup_{z \in \overline{T}} \absj{ \prodscal{\xi}{z} } } 
\detail{
\leq 12 \int_0^{1} \sqrt{\ln\paren{2 + \frac{4 n (Z+1)}{\delta}} } d\delta \detailtag
}
&\leq 12 \sqrt{\ln( 2 + 4 n (Z+1))} + 12 \int_0^{1} \sqrt{\ln\paren{\frac{1}{\delta}} } d\delta 
\\
&\leq 12 \sqrt{\ln( 2 + 4 n (Z+1))} + 6 \sqrt{\pi}
\nodetail{\enspace .}
\end{align*}
\detail{since 
\begin{equation} \detailtag 
\int_0^{1} \sqrt{\ln\paren{\frac{1}{\delta}} } d\delta 
= \int_0^{+\infty} 2 x^2 \exp(-x^2) dx =  \int_0^{+\infty} \exp(-x^2) dx 
= \frac{\sqrt{\pi}}{2} \enspace . \end{equation}   }
Finally, by assumption \eqref{eq.hyp.u.le-conc-lin-Gauss-gal}, $\sup_{z \in \overline{T}} \norm{z} \leq 1\,$, so by Proposition~3.19 in \cite{Mas:2003:St-Flour}, with probability $1 - e^{-x}$, 
\begin{equation} \label{eq.conc.sup-proc-Gauss} 
\sup_{z \in \overline{T}} \absj{ \prodscal{\xi}{z}} \leq \E\croch{ \sup_{z \in \overline{T}} \absj{ \prodscal{\xi}{z}} } + \sqrt{2 x} \enspace ,
\end{equation}
hence the result. 
\end{proof}

\begin{proof}[Proof of Corollary~\ref{cor.le.conc-linear.Gauss.general}]
First, $x(t)$ is well-defined for every $t>0$ since $S_j$ has no positive root for all $j = 1, \ldots,  n\,$. 
Second, $x(t) \neq 0$ for every $t >0$ since the $R_j$ have no common positive root, so that $u$ is well-defined on $(0,+\infty)\,$. 
For every $t>0\,$, let $N(t) \egaldef \norm{x(t)} > 0\,$. 
Each coordinate $x_j$ of $x$ is of class $\mathcal{C}^1$ because it is a well-defined rational fraction, so $N$ also is of class $\mathcal{C}^1$, as well as each coordinate $u_j$ of $u$ and for every $t \in (0,+\infty)$ and $j \in \set{1 , \ldots, n}\,$, 
\begin{align*} 
N^{\prime}(t) &= \frac{\prodscal{x(t)}{x^{\prime}(t)}}{N(t)} 
\\
u_j^{\prime}(t) 
&= \frac{x_j^{\prime}(t) N(t) - x_j(t) N^{\prime}(t)}{N(t)^2} 
= \frac{1}{N(t)^3} \paren{ x_j^{\prime}(t) \paren{N(t)}^2 - \prodscal{x(t)}{x^{\prime}(t)} x_j(t) } 
\detail{
\\ \detailtag &= \frac{1}{N(t)^3} \paren{ x_j^{\prime}(t) \sum_{k=1}^n \paren{ x_k(t)^2 } - x_j(t) \sum_{k=1}^n \paren{ x_k(t) x_k^{\prime}(t) } } 
}
\\ &= \frac{1}{N(t)^3} \croch{ \frac{R_j^{\prime}(t) S_j(t) - R_j(t) S_j^{\prime}(t)}{S_j(t)^2} \sum_{k=1}^n \paren{ \frac{R_k(t)^2}{S_k(t)^2} } - \frac{R_j(t)}{S_j(t)} \sum_{k=1}^n \paren{ \frac{R_k(t) R_k^{\prime}(t) S_k(t) - \paren{R_k(t)}^2 S_k^{\prime}(t)}{S_k(t)^3} } } 
\\ &= \frac{P_j(t)}{ \paren{N(t)}^3 S_j(t) \paren{\prod_{k=1}^n S_k(t)}^3 }
\end{align*}
for some polynomial $P_j \in \R[X]$, either equal to the null polynomial, or of degree smaller or equal to 
\[ 
Z = \max\set{ 0 \, , \, 3 D_R + (3n-2) D_S - 1 } \enspace , \]
which proves \eqref{eq.hyp.u.le-conc-lin-Gauss-gal} holds true.
\end{proof}

\subsection{Proof of Proposition~\ref{A.pro.conc.quad.gal.Gauss.momexp.improved}}
\label{A.sec.proof.pro.conc.quad.gal.Gauss.momexp.improved}
The case where $M$ is a diagonal matrix is Lemma~1 in \cite{Lau_Mas:2000}. Let us prove how the general case can be reduced to the diagonal case. 
Let $B = \frac{1}{2} \paren{ \transpose{M} + M}$ so that $Z = \prodscal{X}{BX}-\tr(B)\,$. 
Since $B$ is symmetric, it can be diagonalized in an orthonormal basis:
\[ \exists P \in O(n) \telque B = P^{\top} D P \quad\mbox{with}\quad D = \diag\paren{d_1, \ldots, d_n} \eqpoint \]
Hence, 
\[ Z = \transpose{X}BX - \tr(B) = \transpose{(PX)} D (PX) - \tr(D) = \sum_{i=1}^n d_i \paren{ \xi_i^2  - 1} \]
where  $\xi = PX \in \R^n$ is a standard Gaussian vector.
By Lemma~1 in \cite{Lau_Mas:2000}, we get that for every $x \geq 0$,
\[ \Proba\paren{Z \geq  2 \sqrt{\tr(\transpose{B}B) x} + 2 \normmat{B}x} \leq e^{-x} \enspace . \]
(It is assumed that $d_i \geq 0$ in \cite{Lau_Mas:2000}, but the proof actually does not use it for proving the above inequality.)
The result follows since 
\begin{align*}
\normmat{B} &= \normmat{\frac{1}{2} \paren{ \transpose{M} + M}} \leq \frac{1}{2} \paren{\normmat{\transpose{M}} + \normmat{M}} = \normmat{M} \\
\mbox{and} \quad \tr(\transpose{B}B) &= \frac{1}{4} \tr\croch{ \paren{ \transpose{M} + M}\paren{ \transpose{M} + M} } \\
\detail{ \detailtag
&= \frac{1}{4} \croch{\tr\paren{ \paren{\transpose{M}}^2} + \tr\paren{M^2} + 2 \tr\paren{\transpose{M} M}} \\
}
&= \frac{1}{2} \croch{\tr\paren{M^2} + \tr\paren{\transpose{M} M}} \eqpoint \qed
\end{align*}

\subsection{Proof of Proposition~\ref{pro.Omega.ridge.del12}} \label{A.sec.proof.le.Omega.ridge.del12}
%
%
Proposition~\ref{pro.hyp.Al} shows that \eqref{hyp.ridge} implies \eqref{hyp.Al} with $\majnormAl=\cAl=1$.
Let us now consider the concentration inequalities \eqref{eq.conc.del1.ridge-unif} and~\eqref{eq.conc.del2.ridge-unif} for $\delun$ and $\deldeux\,$. 
%
%
By Lemma~\ref{le.ridge.monotone}, a sequence $\lambda_1^a > \dots > \lambda_{n-1}^a \in \Lambda$ exists such that for every $j \in \set{1, \dots, n-1}\,$, $\tr(A_{\lambda_j^a}) = j\,$.
Similarly, a sequence $\lambda_1^b > \dots > \lambda_{n-1}^b \in \Lambda$ exists such that for every $j \in \set{1, \ldots, n-1}\,$, $\tr(A_{\lambda_j^b}^{\top} A_{\lambda_j^b}) = j\,$.
Let 
\[ \Lambda_1 = \set{ 0, \lambda_1^a , \dots, \lambda_{n-1}^a , \lambda_1^b , \dots , \lambda_{n-1} ^b } \enspace , \] 
so that $\card(\Lambda_1) \leq 2 n$ and for every  $\lamm\,$, some $\lambda_+, \lambda_- \in \Lambda_1 $ exist such that 
\begin{align}  \label{eq.le.Omega.ridge.pr1}
\lambda_- &\leq \lambda \leq \lambda_+ \\
\label{eq.le.Omega.ridge.pr2}
\tr(A_{\lambda_-}) - 1 \leq \tr(A_{\lambda_+}) &\leq \tr(\Al) \leq \tr(A_{\lambda_-}) \leq \tr(A_{\lambda_+}) + 1 \\
\label{eq.le.Omega.ridge.pr3}
\tr(A_{\lambda_-}^{\top} A_{\lambda_-}) - 1 \leq \tr(A_{\lambda_+}^{\top} A_{\lambda_+}) &\leq \tr(\Al^{\top} \Al) \leq \tr(A_{\lambda_-}^{\top} A_{\lambda_-}) \leq \tr(A_{\lambda_+}^{\top} A_{\lambda_+}) + 1 \enspace . 
\end{align}  


Using the notation introduced in Section~\ref{A.sec.ridge}, let $\xi = P \varepsilon \sim \mathcal{N}(0,\sigma^2 \Id_n)\,$, so that 
\begin{align*}
\norm{\Al \varepsilon}^2 = \sum_{j=1}^n \paren{ \xi_j^2 \frac{\mu_j^2}{(\mu_j + n \lambda)^2}} 
\quad \mbox{and} \quad 
\prodscal{\varepsilon}{\Al \varepsilon} = \sum_{j=1}^n \paren{ \xi_j^2 \frac{\mu_j}{\mu_j + n \lambda}} 
\end{align*}
both are non-increasing functions of $\lambda$ since $\mu_j \geq 0\,$.


Let us now assume the event $\Omega_x \paren{ \Lambda_1, \Comega }$ is realized, so that Eq.~\eqref{A.eq.conc.del1} and~\eqref{A.eq.conc.del2} hold for any $\lamm_1$. 
For any $\lamm\,$, let $\lambda_-, \lambda_+\in\Lambda_1$ such that Eq.~\eqref{eq.le.Omega.ridge.pr1}, \eqref{eq.le.Omega.ridge.pr2} and~\eqref{eq.le.Omega.ridge.pr3} hold true. 
Then, since $\lambda \mapsto \norm{\Al \varepsilon}^2$ is non-increasing, for every $\theta_1 \in (0,1]\,$,
\begin{align*}
\delun(\lambda) 
\detail{ \detailtag
&= \norm{\Al \varepsilon}^2 - \sigma^2 \tr(\Al^{\top}\Al) \\
}
&\geq \norm{A_{\lambda_+} \varepsilon}^2 - \sigma^2 \tr(\Al^{\top}\Al) 
\qquad \mbox{by Eq.~\eqref{eq.le.Omega.ridge.pr1}} \\
&\geq \norm{A_{\lambda_+} \varepsilon}^2 - \sigma^2 \tr(A_{\lambda_+}^{\top}A_{\lambda_+}) - \sigma^2 
\qquad \mbox{by Eq.~\eqref{eq.le.Omega.ridge.pr3}} \\
\detail{ \detailtag
&= \delun(\lambda_+) - \sigma^2  \\
}
&\geq - \theta_1  \sigma^2 \tr(A_{\lambda_+}^{\top}A_{\lambda_+}) - \paren{ \ComegaA + \ComegaB \theta_1^{-1}} x \sigma^2 - \sigma^2 
\qquad \mbox{by Eq.~\eqref{A.eq.conc.del1}} \\
&\geq - \theta_1  \sigma^2 \tr(\Al^{\top}\Al) - \paren{ \ComegaA + \ComegaB \theta_1^{-1}} x \sigma^2 - \sigma^2 
\qquad \mbox{by Eq.~\eqref{eq.le.Omega.ridge.pr3}}
\end{align*}
and 
\begin{align*}
\delun(\lambda) 
\detail{ \detailtag
&= \norm{\Al \varepsilon}^2 - \sigma^2 \tr(\Al^{\top}\Al) 
\\ 
}
&\leq \norm{A_{\lambda_-} \varepsilon}^2 - \sigma^2 \tr(\Al^{\top}\Al) 
\qquad \mbox{by Eq.~\eqref{eq.le.Omega.ridge.pr1}} 
\\ &\leq \norm{A_{\lambda_-} \varepsilon}^2 - \sigma^2 \tr(A_{\lambda_-}^{\top}A_{\lambda_-}) + \sigma^2 \qquad \mbox{by Eq.~\eqref{eq.le.Omega.ridge.pr3}} 
\detail{
\\ \detailtag &= \delun(\lambda_-) + \sigma^2 
}
\\ &\leq \theta_1  \sigma^2 \tr(A_{\lambda_-}^{\top}A_{\lambda_-}) + \paren{ \ComegaA + \ComegaB \theta_1^{-1}} x \sigma^2 + \sigma^2
\qquad \mbox{by Eq.~\eqref{A.eq.conc.del1}} 
\\ &\leq \theta_1  \sigma^2 \tr(\Al^{\top}\Al) + (1+\theta_1) \sigma^2 + \paren{ \ComegaA + \ComegaB \theta_1^{-1}} x \sigma^2 
\qquad \mbox{by Eq.~\eqref{eq.le.Omega.ridge.pr3}.}
\end{align*}
So, for every $\theta_1\in (0,1]\,$, Eq.~\eqref{eq.conc.del1.ridge-unif} holds true.

Similarly, since $\lambda \mapsto \prodscal{\varepsilon}{\Al \varepsilon}$ is non-increasing, 
\nodetail{
by Eq.~\eqref{eq.le.Omega.ridge.pr1}, \eqref{eq.le.Omega.ridge.pr2}, \eqref{eq.le.Omega.ridge.pr3} and~\eqref{A.eq.conc.del2}, 
}
for every $\theta_2 \in (0,1] \,$,
\begin{align*}
\deldeux(\lambda) 
\superdetail{ \detailtag
&= \prodscal{\varepsilon}{\Al \varepsilon} - \sigma^2 \tr(\Al) \\
\detailtag &\geq \prodscal{\varepsilon}{A_{\lambda_+} \varepsilon} - \sigma^2 \tr(\Al) 
\qquad \mbox{by Eq.~\eqref{eq.le.Omega.ridge.pr1}} \\
\detailtag &\geq \prodscal{\varepsilon}{A_{\lambda_+} \varepsilon} - \sigma^2 \tr(A_{\lambda_+}) - \sigma^2 
\qquad \mbox{by Eq.~\eqref{eq.le.Omega.ridge.pr2}} \\
\detailtag &= \deldeux(\lambda_+) - \sigma^2 \\
\detailtag &\geq - \theta_2  \sigma^2 \tr(A_{\lambda_+}^{\top} A_{\lambda_+}) - \paren{ \ComegaC + \ComegaD \theta_2^{-1}} x \sigma^2 - \sigma^2 
\qquad \mbox{by Eq.~\eqref{A.eq.conc.del2}} \\
}
&\geq - \theta_2  \sigma^2 \tr(\Al^{\top} \Al) - \paren{ \ComegaC + \ComegaD \theta_2^{-1}} x \sigma^2 - \sigma^2 
\superdetail{ \qquad \mbox{by Eq.~\eqref{eq.le.Omega.ridge.pr3}} }
\\
\mbox{and} \quad  
\deldeux(\lambda) 
\superdetail{ \detailtag
&= \prodscal{\varepsilon}{\Al \varepsilon} - \sigma^2 \tr(\Al)
\\ \detailtag &\leq \prodscal{\varepsilon}{A_{\lambda_-} \varepsilon} - \sigma^2 \tr(\Al) 
\qquad \mbox{by Eq.~\eqref{eq.le.Omega.ridge.pr1}} 
\\ \detailtag &\leq \prodscal{\varepsilon}{A_{\lambda_-} \varepsilon} - \sigma^2 \tr(A_{\lambda_-}) + \sigma^2  \qquad \mbox{by Eq.~\eqref{eq.le.Omega.ridge.pr2}} 
\\ \detailtag &= \deldeux(\lambda_-) + \sigma^2 
\\ \detailtag &\leq \theta_2  \sigma^2 \tr(A_{\lambda_-}^{\top}A_{\lambda_-}) + \paren{ \ComegaC + \ComegaD \theta_2^{-1}} x \sigma^2 + \sigma^2 
\qquad \mbox{by Eq.~\eqref{A.eq.conc.del2}} 
\\ 
}
&\leq \theta_2  \sigma^2 \tr(\Al^{\top}\Al) + (1+\theta_2) \sigma^2 + \paren{ \ComegaC + \ComegaD \theta_2^{-1}} x \sigma^2 
\superdetail{ \qquad \mbox{by Eq.~\eqref{eq.le.Omega.ridge.pr3}.} } \nodetail{\enspace . }
\end{align*}
So, for every $\theta_2 \in (0,1]\,$, Eq.~\eqref{eq.conc.del2.ridge-unif} holds true.
%
\qed

\subsection{Proof of Lemma~\ref{le.Omega.finite}} \label{A.sec.proof.le.Omega.finite}
%
%
Since $\Lambda$ is finite by assumption \eqref{hyp.Lam-dis}, we use a union bound over $\lamm$. 
For each $\lamm$, using assumption \eqref{hyp.eps-Gauss-hom}, $\sigma^{-1} \varepsilon$ is a standard Gaussian vector.
\begin{itemize}
\item By Proposition~\ref{A.pro.conc.quad.gal.Gauss.momexp.improved} in Section~\ref{A.sec.conc.quad-Gauss} with $M= \pm \sigma^2 \Al$ and $M= \pm \sigma^2 \Al^{\top} \Al\,$, we deduce that for every $x \geq 0$, 
\begin{align*}
\Prob\paren{ \forall \theta>0, \,\, \absj{ \delun(\lambda) } \leq \theta \sigma^2 \tr(\Al^\top \Al) + \paren{2 \normmat{\Al} + \theta^{-1} } x \sigma^2 } &\geq 1 - 2e^{-x} 
\\
\Prob\paren{ \forall \theta>0, \,\, \absj{ \deldeux(\lambda) } \leq \theta \sigma^2  \tr(\Al^\top \Al) + \paren{2 + \theta^{-1} } \normmat{\Al}^2 x \sigma^2  } &\geq 1 - 2 e^{-x} 
\enspace ,
\end{align*}
where we used Eq.~\eqref{eq.maj-theta}, $\normmat{\Al^{\top} \Al} \leq \normmat{\Al}^2\,$, and that $\tr((\Al^\top \Al)^2) \leq \normmat{\Al}^2 \tr(\Al^\top \Al)\,$.
Eq.~\eqref{A.eq.conc.del1} and~\eqref{A.eq.conc.del2} follow, using that $\normmat{\Al} \leq \majnormAl$ by assumption \eqref{hyp.Al}.
\item Since $ \delta_3(\lambda) = \prodscal{ \sigma^{-1} \varepsilon} { 2 \sigma \Al^\top (\Id_n - \Al) F }$  and $\delta_4(\lambda) = \prodscal{ \sigma^{-1} \varepsilon} { 2 \sigma  (\Id_n - \Al) F}  \,$, Proposition~\ref{A.pro.conc.linear.gal.Gauss} in Section~\ref{A.sec.conc.lin-Gauss} shows, for every $x\geq 0\,$,
\begin{gather*}
\Prob\paren{ \forall \theta>0, \,\, \absj{\delta_3(\lambda)} \leq \theta \norm{(\Id_n - \Al) F}_2^2 +  \frac{2 \normmat{\Al}^2 x \sigma^2}{ \theta} } \geq 1 - e^{-x} 
\\
\mbox{and} \quad 
\Prob\paren{ \forall \theta>0, \,\, \absj{\delta_4(\lambda)} \leq \theta \norm{(\Id_n - \Al) F}_2^2 + \frac{2 x \sigma^2}{\theta} } \geq 1 - e^{-x} \enspace ,
\end{gather*}
using Eq.~\eqref{eq.maj-theta}.
Eq.~\eqref{A.eq.conc.del3} and~\eqref{A.eq.conc.del4} follow, using that $\normmat{\Al} \leq \majnormAl$ by assumption~\eqref{hyp.Al}.
\end{itemize}
\qed

\subsection{Proof of Lemma~\ref{le.Omega.cont-ridge}} \label{A.sec.proof.le.Omega.cont-ridge}
%
%
Since \eqref{hyp.ridge} holds true, we can apply Proposition~\ref{pro.Omega.ridge.del12}: some $\Lambda_1$ exists such that $\card(\Lambda_1) \leq 2 n$ and for every $x \geq 2/3$, 
\begin{equation} \label{eq.proof.le.Omega.cont-ridge.1} \Omega_{x-2/3} \paren{\Lambda_1 , \Comega } \subset \Omega_{x} \paren{\Lambda , \Comega} \enspace .  \end{equation}
Since $\Lambda_1$ is finite, and assumptions \eqref{hyp.eps-Gauss-hom} and \eqref{hyp.Al} hold true with $\majnormAl=\cAl=1$, we get by the first part of proof of Lemma~\ref{le.Omega.finite} that for every $x \geq 2/3$, 
\begin{equation} \label{eq.proof.le.Omega.cont-ridge.2} \Proba \paren{ \Omega_{x-2/3} \paren{\Lambda_1 , \paren{2, 1, 2, 1, +\infty, +\infty}} } \geq 1 - 6 e^{-x + 2/3} \enspace . \end{equation}
Now, by Proposition~\ref{pro.Omega.ridge.del34.Gauss} and Eq.~\eqref{eq.maj-theta}, 
\begin{equation} \label{eq.proof.le.Omega.cont-ridge.3} \Proba \paren{ \Omega_x \paren{\Lambda_1 , \paren{+\infty, +\infty, +\infty, +\infty, \cteC , \cteC }} } \geq 1 - e^{-x + \betaLRUc + \ln(n)} \enspace . \end{equation}
Combining Eq.~\eqref{eq.proof.le.Omega.cont-ridge.1}, \eqref{eq.proof.le.Omega.cont-ridge.2}, \eqref{eq.proof.le.Omega.cont-ridge.3} and an union bound, we get for every $x \geq 2/3$, 
\[ 
\Proba \paren{ \Omega_x \paren{\Lambda_1 , \paren{2, 1, 2, 1, \cteC, \cteC}} } \geq 1 - e^{-x + \ln(n)} \paren{ e^{\betaLRUc} + 6 e^{2/3} } \geq 1 - e^{-x + \ln(n) + \betaLRUd} \enspace ,
\]
a bound also valid when $0 \leq x < 2/3$ (since it is negative).  \qed

\section*{Acknowledgments}
We would like to thank Matthieu Solnon for helping us to improve an earlier version of the paper.
We also acknowledge the support of the French Agence Nationale de la Recherche (ANR) under reference ANR-09-JCJC-0027-01 ({\sc Detect} project) as well as the European Research Council (SIERRA starting grant 239993).

\bibliographystyle{plain}
\bibliography{minikernel}

\clearpage
\section{Supplementary material} \label{sec.supmat}

\subsection{Technical lemmas}

\begin{lemma} \label{le.trM^2<trMTM}
For any $M  \in \M_n(\R)$, 
\begin{equation}  \label{eq.le.trM^2<trMTM}
 \tr(M^2) \leq \tr(M^{\top} M)
\end{equation}
\end{lemma}
\begin{proof}[Proof of Lemma~\ref{le.trM^2<trMTM}]
\begin{align*}
 \tr(M^2) = \sum_{i=1}^n \sum_{j=1}^n M_{i,j} M_{j,i} 
\leq \sum_{i=1}^n \sum_{j=1}^n \frac{M_{i,j}^2 +  M_{j,i}^2}{2} 
=  \sum_{i=1}^n \sum_{j=1}^n M_{i,j}^2 
= \tr(M^{\top} M)
 \end{align*}
\end{proof}

\subsection{About the concentration of quadratic forms of Gaussian vectors} \label{sec.supmat.var-quad-form}

\begin{remark}
By Lemma~\ref{le.var.quad} below \textup{(}with $m=2\,$\textup{)}, when $\xi_1, \ldots, \xi_n$ are i.i.d. standard Gaussian variables, 
\begin{equation*} 
\var(\prodscal{\xi}{M \xi}) = \tr(M^2) + \tr(M^{\top} M) \leq 2 \tr(M^{\top} M) \enspace . 
\end{equation*}
Therefore, the deviation term $\sqrt{2 x (\tr(M^2) + \tr(M^{\top}M))}$ cannot be improved in Eq.~\eqref{A.eq.pro.conc.quad.gal.Gauss.momexp.maj.precis.improved}.
\end{remark}

\begin{lemma} \label{le.var.quad}
Let $\xi_1, \ldots, \xi_n$ be independent random variables such that for every $i \in \set{1, \ldots, n}\,$, $\E\croch{\xi_i}=0\,$, $\E\croch{\xi_i^2}=1$ and $\var(\xi_i^2) = \E\croch{\xi_i^4}-1 = m\,$.
Let $M \in \mathcal{M}_n(\R)\,$. 
Then, 
\begin{equation} \label{eq.var-xiMxi}
\var\paren{\prodscal{\xi}{M \xi}} = (m-2) \sum_{i=1}^n M_{i,i}^2  + \tr(M^2) + \tr(M^{\top} M) 
\leq \paren{ 2 + (m-2)_+ } \tr(M^{\top} M) \enspace . 
\end{equation}
\end{lemma}

\begin{proof}[Proof of Lemma~\ref{le.var.quad}]
On the one hand, 
\[ \prodscal{\xi}{M \xi} = \sum_{1 \leq i,j \leq n} M_{i,j} \xi_i \xi_j \] 
so that $ \E\croch{ \prodscal{\xi}{M \xi} } = \tr(M) \,$. 
On the other hand, 
\begin{align*}
\E\croch{ \prodscal{\xi}{M \xi}^2 } 
\detail{
&= \sum_{1 \leq i,j,k,\ell \leq n} M_{i,j} M_{k,\ell} \E\croch{ \xi_i \xi_j \xi_k \xi_{\ell} } \\
}
&= \sum_{i=1}^n M_{i,i}^2 \E\croch{\xi_i^4} 
	+ \sum_{1 \leq i \neq j \leq n} \paren{ M_{i,j} M_{j,i} + M_{i,j}^2 + M_{i,i} M_{j,j}} \\
&= (m-2) \sum_{i=1}^n M_{i,i}^2 + \tr\paren{M^2} + \tr\paren{M^{\top} M} + \paren{\tr(M)}^2 
\end{align*}
so that Eq.~\eqref{eq.var-xiMxi} holds true (using Lemma~\ref{le.trM^2<trMTM}).
\end{proof}

\end{document}